\documentclass[12pt,oneside]{amsart}
\usepackage{psfig}
\usepackage{amsxtra}
\usepackage{amsfonts}
\usepackage{amssymb}
\makeatletter\def\section{\@startsection{section}{1}%
  \z@{.7\linespacing\@plus\linespacing}{.5\linespacing}%
  {\normalfont\bfseries\indent}}
\makeatother
\newtheorem{theorem}{\mdseries\indent THEOREM}[section]
\newtheorem{corollary}[theorem]{\mdseries\indent COROLLARY}
\newtheorem{lemma}[theorem]{\mdseries\indent LEMMA}
\newtheorem{proposition}[theorem]{\mdseries\indent PROPOSITION}
\theoremstyle{definition}
\newtheorem{definition}[theorem]{\mdseries\indent DEFINITION}
\newtheorem{remark}[theorem]{\mdseries\indent REMARK}

\newtheorem{example}[theorem]{\mdseries\indent EXAMPLE}
\theoremstyle{remark}
\newtheorem{claim}{\indent Claim}
\newtheorem*{acknowledgements}{\indent Acknowledgements}
\renewcommand{\theenumi}{\roman{enumi}}
\renewcommand{\labelenumi}{(\theenumi)}

\numberwithin{equation}{section}
\newlength{\qedskip}

\newlength{\betweenpartsskip}
\newlength{\afterabstractskip}
\def\openone
{\mathchoice
{\hbox{\upshape \small1\kern-3.3pt\normalsize1}}
{\hbox{\upshape \small1\kern-3.3pt\normalsize1}}
{\hbox{\upshape \tiny1\kern-2.3pt\SMALL1}}
{\hbox{\upshape \Tiny1\kern-2pt\tiny1}}}
\makeatletter
\newbox\ipbox
\newcommand{\ip}[2]{\left\langle #1\mathrel{\mathchoice
{\setbox\ipbox=\hbox{$\displaystyle \left\langle\mathstrut #1#2\right\rangle$}
\vrule height\ht\ipbox width0.25pt depth\dp\ipbox}
{\setbox\ipbox=\hbox{$\textstyle \left\langle\mathstrut #1#2\right\rangle$}
\vrule height\ht\ipbox width0.25pt depth\dp\ipbox}
{\setbox\ipbox=\hbox{$\scriptstyle \left\langle\mathstrut #1#2\right\rangle$}
\vrule height\ht\ipbox width0.25pt depth\dp\ipbox}
{\setbox\ipbox=\hbox{$\scriptscriptstyle \left\langle\mathstrut #1#2\right\rangle$}
\vrule height\ht\ipbox width0.25pt depth\dp\ipbox}
} #2\right\rangle}
\newcommand{\diracb}[1]{\left\langle #1\mathrel{\mathchoice
{\setbox\ipbox=\hbox{$\displaystyle \left\langle\mathstrut #1\right.$}
\vrule height\ht\ipbox width0.25pt depth\dp\ipbox}
{\setbox\ipbox=\hbox{$\textstyle \left\langle\mathstrut #1\right.$}
\vrule height\ht\ipbox width0.25pt depth\dp\ipbox}
{\setbox\ipbox=\hbox{$\scriptstyle \left\langle\mathstrut #1\right.$}
\vrule height\ht\ipbox width0.25pt depth\dp\ipbox}
{\setbox\ipbox=\hbox{$\scriptscriptstyle \left\langle\mathstrut #1\right.$}
\vrule height\ht\ipbox width0.25pt depth\dp\ipbox}
}\right. }
\newcommand{\dirack}[1]{\left. \mathrel{\mathchoice
{\setbox\ipbox=\hbox{$\displaystyle \left.\mathstrut #1\right\rangle$}
\vrule height\ht\ipbox width0.25pt depth\dp\ipbox}
{\setbox\ipbox=\hbox{$\textstyle \left.\mathstrut #1\right\rangle$}
\vrule height\ht\ipbox width0.25pt depth\dp\ipbox}
{\setbox\ipbox=\hbox{$\scriptstyle \left.\mathstrut #1\right\rangle$}
\vrule height\ht\ipbox width0.25pt depth\dp\ipbox}
{\setbox\ipbox=\hbox{$\scriptscriptstyle \left.\mathstrut #1\right\rangle$}
\vrule height\ht\ipbox width0.25pt depth\dp\ipbox}
} #1\right\rangle}
\newbox\froglongup
\long\def\hooklonguparrow{\setbox\froglongup=\hbox to 0pt{\hss \raisebox{4pt}{\makebox[0pt]{\hss $\displaystyle \uparrow $\hss }}\raisebox{-2pt}{\makebox[0pt]{\hss $\displaystyle |$\hss }}\hss }\raisebox{-5pt}{\raisebox{-2pt}{\setlength{\unitlength}{0.4pt}\begin{picture}(10,5)
\put(5,5){\oval(10,10)[b]}
\end{picture}}\raise\dp\froglongup\box\froglongup}}

\begin{document}
\thanks{Work supported in part by the U.S. National Science Foundation.}
\subjclass{Primary 46L60, 47D25, 42A16, 43A65; Secondary 46L45, 42A65, 41A15}
\keywords{Wavelet, cascade algorithm, refinement operator, representation, orthogonal
expansion, quadrature mirror filter, isometry in Hilbert space}

\vspace*{3cm}
{\normalsize\bf  
A GEOMETRIC APPROACH TO THE CASCADE 

APPROXIMATION OPERATOR FOR WAVELETS\footnote{%
Work supported in part by the U.S. National Science Foundation.%
}}%
\vspace{\betweenpartsskip}

{\normalsize 
Palle~E.~T.~ Jorgensen
}
\vspace{\betweenpartsskip}

\thispagestyle{empty}
\baselineskip=13pt
\begin{quote}
This paper is devoted to an approximation problem for operators in
Hilbert
space, that appears when one tries to study geometrically the
\emph{cascade algorithm} in wavelet theory.
Let $\mathcal{H}$ be a Hilbert space, and let $\pi$ be a representation of
$L^{\infty}\left(  \mathbb{T}\right)  $ on $\mathcal{H}$. Let $R$ be a
positive operator in $L^{\infty}\left(  \mathbb{T}\right)  $ such that
$R\left(  \openone\right)  =\openone$, where
$\openone$ denotes the constant function $1$. We study operators
$M$ on $\mathcal{H}$ (bounded, but non-contractive) such that
\newline
\begin{minipage}{\textwidth}
\hskip-24pt\begin{minipage}{\textwidth}
\begin{equation*}
\pi\left( f\right) M
=M\pi\left( f\left( z\sp{2}\right) \right) \text{\quad and\quad}
M\sp{\ast}\pi\left( f\right) M
=\pi\left( R\sp{\ast}f\right) , \quad f\in L\sp{\infty}\left( \mathbb
{T}\right) ,
\end{equation*}
\end{minipage}\hskip-24pt\vspace*{\abovedisplayskip}
\end{minipage}
where the $\ast$ refers to Hilbert space adjoint. We give a complete
orthogonal expansion of $\mathcal{H}$ which reduces $\pi$ such that $M$ acts
as a shift on one part, and the residual part is $\mathcal{H}^{\left(
\infty\right)  }=\bigcap_{n}\left[  M^{n}\mathcal{H}\right]  $, where $\left[
M^{n}\mathcal{H}\right]  $ is the closure of the range of $M^{n}$. The shift
part is present, we show, if and only if $\ker\left(  M^{\ast}\right)
\neq\left\{  0\right\}  $. We apply the operator-theoretic results to the
refinement operator (or cascade algorithm) from wavelet theory. Using the
representation $\pi$, we show that, for this wavelet operator $M$, the
components in the decomposition are unitarily, and canonically, equivalent to
spaces $L^{2}\left(  E_{n}\right)  \subset L^{2}\left(  \mathbb{R}\right)  $,
where $E_{n}\subset\mathbb{R}$, $n=0,1,2,\dots,\infty$, are measurable subsets
which form a tiling of $\mathbb{R}$; i.e., the union is $\mathbb{R}$ up to zero
measure, and pairwise intersections of different $E_{n}$'s have measure zero.
We prove two results on the convergence of the cascade algorithm, and identify
singular vectors for the starting point of the algorithm.
\end{quote}
\vspace{13pt}

\noindent
\begin{minipage}{0.5\textwidth}
\renewcommand{\theenumi}{\arabic{enumi}}
\renewcommand{\labelenumi}{\theenumi.}
\raggedright
\textsc{Contents:}
\begin{enumerate}
\item  Introduction

\item  Some Main Results in the Paper

\item  Background, summary, and operator relations

\item  The Zak transform

\item Proof of Theorem \textup {3.1\hbox {}}

\item Sub-isometries

\item  Singular cascade approximations

\item  Singular vectors

\item  Approximation results
\end{enumerate}
\end{minipage}%
\begin{minipage}{0.5\textwidth}
\renewcommand{\theenumi}{\arabic{enumi}}
\renewcommand{\labelenumi}{\theenumi.}
\raggedright
\textsc{List of Tables:}
\begin{enumerate}
\item  Summary of Theorem 4.7\hbox {}: Embedding of the isometric model into $L^{2}\left ( \mathbb {R}\right ) $

\item  Approximation properties of $m_{0}^{\left ( n\right ) }$ and $D_{n}^{{}}\left ( z\right ) $

\item  Operator correspondence between $L^{2}\left ( \mathbb {R}\right ) $ and $\mathcal {H}_{Z}$ (Lemma 5.1\hbox {})
\end{enumerate}
\textsc{List of Figures:}
\begin{enumerate}
\item  $\sigma ^{-2}\left ( N\left ( m_{0}\right ) \right ) $

\item  An ergodic map on $\mathbb {T}$
\end{enumerate}
\end{minipage}

\begin{center}
\begin{tabular}{rl}
\multicolumn{2}{c}{\textsc{Terminology used in the paper:}} \\
$\mathbb{T}$: & the one-torus \\
$\mu $: & Haar measure on the torus $\mathbb{T}$ \\
$Z$: & the Zak transform \\
$\tilde{X}=ZXZ^{-1}$: & transformation of operators \\
$\mathcal{H}$: & a given Hilbert space \\
$\pi $: & a representation of $L^{\infty }\left( \mathbb{T}\right) $ on $%
\mathcal{H}$ \\ 
$R$: & the Ruelle operator on $L^{\infty }\left( \mathbb{T}\right) $ \\ 
$M$: & an operator on $\mathcal{H}$ \\ 
$R^{*},M^{*}$: & adjoint operators\\ 
\end{tabular}%
\end{center}
\vspace{\betweenpartsskip}

\baselineskip=16pt

\section[Introduction]{\label{Intr}INTRODUCTION}

In wavelet theory, one is given a \emph{subband filter},
i.e., a function $m_0$ on the unit circle, satisfying
(\ref{IntrAxiom(1)})--(\ref{IntrAxiom(3)}) from
below and one wants to construct a scaling function $\varphi$
(relative to
$m_0$), i.e., a nonzero function $\varphi$ on $\mathbb{R}$ satisfying the scaling
relation
$\varphi=M\varphi$ (relation (\ref{eqIntr.1}) in the paper).
Here $M=M_{m_0}$ is the so-called
cascade operator defined by
\[
\left( Mh\right) \left( x\right)  = \sqrt{2} \sum_{n\in\mathbb Z} a_n h\left( 2x-n\right) ,
\]
where $a_n$ are Fourier coefficients of the function $m_0$.
The scaling function $\varphi$ is important because its shifts generate
(under
some analytic conditions) the so-called \emph{multiresolution subspace}
$V=V\left( \varphi\right) $, which is used to construct the wavelets.

There are several ways of constructing a scaling function $\varphi$. The
cascade algorithm is one of the possibilities. In this algorithm one
picks
some simple function $h$ and considers its iterations  $M^{n}h$. Clearly,
if
the iterations $M^{n}h$ converge to a nonzero function $\varphi$, the function
$\varphi$
satisfies $\varphi=M\varphi$, so the algorithm gives us a scaling function.
We study the problem of convergence of this algorithm in
the
setting of an abstract Hilbert space. The cascade operator $M$ has some
very
special structure. The Ruelle operator $R$, which appears naturally
in this type of problem, gives us a way to describe this structure.
Namely,
the operator $M$ is
what we shall call a
sub-isometry; see Definition \ref{DefR-is.1}.
In fact the concept depends on a pair $\left( R,\pi\right)$
where $R$ is a Ruelle operator and $\pi$ is a representation.
It turns out that sub-isometries admit an analogue of the Kolmogorov--Wold
decomposition for usual isometries. Using this decomposition, we
obtain results about convergence of the cascade algorithm in an
abstract Hilbert space setting, and then, in the last section, give
some applications  for wavelet theory.

This paper is technically concerned with approximation problems for operators
in Hilbert space, but our initial motivation is the \emph{refinement operator}
(alias the \emph{cascade approximation operator\/}) from wavelet theory. Our
starting point is a given function $m_{0}\in L^{\infty}\left(  \mathbb{T}%
\right)  $, which is assumed to satisfy the following three axioms:
\pagebreak
\begin{enumerate}
\item \label{IntrAxiom(1)}$m_{0}$ is continuous in an open neighborhood of
$z=1$ in $\mathbb{T}$,

\item \label{IntrAxiom(2)}$m_{0}=\sqrt{2}$ at $z=1$,

\item \label{IntrAxiom(3)}$\left|  m_{0}\left(  z\right)  \right|
^{2}+\left|  m_{0}\left(  -z\right)  \right|  ^{2}=2$.
\end{enumerate}
In applications (see \cite{Dau92}), $m_{0}$ is called a
\emph{subband filter,} and (\ref{IntrAxiom(2)}) is called
the \emph{low pass filter property.}
That is because of the substitution $z=e^{-i\omega}$, $\omega$ frequency; and
(\ref{IntrAxiom(2)}) implies that all signals pass at $\omega=0$, while
(\ref{IntrAxiom(2)}) and (\ref{IntrAxiom(1)}) imply that ``almost all'' pass
for small $\omega$, i.e., low frequencies. Finally, axiom (\ref{IntrAxiom(3)})
accounts for the name \emph{quadrature filter,} or ``quadrature mirror
filter''. They are used in wavelet theory for generating
\emph{multiresolutions.} Suppose $m_{0}\left(  z\right)  =\sum_{n\in
\mathbb{Z}}a_{n}z^{n}$ is the Fourier expansion of $m_{0}$. A function
$\varphi\in L^{2}\left(  \mathbb{R}\right)  $ is called a \emph{scaling
function} (relative to $m_{0}$) if%
\begin{equation}
\varphi\left(  x\right)  =\sqrt{2}\sum_{n\in\mathbb{Z}}a_{n}\varphi\left(
2x-n\right)  , \label{eqIntr.1}%
\end{equation}
and we refer to the operator%
\begin{equation}
\left(  Mh\right)  \left(  x\right)  =\sqrt{2}\sum_{n\in\mathbb{Z}}%
a_{n}h\left(  2x-n\right)  \label{eqIntr.2}%
\end{equation}
as the \emph{cascade operator.} The approximation problem is that of finding
conditions on the three functions $m_{0}$, $\varphi$, and $h$ ($h\in
L^{2}\left(  \mathbb{R}\right)  $), such that%
\begin{equation}
\lim_{n\rightarrow\infty}\left\|  \varphi-M^{n}h\right\|  _{L^{2}\left(
\mathbb{R}\right)  }=0. \label{eqIntr.3}%
\end{equation}
This is called the \emph{cascade approximation,}
or the \emph{cascade problem.}
Let $R$ be the Ruelle operator (Section \ref{Back}) corresponding to $m_{0}$.
The basic connection between the two operators $R$ and $M$ will be studied in
a geometric Hilbert-space setting (Section \ref{Zak}), where we give a
quadratic form $h\mapsto p_{2}\left(  h\right)  $ on $L^{2}\left(
\mathbb{R}\right)  $ with values in $L^{1}\left(  \mathbb{T}\right)  $. We
then introduce the concept of \emph{sub-isometries.} (The cascade operators
$M$ are examples.) We show that $R\left(  p_{2}\left(  h\right)  \right)
=p_{2}\left(  Mh\right)  $, so if $M\varphi=\varphi$, then $p_{2}\left(
\varphi\right)  $ is an eigenfunction for $R$. In addition to $L^{2}%
$-convergence questions, there are various classes of pointwise convergence
issues. It turns out that, when $m_{0}$ is given, satisfying
(\ref{IntrAxiom(1)})--(\ref{IntrAxiom(3)}), and $\varphi$, $h$ satisfy natural
criteria, then there is an unexpected obstruction to the approximation
(\ref{eqIntr.3}). If $\ker\left(  M^{\ast}\right)  \neq\left\{  0\right\}  $
in $L^{2}\left(  \mathbb{R}\right)  $, then the cascade approximation
(\ref{eqIntr.3}) is exceedingly ``bad''. We show (Proposition \ref{ProSiV.2})
that $\ker\left(  M^{\ast}\right)  =\left\{  0\right\}  $ \emph{if and only if
}$m_{0}$\emph{ does not vanish on a set of positive measure }$\subset
\mathbb{T}$. So if, for example, $m_{0}$ is a polynomial, then $\ker\left(
M^{\ast}\right)  =\left\{  0\right\}  $.

\section[Some Main Results in the Paper]{\label{Some}SOME MAIN RESULTS IN THE PAPER}

When $\ker\left(  M^{\ast}\right)  $ is \emph{not} an obstruction, we show the
following result. We say that $f\in L^{2}\left(  \mathbb{R}\right)  $ is
\emph{orthogonal} if the translates $\left\{  f\left(  \,\cdot\,-n\right)  \mid
n\in\mathbb{Z}\right\}  $ form an orthonormal family in $L^{2}\left(
\mathbb{R}\right)  $. Let $\varphi$ be a scaling function relative to some
$m_{0}$ satisfying (\ref{IntrAxiom(1)})--(\ref{IntrAxiom(3)}), and let $h\in
L^{2}\left(  \mathbb{R}\right)  $. Suppose both $\varphi$ and $h$ are
orthogonal in this sense.
Then we show that the following two conditions are
equivalent:\renewcommand{\theenumi}{\alph{enumi}}
\pagebreak
\begin{enumerate}
\item \label{IntrCond(1)}$\left\|  \varphi-M^{n}h\right\|  _{2}\underset
{n\rightarrow\infty}{\longrightarrow}0$, and

\item \label{IntrCond(2)}the function%
\[
p\left(  t\right)  =\sum_{n\in\mathbb{Z}}e^{int}\int_{\mathbb{R}}%
\overline{\varphi\left(  x+n\right)  }h\left(  x\right)  \,dx
\]
is continuous near $t=0$, and $p\left(  0\right)  =1$.
\end{enumerate}
The result fails when $\ker\left(  M^{\ast}\right)  \neq\left\{
0\right\}  $; see Section \ref{SiV}.

In studying this, and other related approximation problems in earlier papers
\cite{BrJo97b,BrJo98b}, the following general and geometric Hilbert-space
framework proved useful. Let $m_{0}$ be given subject to (\ref{IntrAxiom(1)}%
)--(\ref{IntrAxiom(3)}). (For much of it, only (\ref{IntrAxiom(3)}) is
needed.) We introduce the \emph{Ruelle operator}%
\begin{equation}
\left(  R\xi\right)  \left(  z\right)  =\frac{1}{2}\sum_{w^{2}=z}\left|
m_{0}\left(  w\right)  \right|  ^{2}\xi\left(  w\right)  ,\qquad\xi\in
L^{\infty}\left(  \mathbb{T}\right)  , \label{eqIntr.4}%
\end{equation}
and its dual,%
\[
\left(  R^{\ast}\xi\right)  \left(  z\right)  =\left|  m_{0}\left(  z\right)
\right|  ^{2}\xi\left(  z^{2}\right)  .
\]
Let $\mathcal{H}$ be an abstract Hilbert space, and let $\pi$ be a
representation of $L^{\infty}\left(  \mathbb{T}\right)  $ on $\mathcal{H}$. An
operator $M$ in $\mathcal{H}$ is said to be an
$\left( R,\pi\right) $\emph{-isometry} if%
\begin{equation}
M^{\ast}\pi\left(  \xi\right)  M=\pi\left(  R^{\ast}\left(  \xi\right)
\right)  \text{\qquad for all }\xi\in L^{\infty}\left(  \mathbb{T}\right)  .
\label{eqIntr.5}%
\end{equation}
We also call $M$ a sub-isometry when $\left( R,\pi\right) $ is understood.
The cascade operators are special cases of
$\left( R,\pi\right) $-isometries (see Section
\ref{Back}), and several of our results for the cascade operator will be
derived from a general result about
$\left( R,\pi\right) $-isometries; see Section \ref{R-is}. We
now summarize briefly our main result for
$\left( R,\pi\right) $-isometries. Let $M$ be an
$\left( R,\pi\right) $-isometry on some Hilbert space
$\mathcal{H}$, and let $\mathcal{L}%
:=\ker\left(  M^{\ast}\right)  $, and $\mathcal{H}^{\left(  \infty\right)
}:=\bigcap_{n=1}^{\infty}\left[  M^{n}\mathcal{H}\right]  $, where $\left[
\,\cdot\,\right]  $ refers to norm-closure in $\mathcal{H}$. We then establish
(Theorem \ref{ThmR-is.2}) the following general \emph{orthogonal}
decomposition,%
\begin{equation}
\mathcal{H}=\sideset{}{^{\smash{\oplus}}}{\sum}\limits_{n=0}^{\infty}\left[
M^{n}\mathcal{L}\right]  \oplus\mathcal{H}^{\left(  \infty\right)  },
\label{eqIntr.6}%
\end{equation}
and, moreover, we show that each of the (mutually orthogonal) component spaces
$\left[  M^{n}\mathcal{L}\right]  $ and $\mathcal{H}^{\left(  \infty\right)
}$ is invariant under $\pi\left(  \xi\right)  $ for all $\xi\in L^{\infty
}\left(  \mathbb{T}\right)  $. In other words, the spaces in the decomposition
reduce the representation $\pi$; if the restricted representations are denoted
$\pi_{n}$, we have $\pi=\pi_{0}\oplus\pi_{1}\oplus\dots\oplus\pi_{\infty}$.
(See \cite{Dix96} for the theory of decompositions of representations.) The
significance of having the terms $\mathcal{H}^{\left(  n\right)  }$ in the sum
(\ref{eqIntr.6}) invariant under the representation $\pi$ is that the
structure of the subspaces $\mathcal{H}^{\left(  n\right)  }$ themselves may
then be determined from the representation. We show that generically the\pagebreak
sequence $\pi_{n}=\pi|_{\mathcal{H}^{\left(  n\right)  }}$ is determined from
$\pi_{0}$, and the latter can be computed from the spectral theorem. We also
show that \emph{every} representation of $L^{\infty}\left(  \mathbb{T}\right)
$ occurs as $\pi_{0}$ in some decomposition (\ref{eqIntr.6}) corresponding to
an $\left( R,\pi\right) $-isometry.

While the operator $M$ plays a central role in the wavelet literature (see,
e.g., \cite{CoDa96}, \cite{CoRy95}, \cite{Dau92}), the approach
(\ref{eqIntr.5})--(\ref{eqIntr.6}) here is new. It is motivated by the need
for a representation-theoretic approach to the classification of wavelets
\cite{BrJo97b}, and also by the need for including in the analysis other
Hilbert spaces $\mathcal{H}$ than just $L^{2}\left(  \mathbb{R}\right)  $; see
\cite{Jor98}. Even if $L^{2}\left(  \mathbb{R}\right)  $ is the final goal,
there is a need for understanding the limiting cases when some different
Hilbert space $\mathcal{H}$ (other than $L^{2}\left(  \mathbb{R}\right)  $) is
dictated by some more general or different framework and analysis; we refer to
\cite{Jor98} for more details on this viewpoint. In \cite{Jor98}, a framework
is adopted which is much more general than the present setting of quadrature
mirror filters. (It includes, for example, the orthogonal harmonic analysis of
$\mathcal{H}=L^{2}\left(  \mu\right)  $ from \cite{JoPe98b} where $\mu$ is
taken to be a self-affine and singular probability measure on $\mathbb{R}^{d}$
of compact support which arises by iteration of a given finite set of affine
mappings in $\mathbb{R}^{d}$.) Our present approach to multiresolution
analysis is inspired by the Lax--Phillips scattering theory \cite{LaPh89} for
the classical wave equation. This is the continuous case; our aim here is to
``discretize'' the Lax--Phillips scattering theory, and (\ref{eqIntr.5})
should be viewed in that light.
The recent papers by Micchelli et al.\ \cite{MiPr89}, \cite{Mic96},
\cite{CDM91} take a somewhat different (but closely related) approach to the
discretization problem: Since the functions in a multiresolution subspace
$V\left(  \varphi\right)  $ may be represented by sequences $\left(  \xi
_{n}\right)  _{n\in\mathbb{Z}}$, via $\sum_{n}\xi_{n}\varphi\left(
\,\cdot\,-n\right)  $, the operator $M$ may therefore be studied alternatively
as acting on one of the sequence spaces $\ell^{p}\left(  \mathbb{Z}\right)  $.
In this guise, it takes the form $M\rightarrow S=\tilde{M}$, where%
\begin{equation}
\left(  S\xi\right)  _{n}=\sum_{k\in\mathbb{Z}}a_{n-2k}\xi_{k}.
\label{eqIntr.7}%
\end{equation}
The connection to the sequence spaces is taken up here in Section \ref{Zak}
below; see especially Theorem \ref{ThmZak.7} and Table \ref{ThmZak.7(table)}.
The operator on $L^{2}\left(  \mathbb{T}\right)  $ which corresponds to
(\ref{eqIntr.7}) is simply $f\mapsto m\left(  z\right)  f\left(  z^{2}\right)
$ when $m\left(  z\right)  =\sum_{n\in\mathbb{Z}}a_{n}z^{n}$ is the Fourier
series of some given filter function $m$.

The motivation for the present paper has several sources. In \cite{BEJ97}, we
showed that certain wavelets, for which the corresponding quadrature mirror
filter $m_{0}$ is a polynomial, may be classified by the labels of a
corresponding family of irreducible representations of the Cuntz algebra
$\mathcal{O}_{2}$. Previously there were known no such clear-cut invariants
that classified wavelet families. But the particular wavelets from
\cite{BEJ97} had in fact been identified earlier (without invariants or
classification) in \cite{Wel93} and \cite{Wic93}. Regularity questions for the
corresponding scaling functions had been studied in \cite{DDL95},
\cite{GMW94}, and \cite{MRV96}. Our references for quadrature mirror filters
are \cite{StNg96} and \cite{Mal98}. Excellent references to wavelets from the
operator-theoretic viewpoint are \cite{HeWe96} and \cite{Hor95}. The
refinement operators of (\ref{eqIntr.7}) are also called \emph{slant Toeplitz
operators,} and their spectral theory was studied previously in \cite{Ho96}. A
standard reference to the ergodic theory which we use in the present paper is
\cite{Wal82}.
\pagebreak

\section[Background, summary, and operator relations]{\label{Back}BACKGROUND, SUMMARY, AND OPERATOR RELATIONS}

Let $\mathbb{T}$ denote the one-dimensional torus
and let $\mu$ be its (usual) Haar measure.
We will use the form
$\mathbb{T}\simeq\mathbb{R}\diagup2\pi\mathbb{Z}$ such that functions on
$\mathbb{T}$ are identified with $2\pi$-periodic functions on $\mathbb{R}$.
For technical reasons, the identification will be made via $z=e^{-i\omega}$,
$\omega\in\mathbb{R}$, and functions on $\mathbb{T}$ will be written,
alternately, as $f\left(  z\right)  $ or $f\left(  \omega\right)  $. Let
$m_{0}\in L^{\infty}\left(  \mathbb{T}\right)  $ be given, and suppose
$m_{0}=\sqrt{2}$ at $z=1$, and further that%
\begin{equation}
\left|  m_{0}\left(  z\right)  \right|  ^{2}+\left|  m_{0}\left(  -z\right)
\right|  ^{2}=2. \label{eqBack.1}%
\end{equation}
One approach to wavelets (see \cite{Dau92}) is to first make precise the
(formal) infinite product (limit as $n\rightarrow\infty$):%
\begin{equation}
\prod_{k=1}^{n}\frac{1}{\sqrt{2}\mathstrut}m_{0}\left(  \frac{\omega
}{2^{k\mathstrut}}\right)  \chi_{\left[  -2^{n}\pi,2^{n}\pi\right]  }^{{}%
}\left(  \omega\right)  . \label{eqBack.2}%
\end{equation}
Suppose this product (\ref{eqBack.2}), in the limit, represents a function
$F\in L^{2}\left(  \mathbb{R}\right)  $. Then%
\begin{equation}
\sqrt{2}F\left(  2\omega\right)  =m_{0}\left(  \omega\right)  F\left(
\omega\right)  ,\qquad\omega\in\mathbb{R}, \label{eqBack.3}%
\end{equation}
and we wish to recover the scaling function $\varphi$ as the inverse
Fourier transform of $F$.
The approach is referred to as the \emph{Mallat algorithm} \cite{Mal89}, but
it is somewhat indirect. Let
\begin{equation}
\varphi\left(  x\right)  =\frac{1}{2\pi}\int_{-\infty}^{\infty}e^{i\omega
x}F\left(  \omega\right)  \,d\omega\label{eqBack.4}%
\end{equation}
be the inverse Fourier transform, $\varphi=F\spcheck$, and let%
\begin{equation}
m_{0}\left(  z\right)  =\sum_{n\in\mathbb{Z}}a_{n}z^{n} \label{eqBack.5}%
\end{equation}
be the Fourier series expansion of $m_{0}$. (The Fourier basis in
(\ref{eqBack.5}) is $e_{n}\left(  z\right)  =z^{n}$.) Then (again formally),%
\begin{equation}
\varphi\left(  x\right)  =\sqrt{2}\sum_{n\in\mathbb{Z}}a_{n}\varphi\left(
2x-n\right)  . \label{eqBack.6}%
\end{equation}
The function $\varphi$ is called the \emph{scaling function,} and the closed
linear span of the translates%
\begin{equation}
\left\{  \varphi\left(  \,\cdot\,-n\right)  \mid n\in\mathbb{Z}\right\}
\label{eqBack.7}%
\end{equation}
generates, under certain analytic conditions \cite{Dau92}, a
\emph{multiresolution subspace} $V=V\left(  \varphi\right)  $ in $L^{2}\left(
\mathbb{R}\right)  $. The identity (\ref{eqBack.6})
allows,
for each $m_{0}$, an operator $M=M_{m_{0}}$ in $L^{2}\left(
\mathbb{R}\right)  $ such that scaling functions $\varphi$ arise as solutions
to a fixed point problem.
If we introduce the operator $M$ by
\begin{equation}
\left(  Mh\right)  \left(  x\right)  =\sqrt{2}\sum_{n\in\mathbb{Z}}%
a_{n}h\left(  2x-n\right)  ,\qquad h\in L^{2}\left(  \mathbb{R}\right)  ,
\label{eqBack.9}%
\end{equation}
\pagebreak
then
\begin{equation}
\varphi=M\varphi , \label{eqBack.8}%
\end{equation}
and we should look for choices of $h$ such that%
\begin{equation}
\lim_{n\rightarrow\infty}M^{n}h=\varphi. \label{eqBack.10}%
\end{equation}
We will be interested here primarily in this as an $L^{2}\left(
\mathbb{R}\right)  $-limit, but other limits (e.g., pointwise) will be
considered as well. The issue is both how preassigned properties of the
\emph{starting function} $h$, and the scaling function $\varphi$, affect the
approximation properties of (\ref{eqBack.10}). We refer to this as the
\emph{cascade approximation,} and the traditional choice for $h$ is
$h=\chi_{\left[  0,1\right]  }^{{}}$, which accidentally is the scaling
function for the Haar wavelet. The cascade algorithm is more direct than the
Mallat one, as it doesn't pass via the Fourier transform.

The
right-hand side in (\ref{eqBack.9}) involves a dyadic scaling, and an action
by $L^{\infty}\left(  \mathbb{T}\right)  $ on $L^{2}\left(  \mathbb{R}\right)
$.
We now introduce a
representation $\alpha\mapsto\pi\left( \alpha\right) $
of $L^{\infty}\left( \mathbb{T}\right) $ in
the
algebra of operators on
$L^{2}\left( \mathbb{R}\right) $ by
\begin{equation}
\left( \pi\left( \alpha\right)  h \right) \left( x\right) = \sum_{n\in\mathbb{Z}} a_{n}h\left( x-n\right) .
\label{eqBack.11}
\end{equation}
Then $M$ can be represented as $M=D\pi\left( m_0\right) $, where $D$ is the operator
of
dyadic scaling, $Dh\left( x\right) = \sqrt{2} h\left( 2x\right) $.
Application of the Fourier transform $\widehat{\,\cdot\,}$ to both sides in
(\ref{eqBack.11}) then yields%
\begin{equation}
\hat{h}\longmapsto\alpha\left(  \omega\right)  \hat{h}\left(  \omega\right)
\label{eqBack.12}%
\end{equation}
via the identification%
\[
\alpha\left(  \omega\right)  \simeq\alpha\left(  z\right)  ,\qquad
z=e^{-i\omega},\;\omega\in\mathbb{R}.
\]
(In applications, $\omega$ is a frequency variable.) In either of its forms,
(\ref{eqBack.11}) or (\ref{eqBack.12}), this representation of $L^{\infty
}\left(  \mathbb{T}\right)  $ will be denoted $\pi\left(  \alpha\right)  h$
with $\pi\left(  \alpha\right)  $ an operator acting on $L^{2}\left(
\mathbb{R}\right)  $, and $\pi\left(  L^{\infty}\left(  \mathbb{T}\right)
\right)  $ the corresponding algebra of operators.

If $m_{0}$ is given, then $L^{\infty}\left(  \mathbb{T}\right)  $ carries two
operators $R$ and $R^{\ast}$, $R$ defined by%
\begin{equation}
\left(  R\alpha\right)  \left(  z\right)  =\frac{1}{2}\sum_{\substack{w\in
\mathbb{T}\\w^{2}=z}}\left|  m_{0}\left(  w\right)  \right|  ^{2}\alpha\left(
w\right)  \label{eqBack.13}%
\end{equation}
and the adjoint operator $R^{\ast}$ by%
\begin{equation}
\left(  R^{\ast}\beta\right)  \left(  z\right)  =\left|  m_{0}\left(
z\right)  \right|  ^{2}\beta\left(  z^{2}\right)  ,\qquad\alpha,\beta\in
L^{\infty}\left(  \mathbb{T}\right)  . \label{eqBack.14}%
\end{equation}
Both of them will also be viewed as $L^{2}\left(  \mathbb{T}\right)
$-operators, and we have, by a simple calculation,%
\begin{equation}
\int_{\mathbb{T}}\left(  R\alpha\right)  \left(  z\right)  \beta\left(
z\right)  \,d\mu\left(  z\right)  =\int_{\mathbb{T}}\alpha\left(  z\right)
R^{\ast}\beta\left(  z\right)  \,d\mu\left(  z\right)  , \label{eqBack.15}%
\end{equation}
\pagebreak
where $d\mu$ is the usual Haar measure on $\mathbb{T}$, i.e., $\frac{1}{2\pi
}\int_{0}^{2\pi}\cdots\,d\omega$, thus justifying the notation $R^{\ast}$. The
operator $R$ is called the \emph{Ruelle operator,} or the \emph{transfer
operator,} for reasons we shall go into later; see also \cite{CoRy95}.

We further have the usual pairing between $M$ in (\ref{eqBack.9}), and its
$L^{2}\left(  \mathbb{R}\right)  $-adjoint $M^{\ast}$, given by
\[
\int_{\mathbb{R}}\overline{Mh_{1}\left(  x\right)  }h_{2}\left(  x\right)
\,dx=\int_{\mathbb{R}}\overline{h_{1}\left(  x\right)  }M^{\ast}h_{2}\left(
x\right)  \,dx,
\]
or equivalently,%
\begin{equation}
\ip{Mh_{1}}{h_{2}}=\ip{h_{1}}{M^{\ast}h_{2}}, \label{eqBack.16}%
\end{equation}
where $\ip{\,\cdot\,}{\,\cdot\,}$ denotes the standard inner product of
$L^{2}\left(  \mathbb{R}\right)  $.

Our first result is

\begin{theorem}
\label{ThmBack.1}Let $m_{0}\in L^{\infty}\left(  \mathbb{T}\right)  $ be
given, and suppose it satisfies \textup{(\ref{eqBack.1})}. Let $M$ be the
corresponding cascade operator, and $R$ the Ruelle operator. The respective
adjoints are $M^{\ast}$ and $R^{\ast}$. Finally, let $\pi$ be the
representation of $L^{\infty}\left(  \mathbb{T}\right)  $ on $L^{2}\left(
\mathbb{R}\right)  $ given in \textup{(\ref{eqBack.12})}. Then we have the
following two commutation relations:

\begin{enumerate}
\item \label{ThmBack.1(1)}$M^{\ast}\pi\left(  \alpha\right)  M=\pi\left(
R^{\ast}\alpha\right)  $ and

\item \label{ThmBack.1(2)}$M\pi\left(  \alpha\right)  M^{\ast}=\pi\left(
R\alpha\right)  +\pi\left(  R\left(  e_{1}\alpha\right)  \right)  T_{\frac
{1}{2}}$ for all $\alpha\in L^{\infty}\left(  \mathbb{T}\right)  $, where
$e_{1}\left(  z\right)  =z$ and $\left(  T_{\frac{1}{2}}h\right)  \left(
x\right)  =h\left(  x+\frac{1}{2}\right)  $, $h\in L^{2}\left(  \mathbb{R}%
\right)  $.
\end{enumerate}
\end{theorem}

The proof will be given in Section \ref{Poof} below. In this paper, we will
study the cascade approximation, and the scaling function $\varphi$, via the
abstract algebraic system which is given by the relations (\ref{ThmBack.1(1)}%
)--(\ref{ThmBack.1(2)}) of the theorem. These two operator commutation
relations will first be studied abstractly (Section \ref{R-is}) and
independently of their origin, and then the results will be specialized, and
applied, to the wavelet problems mentioned above. Our proofs will depend on
some lemmas of a general nature regarding the Zak transform.

\section[The Zak transform]{\label{Zak}THE ZAK TRANSFORM}

The Zak transform $Z$ is known \cite[p.\ 109]{Dau92} to be the \emph{isometric
isomorphism} between $L^{2}\left(  \mathbb{R}\right)  $ and $L^{2}\left(
\mathbb{T}\times\left[  0,1\right]  \right)  $ which is given formally by%
\begin{equation}
\left(  Zh\right)  \left(  z,x\right)  =\sum_{n\in\mathbb{Z}}z^{n}h\left(
x+n\right)  ,\qquad h\in L^{2}\left(  \mathbb{R}\right)  ,\;x\in\mathbb{R}.
\label{eqZak.1}%
\end{equation}
It is studied in \cite[p.\ 109]{Dau92}, \cite{BJR97}, and elsewhere. To make
it precise, it is convenient to identify its range with functions $H$ on
$\mathbb{T}\times\mathbb{R}$ which satisfy the following scaling rule:%
\begin{equation}
H\left(  z,x+n\right)  =z^{-n}H\left(  z,x\right)  \text{\qquad for all }%
z\in\mathbb{T},\;x\in\mathbb{R},\;n\in\mathbb{Z}. \label{eqZak.2}%
\end{equation}
\pagebreak
The norm $\left\|  H\right\|  $, or $\left\|  H\right\|  _{2}$, will be given
by%
\begin{equation}
\left\|  H\right\|  ^{2}=\int_{\mathbb{T}}\int_{0}^{1}\left|  H\left(
z,x\right)  \right|  ^{2}\,d\mu\left(  z\right)  \,dx, \label{eqZak.3}%
\end{equation}
and it can be checked, see \cite[p.\ 109]{Dau92}, that%
\begin{align}
\left\|  H\right\|  ^{2}  &  =\left\|  h\right\|  ^{2}\label{eqZak.4}\\
&  =\int_{\mathbb{R}}\left|  h\left(  x\right)  \right|  ^{2}\,dx,\nonumber
\end{align}
and that $H=Zh$ is an isometric isomorphism of $L^{2}\left(  \mathbb{R}%
\right)  \emph{onto}$ the Hilbert space $\mathcal{H}_{Z}$ of functions which
are defined by the scaling formula (\ref{eqZak.2}), and completion in the norm
$\left\|  \,\cdot\,\right\|  _{2}$ of (\ref{eqZak.3}).
The simplest wavelet scaling function is $\varphi_{H}:=\chi_{\left[
0,1\right]  }^{{}}$ of the Haar wavelet. (The Haar wavelet itself is generated
by $\psi_{H}=\chi_{\left[  0,\frac{1}{2}\right]  }^{{}}-\chi_{\left[  \frac
{1}{2},1\right]  }^{{}}$.) The Zak transform $Z$ has the pleasant feature that
$Z\left(  \varphi_{H}\right)  $ is the constant function $\openone$ in
$L^{2}\left(  \mathbb{T}\times\left[  0,1\right]  \right)  $.

A second advantage of the $\mathcal{H}_{Z}$ formulation (\ref{eqZak.2}) is
that it makes clear a useful \emph{direct integral representation}%
\begin{equation}
\mathcal{H}_{Z}=\int_{\mathbb{T}}^{\oplus}\mathcal{H}\left(  z\right)
\,d\mu\left(  z\right)  , \label{eqZak.5}%
\end{equation}
where each $\mathcal{H}\left(  z\right)  $ is the Hilbert space $\left\{
h_{z}\right\}  $ of functions $h_{z}$ on $\mathbb{R}$ satisfying the
$z$-scaling rule,%
\[
h_{z}\left(  x+n\right)  =z^{-n}h_{z}\left(  x\right)  \text{,\qquad for all
}x\in\mathbb{R},\;n\in\mathbb{Z}.
\]
We refer to \cite{Dix69} for details on direct integrals, and note that each
$\mathcal{H}\left(  z\right)  $ is naturally and isometrically isomorphic to
$L^{2}\left(  0,1\right)  $. Thus $Z$ is a \emph{transform}
representing $L^{2}\left(  \mathbb{R}\right)  $-functions as direct integrals
(over $\mathbb{T}$ with Haar measure) of copies of $L^{2}\left(  0,1\right)  $.

The following lemmas about the Zak transform are new. They are needed in the
sequel, and are also of independent interest, we feel.

\begin{lemma}
\label{LemZak.1}Let $m_{0}\in L^{\infty}\left(  \mathbb{T}\right)  $ be given
satisfying \textup{(\ref{eqBack.1})}, and let $M$ and $R$ denote the
corresponding cascade operator and Ruelle operator. Let $h_{i}$, $i=1,2$, be
functions in $L^{2}\left(  \mathbb{R}\right)  $ with Zak transforms
$H_{i}=Zh_{i}$. Then%
\begin{equation}
R\left(  \ip{H_{1}\left( z\right) }{H_{2}\left( z\right) }\right)  \left(
z\right)  =\ip{Z_{z}Mh_{1}}{Z_{z}Mh_{2}}, \label{eqZak.6}%
\end{equation}
where the inner product $\ip{\,\cdot\,}{\,\cdot\,}$ on both sides of
\textup{(\ref{eqZak.6})} is that of $\mathcal{H}\left(  z\right)  $, i.e.,%
\begin{equation}
\ip{H_{1}\left( z\right) }{H_{2}\left( z\right) }=\int_{0}^{1}\overline
{H_{1}\left(  z,x\right)  }H_{2}\left(  z,x\right)  \,dx, \label{eqZak.7}%
\end{equation}
and%
\begin{equation}
R\left(  \ip{H_{1}\left( \,\cdot\,\right) }{H_{2}\left( \,\cdot\,\right
) }\right)  \left(  z\right)  =\frac{1}{2}\sum_{w^{2}=z}\left|  m_{0}\left(
w\right)  \right|  ^{2}\ip{H_{1}\left( w\right) }{H_{2}\left( w\right) },
\label{eqZak.8}%
\end{equation}
and the summation in \textup{(\ref{eqZak.8})} is over the two roots,
$w\in\left\{  \pm\sqrt{z}\right\}  \subset\mathbb{T}$.
\end{lemma}
\pagebreak

\begin{remark}
\label{RemZak.2}The formula \textup{(\ref{eqZak.6})} is crucial for the use of
$R$ as a transfer operator in the iteration of the cascade algorithm on a
given starting vector $h\in L^{2}\left(  \mathbb{R}\right)  $, i.e.,
$h\rightarrow Mh\rightarrow M^{2}h\rightarrow\cdots$. Secondly, we will show
in Section \textup{\ref{R-is}} that \textup{(\ref{eqZak.6})} is equivalent to
the first of the two commutation relations, \textup{(\ref{ThmBack.1(1)}),} in
Theorem \textup{\ref{ThmBack.1}.}
\end{remark}

\begin{proof}
Let the functions be as described in the lemma, $H_{i}=Zh_{i}$, $i=1,2$. We
first calculate the term $Z_{z}Mh_{1}$ from the desired equation
(\ref{eqZak.6}). Keep in mind that all integrals and summations are convergent
relative to the respective Hilbert-space norms (due to the isometries which we
described before the statement of the lemma). Details are in \cite[p.\ 109]%
{Dau92} and \cite{BJR97}.

Let%
\begin{equation}
\left(  Uh_{1}\right)  \left(  x\right)  :=\frac{1}{\sqrt{2}\mathstrut}%
\,h_{1}\left(  \frac{x}{2}\right)  . \label{eqNewZak.8}%
\end{equation}
Then $Mh_{1}=U^{-1}\pi\left(  m_{0}\right)  h_{1}$, and we now calculate the
Zak transform of the two individual operators making up $M$. First,
$h_{1}\mapsto 
\pi\left(  m_{0}\right)  h_{1}$ transforms into $H_{1}\mapsto m_{0}\left(
z\right)  H_{1}\left(  z,\,\cdot\,\right)  $, since
\begin{align*}
Z\left(  
\pi\left(  m_{0}\right)  h_{1}\right)   &  =\sum_{n}z^{n}\left(  
\pi\left(  m_{0}\right)  h_{1}\right)  \left(  x+n\right) \\
&  =\sum_{n}\sum_{k}z^{n}a_{k}h_{1}\left(  x+n-k\right) \\
&  =\sum_{k}a_{k}z^{k}\sum_{n}z^{n-k}h_{1}\left(  x+n-k\right) \\
&  =\sum_{k}a_{k}z^{k}\sum_{n}z^{n}h_{1}\left(  x+n\right) \\
&  =m_{0}\left(  z\right)  H_{1}\left(  z,x\right)
\end{align*}
as claimed. Note that all summation indices range over $n,k\in\mathbb{Z}$. We
have $\ell^{2}\left(  \mathbb{Z}\right)  $-convergent summations, i.e.,
relative to the respective $\ell^{2}$-norms, and the exchange of the
summations is justified by the norm-isomorphism property of the Zak transform.

We now turn to the operator $\tilde{U}^{-1}$ given by $\tilde{U}^{-1}%
=ZU^{-1}Z^{\ast}$, where $Z^{\ast}$ is the adjoint of $Z$. Since $Z$ is a
norm-isomorphism, we have%
\[
ZZ^{\ast}=\operatorname*{id}\nolimits_{L^{2}\left(  \mathbb{T}\times\left[
0,1\right]  \right)  }\,,
\text{\qquad and\qquad}
Z^{\ast}Z=\operatorname*{id}\nolimits_{L^{2}\left(  \mathbb{R}\right)  }\,,
\]
where $\operatorname*{id}$ refers to the respective identity operators.

We claim that%
\begin{equation}
\left(  Z^{\ast}H_{1}\right)  \left(  x\right)  =\int_{\mathbb{T}}H_{1}\left(
z,x\right)  \,d\mu\left(  z\right)  . \label{eqZak.9}%
\end{equation}
\pagebreak
The proof is the following calculation, for $h\in L^{2}\left(  \mathbb{R}%
\right)  $:%
\begin{align*}
\ip{H_{1}}{Zh}  &  =\int_{\mathbb{T}}\int_{0}^{1}\overline{H_{1}\left(
z,x\right)  }Zh\left(  z,x\right)  \,d\mu\left(  z\right)  \,dx\\
&  =\int_{\mathbb{T}}\int_{0}^{1}\overline{H_{1}\left(  z,x\right)  }\sum
_{n}z^{n}h\left(  x+n\right)  \,d\mu\left(  z\right)  \,dx\\
&  =\int_{\mathbb{T}}\int_{0}^{1}\sum_{n}\overline{H_{1}\left(  z,x+n\right)
}h\left(  x+n\right)  \,dx\,d\mu\left(  z\right) \\
&  =\int_{\mathbb{T}}\int_{-\infty}^{\infty}\overline{H_{1}\left(  z,x\right)
}h\left(  x\right)  \,dx\,d\mu\left(  z\right) \\
&  =\int_{\mathbb{R}}\overline{\int_{\mathbb{T}}H_{1}\left(  z,x\right)
\,d\mu\left(  z\right)  }\,h\left(  x\right)  \,dx\\
&  =\vphantom{\int}\ip{Z^{\ast}H_{1}}{h},
\end{align*}
which proves the stated formula (\ref{eqZak.9}) for $Z^{\ast}$.

We now calculate%
\begin{align*}
\left(  \tilde{U}^{-1}H_{1}\right)  \left(  z,x\right)   &  =\left(
ZU^{-1}Z^{\ast}\right)  H_{1}\left(  z,x\right) \\
&  =\sum_{n}z^{n}\left(  U^{-1}Z^{\ast}H_{1}\right)  \left(  x+n\right) \\
&  =\sum_{n}z^{n}\sqrt{2}\left(  Z^{\ast}H_{1}\right)  \left(  2x+2n\right) \\
&  =\sum_{n}z^{n}\sqrt{2}\int_{\mathbb{T}}H_{1}\left(  \zeta,2x+2n\right)
\,d\mu\left(  \zeta\right) \\
&  =\sqrt{2}\sum_{n}z^{n}\int_{\mathbb{T}}\zeta^{-2n}H_{1}\left(
\zeta,2x\right)  \,d\mu\left(  \zeta\right) \\
&  =\sqrt{2}\sum_{n}z^{n}\int_{\mathbb{T}}\zeta^{-n}\frac{1}{2}\sum
_{\substack{\xi\in\mathbb{T}\\\xi^{2}=\zeta}}H_{1}\left(  \xi,2x\right)
\,d\mu\left(  \zeta\right) \\
&  =\frac{1}{\sqrt{2}\mathstrut}\sum_{w^{2}=z}H_{1}\left(  w,2x\right)  ,
\end{align*}
where the last step simply represents the Fourier series of the final
function. The formula%
\begin{equation}
\tilde{U}^{-1}H_{1}\left(  z,x\right)  =\frac{1}{\sqrt{2}\mathstrut}%
\sum_{w^{2}=z}H_{1}\left(  w,2x\right)  \label{eqZak.10}%
\end{equation}
is basic in later proofs. Combining the three formulas (\ref{eqZak.8}),
(\ref{eqZak.9}), and (\ref{eqZak.10}), we arrive at%
\begin{equation}
\left(  Z_{z}Mh_{1}\right)  \left(  x\right)  =\frac{1}{\sqrt{2}\mathstrut
}\sum_{w^{2}=z}m_{0}\left(  w\right)  H_{1}\left(  w,2x\right)  .
\label{eqZak.11}%
\end{equation}
This is needed in the calculation of the right-hand side in (\ref{eqZak.6})
from the lemma as follows:%
\begin{align}
&  \ip{Z_{z}Mh_{1}}{Z_{z}Mh_{2}}\label{eqZak.12}\\
&  \quad=\int_{0}^{1}\overline{\left(  Z_{z}Mh_{1}\right)  \left(  x\right)
}\left(  Z_{z}Mh_{2}\right)  \left(  x\right)  \,dx\nonumber\\
&  \quad=\frac{1}{2}\int_{0}^{1}\sum_{w_{1}^{2}=z}\overline{m_{0}\left(
w_{1}\right)  }\,\overline{H_{1}\left(  w_{1},2x\right)  }\sum_{w_{2}^{2}%
=z}m_{0}\left(  w_{2}\right)  H_{2}\left(  w_{2},2x\right)  \,dx\nonumber\\
&  \quad=\frac{1}{2}\mathop{\sum\sum}_{w_{1}^{2}=w_{2}^{2}=z}\overline
{m_{0}\left(  w_{1}\right)  }m_{0}\left(  w_{2}\right)  \int_{0}^{1}%
\overline{H_{1}\left(  w_{1},2x\right)  }H_{2}\left(  w_{2},2x\right)
\,dx\nonumber\\[4pt]
&  \quad=\frac{1}{4}\mathop{\sum\sum}_{w_{1}^{2}=w_{2}^{2}=z}\overline
{m_{0}\left(  w_{1}\right)  }m_{0}\left(  w_{2}\right) \nonumber\\
&  \qquad\cdot\int_{0}^{1}\left[  \overline{H_{1}\left(  w_{1},x\right)
}H_{2}\left(  w_{2},x\right)  +\overline{H_{1}\left(  w_{1},x+1\right)  }%
H_{2}\left(  w_{2},x+1\right)  \right]  \,dx\nonumber\\[4pt]
&  \quad=\frac{1}{4}\mathop{\sum\sum}_{w_{1}^{2}=w_{2}^{2}=z}\overline
{m_{0}\left(  w_{1}\right)  }m_{0}\left(  w_{2}\right) \nonumber\\
&  \qquad\cdot\int_{0}^{1}\left[  \overline{H_{1}\left(  w_{1},x\right)
}H_{2}\left(  w_{2},x\right)  +w_{1}^{{}}w_{2}^{-1}\overline{H_{1}\left(
w_{1},x\right)  }H_{2}\left(  w_{2},x\right)  \right]  \,dx\nonumber\\[4pt]
&  \quad=\frac{1}{4}\mathop{\sum\sum}_{w_{1}^{2}=w_{2}^{2}=z}\overline
{m_{0}\left(  w_{1}\right)  }m_{0}\left(  w_{2}\right)  \left(  1+w_{1}^{{}%
}w_{2}^{-1}\right)  \int_{0}^{1}\overline{H_{1}\left(  w_{1},x\right)  }%
H_{2}\left(  w_{2},x\right)  \,dx.\nonumber
\end{align}
But if $w_{1}\neq w_{2}$ in the summation, then the factor $\left(
1+w_{1}^{{}}w_{2}^{-1}\right)  =0$. If $w_{1}\neq w_{2}$, then $w_{2}=-w_{1}$
and $w_{1}^{{}}w_{2}^{-1}=-1$. If $w_{1}=w_{2}$, then $\left(  1+w_{1}^{{}%
}w_{2}^{-1}\right)  =2$. Continuing the calculation, we get%
\begin{align*}
\ip{Z_{z}Mh_{1}}{Z_{z}Mh_{2}}  &  =\frac{1}{2}\sum_{w^{2}=z}\left|
m_{0}\left(  w\right)  \right|  ^{2}\int_{0}^{1}\overline{H_{1}\left(
w,x\right)  }H_{2}\left(  w,x\right)  \,dx\\
&  =\frac{1}{2}\sum_{w^{2}=z}\left|  m_{0}\left(  w\right)  \right|  ^{2}%
\ip{H_{1}\left( w\right) }{H_{2}\left( w\right) }\\
&  =R\left(  \ip{H_{1}\left( \,\cdot\,\right) }{H_{2}\left( \,\cdot\,\right
) }\right)  \left(  z\right)  ,
\end{align*}
which is the final conclusion of the lemma.
\end{proof}

Having formula (\ref{eqZak.10}) for $\tilde{U}^{-1}$, we may derive the
corresponding formula for $\tilde{U}=ZUZ^{\ast}$, but we isolate it here in a
separate lemma.

\begin{lemma}
\label{LemZak.3}Let $U$ be the scaling operator \textup{(\ref{eqNewZak.8})} on
$L^{2}\left(  \mathbb{R}\right)  $, $Uh\left(  x\right)  =\frac{1}{\sqrt
{2}\mathstrut}\,h\left(  \frac{x}{2}\right)  $, and let $\tilde{U}=ZUZ^{\ast}$
be the corresponding operator on $\mathcal{H}_{Z}\simeq L^{2}\left(
\mathbb{T}\times\left[  0,1\right]  \right)  $ of the Zak transform. Then%
\begin{equation}
\left(  \tilde{U}H\right)  \left(  z,x\right)  =\frac{1}{\sqrt{2}\mathstrut
}\left(  H\left(  z^{2},\frac{x}{2}\right)  +zH\left(  z^{2},\frac{x+1}%
{2}\right)  \right)  . \label{eqZak.13}%
\end{equation}
\end{lemma}

\begin{proof}%
\begin{align*}
&  ZUZ^{\ast}H\left(  z,x\right) \\
&  \quad=\sum_{n}z^{n}\left(  UZ^{\ast}H\right)  \left(  x+n\right) \\
&  \quad=\frac{1}{\sqrt{2}\mathstrut}\sum_{n}z^{n}\left(  Z^{\ast}H\right)
\left(  \frac{x+n}{2}\right) \\
&  \quad=\frac{1}{\sqrt{2}\mathstrut}\sum_{n}z^{n}\int_{\mathbb{T}}H\left(
\zeta,\frac{x+n}{2}\right)  \,d\mu\left(  \zeta\right) \\
&  \quad=\frac{1}{\sqrt{2}\mathstrut}\sum_{k}z^{2k}\int_{\mathbb{T}}H\left(
\zeta,\frac{x}{2}+k\right)  \,d\mu\left(  \zeta\right)  +\frac{1}{\sqrt
{2}\mathstrut}\sum_{k}z^{2k+1}\int_{\mathbb{T}}H\left(  \zeta,\frac{x+1}%
{2}+k\right)  \,d\mu\left(  \zeta\right) \\
&  \quad=\frac{1}{\sqrt{2}\mathstrut}\sum_{k}\left(  z^{2k}\int_{\mathbb{T}%
}\zeta^{-k}H\left(  \zeta,\frac{x}{2}\right)  \,d\mu\left(  \zeta\right)
+z^{2k+1}\int_{\mathbb{T}}\zeta^{-k}H\left(  \zeta,\frac{x+1}{2}\right)
\,d\mu\left(  \zeta\right)  \right) \\
&  \quad=\frac{1}{\sqrt{2}\mathstrut}\left(  H\left(  z^{2},\frac{x}%
{2}\right)  +zH\left(  z^{2},\frac{x+1}{2}\right)  \right)  ,
\end{align*}
where, in the last step, we used the standard Fourier series representation of
the respective functions.
\end{proof}

Our interest in the use of the Zak-transform approach to cascade approximation
started in an earlier paper \cite{BrJo98b} on the special case of compactly
supported scaling functions. It is known \cite{Dau92} that compact support of
$\varphi$ (as a scaling function in $L^{2}\left(  \mathbb{R}\right)  $) is
equivalent to the filter function $m_{0}\left(  \,\cdot\,\right)  $ being a
polynomial, i.e., the sum in (\ref{eqBack.5}) being finite.

We encountered there the following two sesquilinear forms, both indexed by
$\mathbb{T}$:%
\begin{equation}
p_{1}\left(  h_{1},h_{2}\right)  \left(  e^{-i\omega}\right)  =\sum
_{n\in\mathbb{Z}}\overline{\hat{h}_{1}\left(  \omega+2\pi n\right)  }\hat
{h}_{2}\left(  \omega+2\pi n\right)  , \label{eqZak.14}%
\end{equation}
where the functions $h_{1}$, $h_{2}$ are in $L^{2}\left(  \mathbb{R}\right)  $
\emph{and} of compact support in the $x$-variable. This puts the Fourier
transform $\hat{h}_{i}$, $i=1,2$, in the Schwartz class so that the sum in
(\ref{eqZak.14}) is well defined. But problem was how to most efficiently
remove the compact-support restriction on the functions $h_{i}$. (Note the
individual terms in the sum on the right-hand side in (\ref{eqZak.14}) are not
periodic. But the sum $\sum_{n}$ serves to ``periodize'' the function
$\omega\mapsto\left(  \Bar{\Hat{h}}_{1}\hat{h}_{2}\right)  \left(
\omega\right)  $, $\omega\in\mathbb{R}$.)

The other sesquilinear form was%
\begin{equation}
p_{2}\left(  h_{1},h_{2}\right)  \left(  z\right)  =\sum_{k\in\mathbb{Z}}%
z^{k}\int_{\mathbb{R}}\overline{h_{1}\left(  x-k\right)  }h_{2}\left(
x\right)  \,dx. \label{eqZak.15}%
\end{equation}
With the compact-support restriction, this is even a \emph{finite} sum. But
introducing the Poisson summation formula, or by a direct Fourier series
expansion of $p_{1}\left(  h_{1},h_{2}\right)  $, we note that, with
$z=e^{-i\omega}$, we have%
\begin{equation}
p_{1}\left(  h_{1},h_{2}\right)  \left(  z\right)  =p_{2}\left(  h_{1}%
,h_{2}\right)  \left(  z\right)  . \label{eqZak.16}%
\end{equation}
When $p_{1}\left(  h_{1},h_{2}\right)  $ is viewed as a function on
$\mathbb{T}$, its $k$'th Fourier series coefficient computes out directly to
be $\int_{\mathbb{R}}\overline{h_{1}\left(  x-k\right)  }h_{2}\left(
x\right)  \,dx$, and (\ref{eqZak.16}) follows from this.

The following result makes it clear that the compact-support restriction can
be removed by use of the Zak-transform approach, and as a bonus, we get some
\emph{a priori} estimates that are needed later.

\begin{proposition}
\label{ProZak.4}Let $h_{i}\in L^{2}\left(  \mathbb{R}\right)  $ be of compact
support. Then the two forms $p_{1}$ and $p_{2}$ coincide with $p_{3}$, where%
\begin{equation}
p_{3}\left(  h_{1},h_{2}\right)  \left(  z\right)  =\ip{Z_{z}h_{1}}{Z_{z}%
h_{2}}=\int_{0}^{1}\overline{Zh_{1}\left(  z,x\right)  }Zh_{2}\left(
z,x\right)  \,dx. \label{eqZak.17}%
\end{equation}
\end{proposition}

\begin{proof}
In the following calculation, convergence is governed by the norm-isometric
property of $Z$, and this also justifies the exchange of summations and
integration:%
\begin{align*}
\ip{Z_{z}h_{1}}{Z_{z}h_{2}}  &  =\int_{0}^{1}\sum_{k}\overline{z^{k}%
h_{1}\left(  x+k\right)  }\sum_{l}z^{l}h_{2}\left(  x+l\right)  \,dx\\
&  =\sum_{k}\sum_{l}z^{l-k}\int_{0}^{1}\overline{h_{1}\left(  x+k\right)
}h_{2}\left(  x+l\right)  \,dx\\
&  =\sum_{n}z^{n}\sum_{l}\int_{0}^{1}\overline{h_{1}\left(  x+l-n\right)
}h_{2}\left(  x+l\right)  \,dx\\
&  =\sum_{n}z^{n}\int_{-\infty}^{\infty}\overline{h_{1}\left(  x-n\right)
}h_{2}\left(  x\right)  \,dx\\
&  =p_{2}\left(  h_{1},h_{2}\right)  \left(  z\right)  ;
\end{align*}
since we already proved the identity $p_{1}=p_{2}$ (in (\ref{eqZak.16})), the
proposition follows.
\end{proof}

Having established $p_{1}=p_{2}=p_{3}$, we will use $p$ to denote the common
form. Since $p_{3}$ is defined for \emph{all} pairs $h_{i}$ in $L^{2}\left(
\mathbb{R}\right)  $, the compact-support restriction involved in the
formulation of $p_{1}$ and $p_{2}$ has been removed.

\begin{corollary}
\label{CorZak.5}Let $p$ be the form on $L^{2}\left(  \mathbb{R}\right)  \times
L^{2}\left(  \mathbb{R}\right)  $ which is defined in \textup{(\ref{eqZak.17}%
),} and taking values in functions on $\mathbb{T}$. Then in fact $p$ takes
values in $L^{1}\left(  \mathbb{T}\right)  $, i.e., the left-hand side of
\textup{(\ref{eqZak.20})} below is finite if $h_{1},h_{2}\in L^{2}$. For
restricted pairs $h_{1},h_{2}$ of $L^{2}\left(  \mathbb{R}\right)
$-functions, $p\left(  h_{1},h_{2}\right)  $ can also be checked to take
values in $L^{2}\left(  \mathbb{T}\right)  $, and with one more restriction in
$L^{\infty}\left(  \mathbb{T}\right)  $, i.e., the left-hand sides of
\textup{(\ref{eqZak.19})} and \textup{(\ref{eqZak.18}),} below, are finite
with suitable restrictions on $h_{1}$ and $h_{2}$. The restrictions are those
which make the right-hand sides of \textup{(\ref{eqZak.18})} and
\textup{(\ref{eqZak.19})} finite. The respective bounds
\textup{(\ref{eqZak.18}),} \textup{(\ref{eqZak.19}),} and
\textup{(\ref{eqZak.20})} are as follows:%
\begin{align}
\left\|  p\left(  h_{1},h_{2}\right)  \right\|  _{\infty}  &  \leq\mathop{\rm
ess\,sup}_{z}\left\|  Z_{z}h_{1}\right\|  \cdot\left\|  Z_{z}h_{2}\right\|
,\label{eqZak.18}\\
\left\|  p\left(  h_{1},h_{2}\right)  \right\|  _{2}  &  \leq\left\|  p\left(
h_{1},h_{1}\right)  \right\|  _{\infty}^{\frac{1}{2}}\cdot\left\|
h_{2}\right\|  _{2}^{{}},\label{eqZak.19}\\
\left\|  p\left(  h_{1},h_{2}\right)  \right\|  _{1}  &  \leq\left\|
h_{1}\right\|  _{2}\cdot\left\|  h_{2}\right\|  _{2}. \label{eqZak.20}%
\end{align}
\end{corollary}

\begin{proof}
We begin with (\ref{eqZak.20}) since it is universal. Before the estimate, we
may restrict to $h_{i}$ of compact support, and then, after the fact, this
restriction is removed by completion.

In the following estimates, we use the Cauchy--Schwarz inequality two times,
for the respective Hilbert inner products involved:%
\begin{align*}
\left\|  p\left(  h_{1},h_{2}\right)  \right\|  _{1}  &  =\int_{\mathbb{T}%
}\left|  \ip{Z_{z}h_{1}}{Z_{z}h_{2}}\right|  \,d\mu\left(  z\right) \\
&  \leq\int_{\mathbb{T}}\left\|  Z_{z}h_{1}\right\|  \cdot\left\|  Z_{z}%
h_{2}\right\|  \,d\mu\left(  z\right) \\
&  \leq\left(  \int_{\mathbb{T}}\left\|  Z_{z}h_{1}\right\|  ^{2}\,d\mu\left(
z\right)  \cdot\int_{\mathbb{T}}\left\|  Z_{w}h_{2}\right\|  ^{2}\,d\mu\left(
w\right)  \right)  _{\vphantom{\frac12}}^{\frac{1}{2}}\\
&  =\left\|  Zh_{1}\right\|  _{2}\cdot\left\|  Zh_{2}\right\|  _{2}=\left\|
h_{1}\right\|  _{2}\cdot\left\|  h_{2}\right\|  _{2},
\end{align*}
proving (\ref{eqZak.20}).

Formula (\ref{eqZak.18}) for $\left\|  p\left(  h_{1},h_{2}\right)  \right\|
_{\infty}$ is trivial to check, but it is not clear which conditions on the
$h_{i}$'s make the right-hand side finite.

Formula (\ref{eqZak.19}) for $\left\|  p\left(  h_{1},h_{2}\right)  \right\|
_{2}$ is useful later since we will be able to check finiteness of the factor
$\left\|  p\left(  h_{1},h_{1}\right)  \right\|  _{\infty}^{\frac{1}{2}}$. In
fact in many cases, we will end up with $p\left(  h_{1},h_{1}\right)  $ a
constant function on $\mathbb{T}$. Checking (\ref{eqZak.19}) goes as follows:
Let $h_{i}\in L^{2}\left(  \mathbb{R}\right)  $, and suppose $p\left(
h_{1},h_{1}\right)  \in L^{\infty}\left(  \mathbb{T}\right)  $. Then
\begin{align*}
\left\|  p\left(  h_{1},h_{2}\right)  \right\|  _{2}^{2}  &  =\int
_{\mathbb{T}}\left|  \ip{Z_{z}h_{1}}{Z_{z}h_{2}}\right|  ^{2}\,d\mu\left(
z\right) \\
&  \leq\int_{\mathbb{T}}\left\|  Z_{z}h_{1}\right\|  ^{2}\cdot\left\|
Z_{z}h_{2}\right\|  ^{2}\,d\mu\left(  z\right) \\
&  \leq\mathop{\rm ess\,sup}_{z}\left|  p\left(  h_{1},h_{1}\right)  \left(
z\right)  \right|  \cdot\int_{\mathbb{T}}\left\|  Z_{z}h_{2}\right\|
^{2}\,d\mu\left(  z\right) \\
&  =\left\|  p\left(  h_{1},h_{1}\right)  \right\|  _{\infty}^{{}}%
\cdot\left\|  Zh_{2}\right\|  _{2}^{2}\\
&  =\left\|  p\left(  h_{1},h_{1}\right)  \right\|  _{\infty}^{{}}%
\cdot\left\|  h_{2}\right\|  _{2}^{2},
\end{align*}
which is exactly (\ref{eqZak.19}).
\end{proof}

We now use the estimates from Corollary \ref{CorZak.5} to examine boundedness
properties of the representation $\pi$ of $L^{\infty}\left(  \mathbb{T}%
\right)  $ on $L^{2}\left(  \mathbb{R}\right)  $ which was introduced in
(\ref{eqBack.11}), i.e., $\pi\left(  \alpha\right)  h=\alpha\ast h$,
$\alpha\in L^{\infty}\left(  \mathbb{T}\right)  $, $h\in L^{2}\left(
\mathbb{R}\right)  $. When using $L^{2}\left(  \mathbb{R}\right)
\simeq\mathcal{H}_{Z}$ (via the Zak transform), a unitarily equivalent form of
the representation is (also denoted $\pi$) the following:%
\begin{equation}
\pi\left(  \alpha\right)  H\left(  z,x\right)  =\alpha\left(  z\right)
H\left(  z,x\right)  \text{,\qquad for }\alpha\in L^{\infty}\left(
\mathbb{T}\right)  ,\;H\in\mathcal{H}_{Z}. \label{eqZak.21}%
\end{equation}
The idea from the proof of Corollary \ref{CorZak.5} yields imediately%
\[
\left\|  \pi\left(  \alpha\right)  H\right\|  _{\mathcal{H}_{Z}}\leq\left\|
\alpha\right\|  _{\infty}\cdot\left\|  H\right\|  _{\mathcal{H}_{Z}},
\]
which is the bound which is required for a representation. In studying scaling
vectors, however, we shall also need, for fixed $H\in\mathcal{H}_{Z}$,
boundedness properties of the map%
\begin{equation}
C_{H}\colon\alpha\longmapsto\pi\left(  \alpha\right)  H \label{eqZak.22}%
\end{equation}
from $L^{2}\left(  \mathbb{T}\right)  $ to $\mathcal{H}_{Z}$.

\begin{proposition}
\label{ProZak.6}The mapping $C_{H}\colon\alpha\mapsto\pi\left(  \alpha\right)
H$ in \textup{(\ref{eqZak.22})} is bounded from $L^{2}\left(  \mathbb{T}%
\right)  $ to $\mathcal{H}_{Z}$ if and only if $p_{2}\left(  H\right)
=p\left(  H,H\right)  $ is in $L^{\infty}\left(  \mathbb{T}\right)  $, and
then the norm of $C_{H}$ is $\left\|  p_{2}\left(  H\right)  \right\|
_{\infty}^{\frac{1}{2}}$. Moreover, $C_{H}$ has a bounded inverse if and only
if $p_{2}\left(  H\right)  $ has an $L^{\infty}\left(  \mathbb{T}\right)  $
inverse, i.e., there is some $\varepsilon\in\mathbb{R}_{+}$ such that
$p_{2}\left(  H\right)  \left(  z\right)  \geq\varepsilon$ a.e.\ on
$\mathbb{T}$.
\end{proposition}

\begin{proof}
We compute%
\begin{align*}
\left\|  C_{H}\left(  \alpha\right)  \right\|  _{\mathcal{H}_{Z}}^{2}  &
=\int_{\mathbb{T}}\int_{0}^{1}\left|  \alpha\left(  z\right)  H\left(
z,x\right)  \right|  ^{2}\,dx\,d\mu\left(  z\right) \\
&  =\int_{\mathbb{T}}\left|  \alpha\left(  z\right)  \right|  ^{2}\int_{0}%
^{1}\left|  H\left(  z,x\right)  \right|  ^{2}\,dx\,d\mu\left(  z\right) \\
&  =\int_{\mathbb{T}}\left|  \alpha\left(  z\right)  \right|  ^{2}p_{2}\left(
H\right)  \left(  z\right)  \,d\mu\left(  z\right)  .
\end{align*}
If $\alpha\in L^{\infty}\left(  \mathbb{T}\right)  $, then%
\begin{align*}
\left\|  C_{H}\left(  \alpha\right)  \right\|  _{\mathcal{H}_{Z}}^{2}  &
\leq\left\|  \alpha\right\|  _{\infty}^{2}\int_{\mathbb{T}}p_{2}\left(
H\right)  \left(  z\right)  \,d\mu\left(  z\right) \\
&  =\left\|  \alpha\right\|  _{\infty}^{2}\cdot\left\|  H\right\|
_{\mathcal{H}_{Z}}^{2},
\end{align*}
and the assertion follows from a standard fact on multiplication operators.
The same argument also yields the condition for invertibility of $C_{H}$.
\end{proof}

The significance of the operator $C_{H}$, $H\in\mathcal{H}_{Z}$, is that if
$H=F$ is a scaling function,
$\tilde{M}\left(  F\right)  =F$,
where $\tilde{M}=ZMZ^{-1}$.
Then $C_{F}$
intertwines $U$ with a special isometry $S_{0}$ in $L^{2}\left(
\mathbb{T}\right)  \approx\ell^{2}\left(  \mathbb{Z}\right)  $. Let $m_{0}$ be
a filter satisfying (\ref{IntrAxiom(1)})--(\ref{IntrAxiom(3)}) in the
Introduction, and let $M$ be the corresponding cascade operator. Let $F$ be a
scaling function. Let $m_{1}\left(  z\right)  :=z\overline{m_{0}\left(
-z\right)  }$ (which is a high-pass filter), and define $\left(
S_{j}f\right)  \left(  z\right)  =m_{j}\left(  z\right)  f\left(
z^{2}\right)  $, $j=0,1$, $f\in L^{2}\left(  \mathbb{T}\right)  $,
$z\in\mathbb{T}$. Then it is easy to check that%
\[
S_{i}^{\ast}S_{j}^{{}}=\delta_{ij}^{{}}\operatorname*{id}\nolimits_{L^{2}%
\left(  \mathbb{T}\right)  }^{{}}%
\]
and%
\[
\sum_{i=0}^{1}S_{i}^{{}}S_{i}^{\ast}=\operatorname*{id}\nolimits_{L^{2}\left(
\mathbb{T}\right)  }.
\]
We say that the $S_{i}$'s form a representation of the $C^{\ast}$-algebra
(called $\mathcal{O}_{2}$) on the relations. (It was introduced in
\cite{Dix64}; see also \cite{Cun77}.) Such representations were studied
extensively in recent papers \cite{BrJo97b}. A main result from \cite{BrJo97b}
is that
\[
\bigcap_{n=1}^{\infty}S_{0}^{n}L^{2}\left(  \mathbb{T}\right)  =\left\{
0\right\}  .
\]
An isometry $S_{0}$ with this property is called a \emph{shift;} see, e.g.,
\cite{SzFo70}. Let $\mathcal{L}_{0}:=\ker\left(  S_{0}^{\ast}\right)
=S_{1}^{{}}L_{{}}^{2}\left(  \mathbb{T}\right)  $. Then the shift property may
be restated as%
\[
L^{2}\left(  \mathbb{T}\right)  =\sideset{}{^{\smash{\oplus}}}{\sum
}\limits_{n=0}^{\infty}S_{0}^{n}\mathcal{L}_{0}^{{}},
\]
and moreover the spaces $S_{0}^{n}\mathcal{L}_{0}^{{}}$ are pairwise mutually
orthogonal.\renewcommand{\theenumi}{\roman{enumi}}

\begin{theorem}
\label{ThmZak.7}Let the setting be as above, and consider the cascade problem
in $\mathcal{H}_{Z}\simeq L^{2}\left(  \mathbb{R}\right)  $. Then the
following two conditions are equivalent:

\begin{enumerate}
\item \label{ThmZak.7(1)}$\tilde{M}\left(  F\right)  =F$, and

\item \label{ThmZak.7(2)}$C_{F}$ intertwines $S_{0}$ and $\tilde{U}$, i.e.,
$C_{F}S_{0}=\tilde{U}C_{F}$.
\end{enumerate}

\noindent Let $V_{0}\left(  F\right)  =\left[  \left\{  \pi\left(
\alpha\right)  F\mid\alpha\in L^{2}\left(  \mathbb{T}\right)  \right\}
\right]  $ where $\left[  \,\cdot\,\right]  $ is norm-closure. Let
$V_{n}\left(  F\right)  :=\tilde{U}^{n}\left(  V_{0}\left(  F\right)  \right)
$. If \textup{(\ref{ThmZak.7(1)})} holds, then $V_{n+1}\left(  F\right)
\subset V_{n}\left(  F\right)  $, and $C_{F}^{{}}\left(  S_{0}^{n}L_{{}}%
^{2}\left(  \mathbb{T}\right)  \right)  =V_{n}\left(  F\right)  $. Let
$W_{n}\left(  F\right)  :=V_{n}\left(  F\right)  \ominus V_{n+1}\left(
F\right)  $. If further $p_{2}\left(  F\right)  \equiv1$, then $C_{F}^{{}%
}\left(  S_{0}^{n}\mathcal{L}_{0}^{{}}\right)  =W_{n}^{{}}\left(  F\right)  $.
\end{theorem}

The results of the theorem may be summarized as in Table \ref{ThmZak.7(table)} above.

\begin{table}[ptb]
\caption{Summary of Theorem \ref{ThmZak.7}: Embedding of the isometric model
into $L^{2}\left(  \mathbb{R}\right)  $}%
\label{ThmZak.7(table)}%
\renewcommand{\arraystretch}{0.25}
\begin{tabular}
[c]{ccccrrrcl}%
$\left\{  0\right\}  $ & $\longleftarrow$ & $\cdots$ & $\longleftarrow$ &
$V_{2}\left(  F\right)  \raisebox{-12pt}{$\searrow\vphantom{\raisebox
{-2pt}{$\searrow$}}$\hskip-8pt}$ & $V_{1}\left(  F\right)  \raisebox
{-12pt}{$\searrow\vphantom{\raisebox{-2pt}{$\searrow$}}$\hskip-8pt}$ &
$V_{0}\left(  F\right)  \raisebox{-12pt}{$\searrow\vphantom{\raisebox
{-2pt}{$\searrow$}}$\hskip-8pt}$ &  & finer scales\\\hline
&  &  &  & \multicolumn{1}{c}{} & \multicolumn{1}{c}{} & \multicolumn{1}{c}{}%
&  & \\\cline{1-7}%
&  &  &  & \multicolumn{1}{c}{} & \multicolumn{1}{c}{} & \multicolumn{1}{c}{}%
& \multicolumn{1}{|c}{} & \\\cline{1-6}%
&  &  &  & \multicolumn{1}{c}{} & \multicolumn{1}{c}{} & \multicolumn{1}{|c}{}%
& \multicolumn{1}{|c}{} & \\\cline{1-5}%
&  &  &  & \multicolumn{1}{c}{} & \multicolumn{1}{|c}{} &
\multicolumn{1}{|c}{} & \multicolumn{1}{|c}{} & \\\cline{1-4}%
&  & $\cdots$ &  & \multicolumn{1}{|c}{$W_{2}\left(  F\right)  $} &
\multicolumn{1}{|c}{$W_{1}\left(  F\right)  $} & \multicolumn{1}{|c}{$W_{0}%
\left(  F\right)  $} & \multicolumn{1}{|c}{$\cdots$} & rest of $L^{2}\left(
\mathbb{R}\right)  $\\\cline{1-4}%
&  &  &  & \multicolumn{1}{c}{} & \multicolumn{1}{|c}{} &
\multicolumn{1}{|c}{} & \multicolumn{1}{|c}{} & \\\cline{1-5}%
&  &  &  & \multicolumn{1}{c}{} & \multicolumn{1}{c}{} & \multicolumn{1}{|c}{}%
& \multicolumn{1}{|c}{} & \\\cline{1-6}%
&  &  &  & \multicolumn{1}{c}{} & \multicolumn{1}{c}{} & \multicolumn{1}{c}{}%
& \multicolumn{1}{|c}{} & \\\cline{1-7}%
&  &  &  & \multicolumn{1}{c}{} & \multicolumn{1}{c}{} & \multicolumn{1}{c}{}%
&  & \\\hline
&  & $\cdots$ & \multicolumn{1}{r}{$\rlap{$\underset{\textstyle\vphantom
{U}\smash{\tilde{U}}}{\longleftarrow}$}\hskip3pt$} & $\rlap{$\underset
{\textstyle\vphantom{U}\smash{\tilde{U}}}{\longleftarrow}$}\hskip3pt$ &
$\rlap{$\underset{\textstyle\vphantom{U}\smash{\tilde{U}}}{\longleftarrow}%
$}\hskip3pt$ & \multicolumn{1}{c}{} & $\raisebox{-7pt}[10pt]{$L^{2}%
\left( \mathbb{R}\right) \rlap{${}\approx\mathcal{H}_{Z}$}$}$ & \\
&  & $\llap{$C_{F}$}\makebox[6pt]{\raisebox{4pt}{\makebox[0pt]{\hss
$\uparrow$\hss}}\raisebox{-4pt}{\makebox[0pt]{\hss$|$\hss}}}$ &  &
\multicolumn{1}{c}{} & \multicolumn{1}{c}{} & \multicolumn{1}{c}{} &
\llap{\hooklonguparrow\vphantom{\raisebox{2pt}{\hooklonguparrow}}%
\vphantom{\raisebox{-2pt}{\hooklonguparrow}}}\hskip5pt\rlap{$C_{F}\colon
\alpha\mapsto\pi
\left( \alpha\right) F$} & \\
\rule[-6pt]{0pt}{6pt}$\left\{  0\right\}  $ & $\longleftarrow$ & $\cdots$ &
\multicolumn{1}{r}{$\rlap{$\overset{\textstyle
S_{0}}{\longleftarrow}$}\hskip3pt$} & $\rlap{$\overset{\textstyle S_{0}%
}{\longleftarrow}$}\hskip3pt$ & $\rlap{$\overset{\textstyle S_{0}%
}{\longleftarrow}$}\hskip3pt$ & \multicolumn{1}{c}{} & $\raisebox{5pt}%
{$L^{2}\left( \mathbb{T}\right) $}$ & \\\cline{1-7}%
&  &  &  & \multicolumn{1}{c}{} & \multicolumn{1}{c}{} & \multicolumn{1}{c}{}%
& \multicolumn{1}{|c}{} & \\\cline{1-6}%
&  &  &  & \multicolumn{1}{c}{} & \multicolumn{1}{c}{} & \multicolumn{1}{|c}{}%
& \multicolumn{1}{|c}{} & \\\cline{1-5}%
&  &  &  & \multicolumn{1}{c}{} & \multicolumn{1}{|c}{} &
\multicolumn{1}{|c}{} & \multicolumn{1}{|c}{} & \\\cline{1-4}%
&  & $\cdots$ &  & \multicolumn{1}{|c}{$S_{0}^{2}\mathcal{L}_{0}^{{}}$} &
\multicolumn{1}{|c}{$S_{0}\mathcal{L}_{0}$} & \multicolumn{1}{|c}{$\mathcal{L}%
_{0}=S_{1}L^{2}\left(  \mathbb{T}\right)  $} & \multicolumn{1}{|c}{} &
\\\cline{1-4}%
&  &  &  & \multicolumn{1}{c}{} & \multicolumn{1}{|c}{} &
\multicolumn{1}{|c}{} & \multicolumn{1}{|c}{} & \\\cline{1-5}%
&  &  &  & \multicolumn{1}{c}{} & \multicolumn{1}{c}{} & \multicolumn{1}{|c}{}%
& \multicolumn{1}{|c}{} & \\\cline{1-6}%
&  &  &  & \multicolumn{1}{c}{} & \multicolumn{1}{c}{} & \multicolumn{1}{c}{}%
& \multicolumn{1}{|c}{} & \\\cline{1-7}%
&  &  &  & \llap{$S_{0}^{2}L^{2}\left(  \mathbb{T}\right) $}\raisebox
{12pt}{$\nearrow\vphantom{\raisebox{2pt}{$\nearrow$}}$\hskip-8pt} &
\llap{$S_{0}L^{2}\left(  \mathbb{T}\right) $}\raisebox{12pt}{$\nearrow
\vphantom{\raisebox{2pt}{$\nearrow$}}$\hskip-8pt} & $L^{2}\left(
\mathbb{T}\right)  \raisebox{12pt}{$\nearrow\vphantom{\raisebox{2pt}%
{$\nearrow$}}$\hskip-8pt}$ &  &
\end{tabular}
\end{table}

\begin{proof}
(\ref{ThmZak.7(1)}) $\Rightarrow$ (\ref{ThmZak.7(2)}). Suppose $\tilde{U}%
^{-1}\pi\left(  m_{0}\right)  F=F$; then%
\begin{align*}
\tilde{U}C_{F}\left(  \alpha\right)   &  =\tilde{U}\pi\left(  \alpha\right)
F=\pi\left(  \alpha\left(  z^{2}\right)  \right)  \tilde{U}F=\pi\left(
\alpha\left(  z^{2}\right)  \right)  \pi\left(  m_{0}\right)  F\\
&  =\pi\left(  m_{0}\left(  z\right)  \alpha\left(  z^{2}\right)  \right)
F=\pi\left(  S_{0}\alpha\right)  F=C_{F}S_{0}\left(  \alpha\right)
\text{\qquad for all }\alpha\in L^{2}\left(  \mathbb{T}\right)  ,
\end{align*}
which is (\ref{ThmZak.7(2)}).

(\ref{ThmZak.7(2)}) $\Rightarrow$ (\ref{ThmZak.7(1)}). Suppose
(\ref{ThmZak.7(2)}) holds. Then the previous calculation reverses, and shows
that $\tilde{M}\left(  F\right)  =F$ where $\tilde{M}=\tilde{U}^{-1}\pi\left(
m_{0}\right)  $, and we have (\ref{ThmZak.7(1)}).

Suppose (\ref{ThmZak.7(1)}), and let $\alpha\in L^{2}\left(  \mathbb{T}%
\right)  $; then%
\[
C_{F}^{{}}S_{0}^{n}\alpha=\tilde{U}_{{}}^{n}C_{F}^{{}}\left(  \alpha\right)
=\tilde{U}_{{}}^{n}\pi\left(  \alpha\right)  F\in V_{n}^{{}}\left(  F\right)
,
\]
and we conclude that $C_{F}^{{}}\left(  S_{0}^{n}L^{2}\left(  \mathbb{T}%
\right)  \right)  =V_{n}^{{}}\left(  F\right)  $ as claimed.

Now let $\alpha,\beta\in L^{2}\left(  \mathbb{T}\right)  $, and consider the
following calculation:%
\begin{align*}
\ip{C_{F}S_{1}\alpha}{\smash{\tilde{U}}\pi\left( \beta\right) F}  &
=\ip{C_{F}S_{1}\alpha}{C_{F}S_{0}\alpha}\\
&  =\int_{\mathbb{T}}\overline{m_{1}\left(  z\right)  }m_{0}\left(  z\right)
\overline{\alpha\left(  z^{2}\right)  }\beta\left(  z^{2}\right)  p_{2}\left(
F\right)  \left(  z\right)  \,d\mu\left(  z\right)  .
\end{align*}
If $p_{2}\left(  F\right)  \equiv1$, then this integral is%
\[
\int_{\mathbb{T}}\frac{1}{2}\sum_{w^{2}=z}\overline{m_{1}\left(  w\right)
}m_{0}\left(  w\right)  \overline{\alpha\left(  z\right)  }\beta\left(
z\right)  \,d\mu\left(  z\right)  .
\]
But the choice of $m_{1}\left(  z\right)  =z\overline{m_{0}\left(  -z\right)
}$ makes it zero:%
\[
\frac{1}{2}\sum_{w^{2}=z}\overline{m_{1}\left(  w\right)  }m_{0}\left(
w\right)  =\frac{1}{2}\sum_{w^{2}=z}\bar{w}m_{0}\left(  -w\right)
m_{0}\left(  w\right)  =0,
\]
since $\sum_{w^{2}=z}w=0$.

We have proved that if $p_{2}\left(  F\right)  =\openone$, $C_{F}$ maps
$\mathcal{L}_{0}=S_{1}\left(  L^{2}\left(  \mathbb{T}\right)  \right)  $ onto
$W_{0}\left(  F\right)  =V_{0}\left(  F\right)  \ominus\tilde{U}V_{0}\left(
F\right)  $, and an iteration of the same argument (induction) yields%
\[
C_{F}^{{}}\left(  S_{0}^{n}S_{1}^{{}}L_{{}}^{2}\left(  \mathbb{T}\right)
\right)  =W_{n}\left(  F\right)
\]
as claimed.

We also saw that $C_{F}$ is \emph{isometric} if and only if $p_{2}\left(
F\right)  \equiv1$.
\end{proof}

The significance of the spaces $V_{0}\left(  F\right)  $ and $W_{0}\left(
F\right)  $ in wavelet theory is that, in the $L^{2}\left(  \mathbb{R}\right)
$-picture, $V_{0}\left(  F\right)  $, $F=Z\varphi$, is generated by the father
function, while $W_{1}\left(  F\right)  $ is generated by the mother function.
It is interesting to summarize the approximation properties of the two
function sequences $m_{0}^{\left(  n\right)  }\left(  z\right)  =m_{0}^{{}%
}\left(  z\right)  m_{0}^{{}}\left(  z^{2}\right)  \cdots m_{0}^{{}}\left(
z^{2^{n-1}}\right)  _{{}}$ and $D_{n}^{{}}\left(  z\right)  =\left|
m_{0}^{\left(  n\right)  }\left(  z\right)  \right|  ^{2}$. In Table
\ref{tableZak.2} we include results from \cite{BrJo97b} in the left-hand
column, and results from Meyer and Paiva \cite{MePa93} in the right-hand
column. The filter $m_{0}$ is given as usual, and (\ref{IntrAxiom(1)}%
)--(\ref{IntrAxiom(3)}) in the Introduction are assumed.

\begin{table}[ptb]
\caption{Approximation properties of $m_{0}^{\left(  n\right)  }$ and
$D_{n}^{{}}\left(  z\right)  $}%
\label{tableZak.2}%
\renewcommand{\arraystretch}{2.75}
\begin{tabular}
[c]{p{186pt}|p{186pt}}%
\makebox[186pt]{$m_{0}^{\left( n\right) }$} & \makebox[186pt]{$D_{n}^{{}}%
\left( z\right)
=\left| m_{0}^{\left( n\right) }\left( z\right) \right| ^{2}_{\mathstrut}%
$}\\\hline
\parbox[t]{186pt}{$m_{0}^{\left( n\right) }=S_{0}^{n}\left( \openone
\right) $} & \parbox[t]{186pt}{$D_{n}=R^{\ast\,n}\left
( \openone\right) $}\\
\parbox[t]{186pt}{$\int_{\mathbb{T}}m_{0}^{\left( n\right) }f\,d\mu
\underset{n\rightarrow\infty}{\longrightarrow}0 \quad\forall f\in L^{2}%
\left( \mathbb{T}\right) $} & \parbox[t]{186pt}{$\int_{\mathbb
{T}}D_{n}f\,d\mu\underset{n\rightarrow\infty}{\longrightarrow}f\left
( 1\right) \quad\forall f\in C\left( \mathbb{T}\right) $ \\if $ \left\{ f\mid
Rf=f\right\} $ is one-dimensional \cite{MePa93}}\\
\parbox[t]{186pt}{$\int_{\mathbb{T}}m_{0}^{\left( n\right) }\,d\mu=\left
( a_{0}\right) ^{n} $ \\if $ m_{0}\left( z\right) =\sum_{k=0}^{\infty}%
a_{k}z^{k}$} & \parbox[t]{186pt}{$\int_{\mathbb{T}}D_{n}\left( z\right
) \,d\mu\left( z\right) =1 \quad\forall n=1,2,\dots
$}%
\end{tabular}
\end{table}

We now show that when $m_{0}$ is given, and $M$ is the corresponding cascade
operator, \emph{then the solutions }$M\left(  F\right)  =F$ \emph{may be
identified with a space of intertwining operators.} We shall state the details
in the Hilbert space $\mathcal{H}_{Z}$ of the Zak transform.

\begin{corollary}
\label{CorZak.8}An operator $C\colon L^{2}\left(  \mathbb{T}\right)
\rightarrow\mathcal{H}_{Z}$ is of the form $C_{F}$ for some $F\in
\mathcal{H}_{Z}$ satisfying $\tilde{M}\left(  F\right)  =F$ if and only if it
satisfies the following two intertwining properties:\renewcommand{\theenumi
}{\alph{enumi}}

\begin{enumerate}
\item \label{CorZak.8(1)}$CS_{0}=\tilde{U}C$, and

\item \label{CorZak.8(2)}$C\left(  f\alpha\right)  =\left(  \pi\left(
f\right)  C\right)  \left(  \alpha\right)  $, for all $f\in L^{\infty}\left(
\mathbb{T}\right)  $ and $\alpha\in L^{2}\left(  \mathbb{T}\right)  $.
\end{enumerate}

\noindent If multiplication by $L^{\infty}\left(  \mathbb{T}\right)  $ on
$L^{2}\left(  \mathbb{T}\right)  $ is written $\tau\left(  f\right)
\alpha=f\alpha$, then \textup{(\ref{CorZak.8(2)})} reads:

\begin{enumerate}
\item [(\ref{CorZak.8(2)}$^{\prime}$)]$C\tau\left(  f\right)  =\pi\left(
f\right)  C$.
\end{enumerate}
\end{corollary}

\begin{proof}
By Theorem \ref{ThmZak.7}, it is enough to check that every operator $C\colon
L^{2}\left(  \mathbb{T}\right)  \rightarrow\mathcal{H}_{Z}$ which satisfies
(\ref{CorZak.8(2)}) $\simeq$ (\ref{CorZak.8(2)}$^{\prime}$) must be of the
form $C=C_{F}$ for some $F\in\mathcal{H}_{Z}$. So let $C$ be given, and assume
(\ref{CorZak.8(2)}). Let $C\left(  \openone\right)  :=F$. Then, for $\alpha\in
L^{\infty}\left(  \mathbb{T}\right)  $, we have $C\left(  \alpha\right)
=C\tau\left(  \alpha\right)  \openone=\pi\left(  \alpha\right)  C\openone
=\pi\left(  \alpha\right)  F=C_{F}\left(  \alpha\right)  $. We are considering
only the case when $C$ is bounded. Since $L^{\infty}\left(  \mathbb{T}\right)
$ is dense in $L^{2}\left(  \mathbb{T}\right)  $, the result follows. If
(\ref{CorZak.8(1)}) also holds, we saw in Theorem \ref{ThmZak.7} that then
$\tilde{M}\left(  F\right)  =F$, and the proof is completed.
\end{proof}

\begin{remark}
\label{RemOnCorZak.8} \textbf{An inner product on the intertwining operators.}
From \textup{(\ref{CorZak.8(2)}),} we also get the identities%
\[
C_{H}^{\ast}C_{H}^{{}}=\tau\left(  p_{2}\left(  H\right)  \right)
\text{\qquad and\qquad}C_{H}^{\ast}C_{H^{\prime}}^{{}}=\tau\left(  p\left(
H,H^{\prime}\right)  \right)
\]
for $H,H^{\prime}\in\mathcal{H}_{Z}$, where%
\[
p_{2}\left(  H\right)  \left(  z\right)  =p\left(  H,H\right)  \left(
z\right)  =\int_{0}^{1}\left|  H\left(  z,x\right)  \right|  ^{2}\,dx,
\]
and $\tau\left(  f\right)  $ denotes multiplication by $f$ on $L^{2}\left(
\mathbb{T}\right)  $. Hence \textup{(}as noted\textup{),} $\left\|
C_{H}\right\|  =\left\|  p_{2}\left(  H\right)  \right\|  _{\infty}^{\frac
{1}{2}}$, and $C_{H}$ is bounded if and only if $p_{2}\left(  H\right)  \in
L^{\infty}\left(  \mathbb{T}\right)  $.
\end{remark}

\begin{remark}
\label{RemZak.9} Proposition \textup{\ref{ProZak.6}} will be used in Sections
\textup{\ref{SiC}--}\ref{SiV} in the study of scaling functions $\varphi\in
L^{2}\left(  \mathbb{R}\right)  $, i.e., solutions to $M\varphi=\varphi$ where
$M$ is a cascade operator for some given filter $m_{0}$. If $Z\varphi=F$, then
the scaling equation is equivalent to $\tilde{M}\left(  F\right)  =F$. Then,
generically, $C_{F}$ will be bounded, but will not have a bounded inverse.

If $m_{0}\left(  z\right)  =\frac{1}{\sqrt{2}\mathstrut}\left(  1+z^{3}%
\right)  $, then $\varphi\left(  x\right)  =\frac{1}{3}\chi^{{}}_{\left[
0,3\right]  }\left(  x\right)  $, and for $0\leq x\leq1$,
\[
F\left(  z,x\right)  =Z\varphi\left(  z,x\right)  =\frac{1}{3}\left(  \chi
^{{}}_{\left[  0,3\right]  }\left(  x\right)  +z\chi^{{}}_{\left[
-1,2\right]  }\left(  x\right)  +z^{2}\chi^{{}}_{\left[  -2,1\right]  }\left(
x\right)  \right)  .
\]
It is easy then to check that%
\begin{align*}
p_{2}\left(  H\right)  \left(  z\right)  =p\left(  H,H\right)  \left(
z\right)   &  =\frac{1}{9}\left(  z^{-2}+2z^{-1}+3+2z+z^{2}\right) \\
&  =\frac{1}{9}\left(  3+4\cos\omega+2\cos\left(  2\omega\right)  \right)  ,
\end{align*}
where $z=e^{-i\omega}$. In this case, $C_{F}$ is not invertible as
$p_{2}\left(  H\right)  $ vanishes on $\mathbb{T}$. In fact, $p_{2}\left(
H\right)  \left(  \omega\right)  =\frac{1}{9}\left(  2\cos\omega+1\right)
^{2}$ which, of course, vanishes for $\omega=\pm\frac{2\pi}{3}$.
\end{remark}

\begin{remark}
\label{RemNewZak.11}Using the usual isomorphism $L^{2}\left(  \mathbb{T}%
\right)  \simeq\ell^{2}\left(  \mathbb{Z}\right)  $ defined by the
Plan\-che\-rel theorem for Fourier series, we note that, if $m\left(
z\right)  =\sum_{n\in\mathbb{Z}}a_{n}z^{n}$, then the operator%
\[
f\longmapsto m\left(  z\right)  f\left(  z^{2}\right)  \text{,\qquad on }%
L^{2}\left(  \mathbb{T}\right)  ,
\]
takes the form%
\begin{equation}
\left(  S\xi\right)  _{n}=\sum_{k\in\mathbb{Z}}a_{n-2k}\xi_{k}
\label{eqNewZak.24}%
\end{equation}
when realized as an operator on the sequence space $\ell^{2}\left(
\mathbb{Z}\right)  $ via the Fourier series representation%
\[
f\left(  z\right)  =\sum_{n\in\mathbb{Z}}\xi_{n}z^{n},\qquad\left(  \xi
_{n}\right)  \in\ell^{2},\qquad\sum_{n\in\mathbb{Z}}\left|  \xi_{n}\right|
^{2}=\left\|  f\right\|  _{2}^{2},
\]
and this is the connection to the Micchelli operator \textup{(\ref{eqIntr.7})}
mentioned in the Introduction. We sketch the details of this argument below,
and refer to \cite{Mic96} for more details.

Starting with $2\pi$-periodic functions $m$ and $f$, corresponding to the
Fourier representation%
\begin{align*}
m\left(  z\right)   &  \sim m\left(  \omega\right)  \sim\sum_{n\in\mathbb{Z}%
}a_{n}z^{n},\\%
\intertext{and}%
f\left(  z\right)   &  \sim f\left(  \omega\right)  \sim\sum_{n\in\mathbb{Z}%
}\xi_{n}z^{n},
\end{align*}
with $z=e^{-i\omega}$, $\omega\in\mathbb{R}$, as the $\mathbb{T}%
=\mathbb{R}\diagup2\pi\mathbb{Z}$ convention, we have Micchelli's operator $S$
of \textup{(\ref{eqNewZak.24})} or \textup{(\ref{eqIntr.7})} in the function
form%
\[
\left(  \tilde{S}f\right)  \left(  \omega\right)  =m\left(  \omega\right)
f\left(  2\omega\right)  ,\qquad\omega\in\mathbb{R}.
\]
For each $k$, the iteration $\tilde{S}^{k}f$ is also $2\pi$-periodic, while
$\left(  \tilde{S}^{k}f\right)  \left(  \frac{\omega}{2^{k\mathstrut}}\right)
$ has period $2^{k}\cdot\left(  2\pi\right)  $. Introducing $\left(
UF\right)  \left(  x\right)  =2^{-\frac{1}{2}}F\left(  \frac{x}{2}\right)  $,
$x\in\mathbb{R}$, on functions, or distributions, on $\mathbb{R}$, we arrive
at the representation%
\begin{equation}
U^{k}\tilde{S}^{k}f\sim2^{-\frac{k}{2}}m\left(  \omega\right)  m\left(
\frac{\omega}{2}\right)  \cdots m\left(  \frac{\omega}{2^{k\mathstrut}%
}\right)  f\left(  \omega\right)  . \label{eqNewZak.25}%
\end{equation}
If there is a limit function \textup{(}or distribution\textup{)} $F_{\xi}$, as
$k\rightarrow\infty$, then the difference%
\[
\left(  \Delta^{k}\xi\right)  _{j}=\left(  S^{k}\xi\right)  _{j}-F_{\xi
}\left(  \frac{j}{2^{k\mathstrut}}\right)  ,\qquad j\in\mathbb{Z},
\]
tends to zero in the limit $k\rightarrow\infty$. Note that $\Delta^{k}$, for
each $k$, is acting on sequences, say $\ell^{2}\left(  \mathbb{Z}\right)  $.
Hence we get the corresponding scaling function $F_{\xi}$ from
\textup{(\ref{eqNewZak.25})} at the dyadic rational points $\left\{  \frac
{j}{2^{k\mathstrut}}\bigm|j,k\in\mathbb{Z}\right\}  \subset\mathbb{R}$ this
way, and we have therefore made the connection to the Micchelli approximation
of \cite{Mic96}; see also \textup{(\ref{eqIntr.7})} in the Introduction above.
\end{remark}

\section[Proof of Theorem \textup{\ref{ThmBack.1}}]{\label{Poof}PROOF OF THEOREM \textup{\ref{ThmBack.1}}}

Recall that $Z$ is an isometric isomorphism, viz.:%
\[
L^{2}\left(  \mathbb{R}\right)  \simeq\mathcal{H}_{Z}\simeq L^{2}\left(
\mathbb{T}\times\left[  0,1\right]  \right)  ,
\]
where the second $\simeq$ amounts to restriction from $\mathbb{R}$ to $\left[
0,1\right]  $ in the $x$-variable: if $H\in\mathcal{H}_{Z}$ satisfies
(\ref{eqZak.2}), then the restriction $H\left(  z,x\right)  $, $0\leq x\leq1$,
defines the corresponding element in $L^{2}\left(  \mathbb{T}\times\left[
0,1\right]  \right)  $, and a simple argument shows that this restriction
mapping is indeed an isomorphic \emph{isometry} of $\mathcal{H}_{Z}$ onto
$L^{2}\left(  \mathbb{T}\times\left[  0,1\right]  \right)  $. It follows that
operators in one space identify with corresponding operators in the other. If
$A$ is a given operator in $L^{2}\left(  \mathbb{R}\right)  $, then $\tilde
{A}:=ZAZ^{\ast}$ is the corresponding operator in $\mathcal{H}_{Z}$.

The proof of the following lemma is essentially contained in the previous
section: see especially (\ref{eqZak.10}) and Lemma \ref{LemZak.3}.

\begin{lemma}
\label{LemPoof.1}If $A$ is one of the operators in $L^{2}\left(
\mathbb{R}\right)  $ listed in the first column of Table
\textup{\ref{LemPoof.1(table)}}, then $\tilde{A}=ZAZ^{\ast}$ in $\mathcal{H}%
_{Z}$ is given by the corresponding entry in the second column of Table
\textup{\ref{LemPoof.1(table)}.}
\end{lemma}

\begin{table}[ptb]
\caption{Operator correspondence between $L^{2}\left(  \mathbb{R}\right)  $
and $\mathcal{H}_{Z}$ (Lemma \ref{LemPoof.1})}%
\label{LemPoof.1(table)}%
\renewcommand{\arraystretch}{2.75}
\begin{tabular}
[c]{p{186pt}|p{186pt}}%
\makebox[186pt]{$h\in L^{2}\left( \mathbb{R}\right) $} & \makebox
[186pt]{$H\in\mathcal{H}_{Z\mathstrut}$}\\\hline
\parbox[t]{186pt}{$\left( T_{n}h\right) \left( x\right) =h\left( x+n\right
) $} & \parbox[t]{186pt}{$\tilde{T}_{n}H\left( z,x\right
) =z^{-n}H\left( z,x\right) $}\\
\parbox[t]{186pt}{$\pi\left( \alpha\right) h=\alpha\ast h$} & \parbox
[t]{186pt}{$\widetilde{\pi\left( \alpha\right) }H\left( z,x\right
) =\alpha\left( z\right) H\left( z,x\right) $}\\
\parbox[t]{186pt}{$M=U^{-1}\pi\left( m_{0}\right) $} & \parbox[t]{186pt}%
{$\tilde{M}H\left( z,x\right) =$ \\$\frac{1}{\sqrt{2}\mathstrut
}\sum_{w^{2}=z}m_{0}\left( w\right) H\left( w,2x\right) $}\\
\parbox[t]{186pt}{$M^{\ast}=\pi\left( \bar{m}_{0}\right) U$} & \parbox
[t]{186pt}{$\tilde{M}^{\ast}H\left( z,x\right) =$ \\$\frac{1}{\sqrt
{2}\mathstrut}\overline{m_{0}\left( z\right) }\left( H\left( z^{2},\frac{x}%
{2}\right) +zH\left( z^{2},\frac{x+1}{2}\right) \right) $}\\
\parbox[t]{186pt}{$\left( E_{t}h\right) \left( x\right) =e^{itx}h\left
( x\right) $} & \parbox[t]{186pt}{$\tilde{E}_{t}H\left
( z,x\right) =e^{itx}H\left( ze^{it},x\right) $}\\
\parbox[t]{186pt}{$\left( \mathcal{F}h\right) \left( x\right) =\int_{\mathbb
{R}}e^{-i2\pi xy}h\left( y\right) \,dy$} & \parbox[t]{186pt}{$\left
( \widetilde{\mathcal{F}}H\right) \left( e^{i2\pi\omega},x\right
) =$ \\$e^{-i2\pi x\omega}H\left( e^{-i2\pi x},\omega\right) $ for $\omega
,x\in\mathbb{R}$}%
\end{tabular}
\end{table}

\begin{proof}
In the previous section, we also elaborated on the operators $\pi\left(
\alpha\right)  $, $\alpha\in L^{\infty}\left(  \mathbb{T}\right)  $,
$Uh\left(  x\right)  =\frac{1}{\sqrt{2}\mathstrut}h\left(  \frac{x}{2}\right)
$, and the cascade operator
\[
Mh\left(  x\right)  =\sqrt{2}\sum_{n\in\mathbb{Z}}a_{n}h\left(  2x-n\right)
,
\]
for $m_{0}\left(  z\right)  =\sum_{n\in\mathbb{Z}}a_{n}z^{n}$, representing
the given low-pass filter. The present proof amounts to a combination of the
calculations leading up to Lemma \ref{LemZak.3}, and the argument from the
proof of that lemma.
\end{proof}

\begin{proof}
[Proof of Theorem \textup{\ref{ThmBack.1}}]With the aid of Lemma
\ref{LemPoof.1}, the proof of the two commutation relations
(\ref{ThmBack.1(1)})--(\ref{ThmBack.1(2)}) in Theorem \ref{ThmBack.1} now
amounts to the following computations. They take place in the space
$\mathcal{H}_{Z}$, i.e., the range of the Zak transform, so it is the
right-hand column in Table \ref{LemPoof.1(table)} which is used.

\emph{Ad} (\ref{ThmBack.1(1)}): Let $\alpha\in L^{\infty}\left(
\mathbb{T}\right)  $ and $H\in\mathcal{H}_{Z}$. Then%
\begin{align*}
\left(  \tilde{M}^{\ast}\widetilde{\pi\left(  \alpha\right)  }\tilde
{M}\right)  H\left(  z,x\right)   &  =\frac{1}{\sqrt{2}\mathstrut}%
\,\overline{m_{0}\left(  z\right)  }\,\alpha\left(  z^{2}\right)  \left(
\tilde{M}H\left(  z^{2},\frac{x}{2}\right)  +z\tilde{M}H\left(  z^{2}%
,\frac{x+1}{2}\right)  \right) \\
&  =\frac{1}{2}\,\overline{m_{0}\left(  z\right)  }\, \alpha\left(
z^{2}\right)  \sum_{w^{2}=z^{2}}m_{0}\left(  w\right)  \left(  H\left(
w,x\right)  +zH\left(  w,x+1\right)  \right) \\
&  =\frac{1}{2}\,\overline{m_{0}\left(  z\right)  }\, \alpha\left(
z^{2}\right)  \sum_{w^{2}=z^{2}}m_{0}\left(  w\right)  \left(  H\left(
w,x\right)  +zw^{-1}H\left(  w,x\right)  \right) \\
&  =\frac{1}{2}\,\overline{m_{0}\left(  z\right)  }\, \alpha\left(
z^{2}\right)  \sum_{w^{2}=z^{2}}m_{0}\left(  w\right)  \left(  1+zw^{-1}%
\right)  H\left(  w,x\right)  .
\end{align*}
But the summation is over $w\in\left\{  \pm z\right\}  $, and the term
$1+zw^{-1}=2$ if $w=z$, and $1+zw^{-1}=0$ if $w=-z$. We get%
\[
\left(  \tilde{M}^{\ast}\widetilde{\pi\left(  \alpha\right)  }\tilde
{M}\right)  H\left(  z,x\right)  =\left|  m_{0}\left(  z\right)  \right|
^{2}\alpha\left(  z^{2}\right)  H\left(  z,x\right)  ,
\]
which is the identity (\ref{ThmBack.1(1)}) of Theorem \ref{ThmBack.1}.

\emph{Ad} (\ref{ThmBack.1(2)}): As in (\ref{ThmBack.1(1)}), let $\alpha\in
L^{\infty}\left(  \mathbb{T}\right)  $ and $H\in\mathcal{H}_{Z}$ be given.
Then%
\begin{align*}
\left(  \tilde{M}\widetilde{\pi\left(  \alpha\right)  }\tilde{M}^{\ast
}\right)  H\left(  z,x\right)   &  =\frac{1}{\sqrt{2}\mathstrut}\sum_{w^{2}%
=z}m_{0}\left(  w\right)  \alpha\left(  w\right)  \left(  \tilde{M}^{\ast
}H\right)  \left(  w,2x\right) \\
&  =\frac{1}{2}\sum_{w^{2}=z}\left|  m_{0}\left(  w\right)  \right|  ^{2}%
\cdot\alpha\left(  w\right)  \cdot\left(  H\left(  z,x\right)  +wH\left(
z,x+\frac{1}{2}\right)  \right) \\
&  =R\left(  \alpha\right)  \left(  z\right)  H\left(  z,x\right)  +R\left(
e_{1}\alpha\right)  \left(  z\right)  H\left(  z,x+\frac{1}{2}\right)  ,
\end{align*}
which is precisely the second identity (\ref{ThmBack.1(2)}) of Theorem
\ref{ThmBack.1}. Note that we obtain the identities in $\mathcal{H}_{Z}$, but
since $Z$ is an isomorphism, $Z\colon L^{2}\left(  \mathbb{R}\right)
\underset{\simeq}{\longrightarrow}\mathcal{H}_{Z}$, we automatically get the
same identities in $L^{2}\left(  \mathbb{R}\right)  $ where the translation to
$L^{2}\left(  \mathbb{R}\right)  $ is made via the dictionary of Lemma
\ref{LemPoof.1} (Table \ref{LemPoof.1(table)}).
\end{proof}

\begin{remark}
\label{RemPoof.2}The last line in the dictionary of Lemma
\textup{\ref{LemPoof.1}} \textup{(}Table \textup{\ref{LemPoof.1(table)})} is
the correspondence for the Fourier transform $\mathcal{F}$ in $L^{2}\left(
\mathbb{R}\right)  $, and it shows that the equivalent transform
$\widetilde{\mathcal{F}}$ in $\mathcal{H}_{Z}$ is given by a very simple
formula: it has a phase factor, and otherwise only involves switching of the
two variables $x$, $\omega$, i.e., time and frequency variables,
$z=e^{i2\pi\omega}$. It was included in Table \textup{\ref{LemPoof.1(table)}}
for later use.
\end{remark}

\section[sub-isometries]{\label{R-is}SUB-ISOMETRIES}

Let $m_{0}\in L^{\infty}\left(  \mathbb{T}\right)  $ be a low-pass filter,
i.e., satisfying conditions (\ref{eqBack.1}) and (\ref{IntrAxiom(1)}) in
Section \ref{Intr}, and let $R$, $M$ be the corresponding Ruelle operator and
\emph{cascade refinement operator;} see (\ref{eqIntr.4}) and (\ref{eqBack.9})
for details. Then $L^{\infty}\left(  \mathbb{T}\right)  $ is represented as an
algebra of operators on $L^{2}\left(  \mathbb{R}\right)  $ via $\pi\left(
\alpha\right)  h=\alpha\ast h$, $\alpha\in L^{\infty}\left(  \mathbb{T}%
\right)  $, $h\in L^{2}\left(  \mathbb{R}\right)  $, \emph{cf.}%
\ (\ref{eqBack.12}) or (\ref{eqZak.21}).
Then
we get the following two properties for $M$:%
\begin{align}
M^{\ast}\pi\left(  \alpha\right)  M  &  =\pi\left(  R^{\ast}\left(
\alpha\right)  \right) \label{eqR-is.1}\\%
\intertext{and}%
M^{\ast}\pi\left(  \alpha\right)   &  =\pi\left(  \alpha\left(  z^{2}\right)
\right)  M^{\ast}, \label{eqR-is.2}%
\end{align}
for all $\alpha\in L^{\infty}\left(  \mathbb{T}\right)  $. Taking
$\alpha=\openone$ in (\ref{eqR-is.1}), we get $M^{\ast}M=\pi\left(  \left|
m_{0}\right|  ^{2}\right)  $. So $M$ is not an isometry unless $\left|
m_{0}\right|  ^{2}=\openone$. The last condition is inconsistent with the
low-pass property (\ref{IntrAxiom(1)}) of $m_{0}$, i.e., $m_{0}=\sqrt{2}$ at
$z=1$. But we say that $M$ is a sub-isometry. More generally, let $\pi$ be a
representation on a Hilbert space $\mathcal{H}$.

\begin{definition}
\label{DefR-is.1}
Let $R$ be the Ruelle operator introduced above, i.e.,
$R\colon\mathcal{H}\rightarrow L^{2}\left( \mathbb{R}\right) $,
and let $\pi$ be a
representation of $L^{\infty}\left( \mathbb{T}\right) $ in the algebra of operators on
$\mathcal
H$, such that identities (\ref{DefR-is.1(1)}), (\ref{DefR-is.1(2)})
hold for all $\alpha\in L^{\infty}\left( \mathbb{T}\right) $.
 
We say that an operator $M$ on $\mathcal{H}$ is an
$\left( R,\pi\right) $\emph{-isometry},
or a \emph{sub-isometry} if the data $\left( R,\pi\right) $ is understood,
if

\begin{enumerate}
\item \label{DefR-is.1(1)}$M^{\ast}\pi\left(  \alpha\right)  M=\pi\left(
R^{\ast}\left(  \alpha\right)  \right)  $ on $\mathcal{H}$, and

\item \label{DefR-is.1(2)}$\pi\left(  \alpha\right)  M=M\pi\left(
\alpha\left(  z^{2}\right)  \right)  $ on $\mathcal{H}$.
\end{enumerate}

\noindent\textup{(}In the general case, the two conditions
\textup{(\ref{DefR-is.1(1)})} and \textup{(\ref{DefR-is.1(2)})} are independent.
A discussion of \textup{(\ref{DefR-is.1(2)})} and its variant
\textup{(\ref{eqR-is.2})} will follow.\textup{)}
\end{definition}

If $M$ is an \emph{isometry,} i.e., $M^{\ast}M=\operatorname*{id}%
_{\mathcal{H}}$, then $M\left(  \mathcal{H}\right)  $ is closed, but it may
not be so if $M$ is only a sub-isometry. Then we shall denote the closure
$\left[  M\mathcal{H}\right]  $. Its orthogonal complement is $\ker M^{\ast
}=\left\{  h\in\mathcal{H}\mid M^{\ast}h=0\right\}  $. Let $\mathcal{L}:=\ker
M^{\ast}$, and set%
\begin{equation}
\mathcal{H}^{\left(  \infty\right)  }:=\bigcap_{n=1}^{\infty}\left[
M^{n}\mathcal{H}\right]  . \label{eqR-is.3}%
\end{equation}
Again, if $M$ is an isometry, the classical Wold decomposition of
$\mathcal{H}$ relative to $M$ states the orthogonal decomposition%
\begin{equation}
\mathcal{H}=\sideset{}{^{\smash{\oplus}}}{\sum}\limits_{n=0}^{\infty}%
M^{n}\mathcal{L}\oplus\mathcal{H}^{\left(  \infty\right)  }; \label{eqR-is.4}%
\end{equation}
see \cite{SzFo70} for details. The dimension of $\mathcal{L}$ is then also a
complete invariant for the isometry
in the pure case, i.e., when $\mathcal{H}^{\left( \infty\right) }=\left\{ 0\right\} $.

In this section, we prove an analogue of
this result for general $\left( R,\pi\right) $-isometries. In that form, the components
corresponding to $M^{n}\mathcal{L}$ in (\ref{eqR-is.4}) will instead be
$\left[  M^{n}\mathcal{L}\right]  $. We will still have \emph{orthogonality}
of the subspaces in the decomposition. The important new element for
$\left( R,\pi\right) $-isometries
is that each of the subspaces in the decomposition is invariant
for the representation $\pi$, i.e., $\pi\left(  \alpha\right)  $ maps $\left[
M^{n}\mathcal{L}\right]  $ into itself for all $\alpha\in L^{\infty}\left(
\mathbb{T}\right)  $, and $n=0,1,\dots$. Similarly $\mathcal{H}^{\left(
\infty\right)  }$ is invariant under the operators $\pi\left(  \alpha\right)
$.
The analogy to the classical Wold theorem for isometries raises the question
of whether some \emph{invariant of the representation }$\pi$, when restricted
to $\mathcal{L}=\ker\left(  M^{\ast}\right)  $, is perhaps \emph{a complete
invariant} in the case when $M$ is an
$\left( R,\pi\right) $-isometry. This question is answered
(at least partially) by Theorem \ref{ThmR-is.2}(\ref{ThmR-is.2(3)}) below,
while the first two parts of the theorem give a direct analogue of the Wold
theorem itself in this new representation-theoretic framework.

\begin{theorem}
\label{ThmR-is.2}Let $m_{0}\in L^{\infty}\left(  \mathbb{T}\right)  $ be a
given low-pass filter satisfying the quadratic equation%
\begin{equation}
\left|  m_{0}\left(  z\right)  \right|  ^{2}+\left|  m_{0}\left(  -z\right)
\right|  ^{2}=2\text{\qquad a.e.\ on }\mathbb{T}, \label{eqR-is.5}%
\end{equation}
and let $R$ be the corresponding Ruelle operator. Let $\pi$ be a
representation of $L^{\infty}\left(  \mathbb{T}\right)  $ on a Hilbert space
$\mathcal{H}$, and let $M$ be an associated
$\left( R,\pi\right) $-isometry, i.e., satisfying
conditions \textup{(\ref{DefR-is.1(1)})--(\ref{DefR-is.1(2)})} in Definition
\textup{\ref{DefR-is.1}.}

Then $\mathcal{H}$ has an \emph{orthogonal} decomposition:%
\begin{equation}
\mathcal{H}=\sideset{}{^{\smash{\oplus}}}{\sum}\limits_{n=0}^{\infty}\left[
M^{n}\mathcal{L}\right]  \oplus\mathcal{H}^{\left(  \infty\right)  },
\label{eqR-is.6}%
\end{equation}
where $\mathcal{L}=\ker\left(  M^{\ast}\right)  $, and $\mathcal{H}^{\left(
\infty\right)  }=\bigcap_{n=1}^{\infty}\left[  M^{n}\mathcal{H}\right]  $. It
has the following three properties:\renewcommand{\theenumi}{\alph{enumi}}

\begin{enumerate}
\item \label{ThmR-is.2(1)}the individual closed subspaces in the decomposition
are mutually orthogonal, i.e., $\left[  M^{n}\mathcal{L}\right]  $ is
orthogonal to $\left[  M^{k}\mathcal{L}\right]  $ if $n\neq k$, and they are
all orthogonal to $\mathcal{H}^{\left(  \infty\right)  }$;

\item \label{ThmR-is.2(2)}each of the spaces $\left[  M^{n}\mathcal{L}\right]
$, for $n=0,1,\dots$, and $\mathcal{H}^{\left(  \infty\right)  }$ is invariant
under $\pi\left(  \alpha\right)  $ for all $\alpha\in L^{\infty}\left(
\mathbb{T}\right)  $; and

\item \label{ThmR-is.2(3)}every representation $\pi_{0}$ of $L^{\infty}\left(
\mathbb{T}\right)  $ in a Hilbert space $\mathcal{L}$ arises as the $n=0$ term
of \textup{(\ref{eqR-is.6})} for some $\left( R,\pi\right) $-isometry $M$.
\end{enumerate}
\end{theorem}

\begin{proof}
In the proof, we shall refer to the two properties (\ref{DefR-is.1(1)}%
)--(\ref{DefR-is.1(2)}) in Definition \ref{DefR-is.1}. If $\mathcal{S}%
\subset\mathcal{H}$ is a linear subspace, the orthogonal complement will be
denoted \[\mathcal{S}^{\perp}=\mathcal{H}\ominus\mathcal{S}=\left\{
h\in\mathcal{H}\mid\left\langle h,s\right\rangle =0,\;\forall s\in
\mathcal{S}\right\}  .\] We clearly have $\mathcal{L}=\left(  M\mathcal{H}%
\right)  ^{\perp}$.

\begin{claim}
\label{Cla1}$\mathcal{L}$ is invariant under $\pi\left(  \alpha\right)  $,
$\alpha\in L^{\infty}\left(  \mathbb{T}\right)  $.
\end{claim}

\begin{proof}
Let $l\in\mathcal{L}$. To show that $\pi\left(  \alpha\right)  l\in
\mathcal{L}$, we check that
\[
M^{\ast}\pi\left(  \alpha\right)  l=\pi\left(  \alpha\left(  z^{2}\right)
\right)  M^{\ast}l=0.
\]
We used (\ref{DefR-is.1(2)}) in the calculation, noting that $M^{\ast}l=0$ by definition.
\end{proof}

Our next assertion is this.

\begin{claim}
\label{Cla2}$\mathcal{L}\oplus\left[  M\mathcal{L}\right]  =\left(
M^{2}\mathcal{H}\right)  ^{\perp}$.
\end{claim}

\begin{proof}
We prove the claim by showing that a vector $h_{0}$ which is orthogonal to all
three subspaces $\mathcal{L}$, $M\mathcal{L}$, and $M^{2}\mathcal{H}$, must be
zero. The three subspaces are pairwise mutually orthogonal. This is immediate
from the definitions except for the last pair. Let $l\in\mathcal{L}$ and
$h\in\mathcal{H}$. Then
\begin{align*}
\ip{Ml}{M^{2}h}  &  =\ip{M^{\ast\,2}Ml}{h}\\
&  =\ip{M^{\ast}M^{\ast}Ml}{h}\\
&  =\ip{M^{\ast}\pi\left( \left| m_{0}\right| ^{2}\right) l}{h}\\
&  =\ip{\pi\left( \left| m_{0}\left( z^{2}\right) \right| ^{2}\right) M^{\ast
}l}{h}\\
&  =0,
\end{align*}
where we used properties (\ref{DefR-is.1(1)}) and (\ref{DefR-is.1(2)}), in
that order. In the last step, we used $M^{\ast}l=0$.

The condition that $h_{0}$ is orthogonal to all three subspaces amounts to:
$h_{0}\in\left[  M\mathcal{H}\right]  $, $M^{\ast}h_{0}\in\left[
M\mathcal{H}\right]  $, and $M^{\ast\,2}h_{0}=0$. Since $\mathcal{H}%
=\mathcal{L}\oplus\left[  M\mathcal{H}\right]  $, we have $\left[
M\mathcal{H}\right]  \subset\left[  M\mathcal{L}\right]  \oplus\left[
M^{2}\mathcal{H}\right]  $, so $h_{0}=u+v$, $u\in\left[  M\mathcal{L}\right]
$, $v\in\left[  M^{2}\mathcal{H}\right]  $. Since $M^{\ast}h_{0}\in\left[
M\mathcal{H}\right]  $, we get $h_{0}\in\left(  M\mathcal{L}\right)  ^{\perp}%
$. Indeed there is a sequence $h_{i}\in\mathcal{H}$ such that $M^{\ast}%
h_{0}=\lim_{i}Mh_{i}$. So if $l\in\mathcal{L}$, then%
\[
\ip{Ml}{h_{0}}=\ip{l}{M^{\ast}h_{0}}
 =\lim_{i}\ip{l}{Mh_{i}}
 =\lim_{i}\ip{M^{\ast}l}{h_{i}}
 =0,
\]
since $l\in\ker\left(  M^{\ast}\right)  $. Hence $u=0$ in the decomposition of
$h_{0}$: $h_{0}=0+v$, $v\in\left[  M^{2}\mathcal{H}\right]  $. But $h_{0}%
\in\ker\left(  M^{\ast\,2}\right)  =\left(  M^{2}\mathcal{H}\right)  ^{\perp}%
$, and we conclude that $h_{0}=0$.
\end{proof}

\begin{claim}
\label{Cla3}$\left[  M\mathcal{L}\right]  $ is invariant under $\pi\left(
\alpha\right)  $, $\alpha\in L^{\infty}\left(  \mathbb{T}\right)  $.
\end{claim}

\begin{proof}
We first show that, if $l_{0}\in\mathcal{L}$ and $\alpha\in L^{\infty}\left(
\mathbb{T}\right)  $, then $\pi\left(  \alpha\right)  Ml_{0}\in\mathcal{L}%
\oplus\left[  M\mathcal{L}\right]  $. Using Claim \ref{Cla2}, we do this by
showing that $\pi\left(  \alpha\right)  Ml_{0}\in\ker\left(  M^{\ast
\,2}\right)  $. But
\[
M^{\ast}\pi\left(  \alpha\right)  Ml_{0}=\pi\left(  R^{\ast}\left(
\alpha\right)  \right)  l_{0}\in\mathcal{L},
\]
where we used first (\ref{DefR-is.1(1)}) and then Claim \ref{Cla1}. Then%
\[
M^{\ast\,2}\pi\left(  \alpha\right)  Ml_{0}=\pi\left(  R^{\ast}\left(
\alpha\right)  \left(  z^{2}\right)  \right)  M^{\ast}l_{0}=0,
\]
where we could use (\ref{DefR-is.1(2)}) or, alternatively, Claim \ref{Cla1}.

This means that $\pi\left(  \alpha\right)  |_{\mathcal{L}\oplus\left[
M\mathcal{L}\right]  }$ has an operator block matrix relative to the
orthogonal decomposition $\mathcal{L}\oplus\left[  M\mathcal{L}\right]  $ of
the form%
\[
\pi\left(  \alpha\right)  \sim%
\begin{pmatrix}
A & B\\
0 & D
\end{pmatrix}
,
\]
with $B\colon\left[  M\mathcal{L}\right]  \rightarrow\mathcal{L}$, and the
diagonal operators $A$ and $D$ being endomorphisms of the respective spaces
$\mathcal{L}$ and $\left[  M\mathcal{L}\right]  $. Since $\pi\left(
\alpha\right)  ^{\ast}=\pi\left(  \bar{\alpha}\right)  $, the adjoint
$\left(
\begin{smallmatrix}
A^{\ast} & 0\\
B^{\ast} & D^{\ast}%
\end{smallmatrix}
\right)  $ must be of the same form, and that forces $B=0$, i.e., $\pi\left(
\alpha\right)  =\left(
\begin{smallmatrix}
A & 0\\
0 & D
\end{smallmatrix}
\right)  $, and each of the spaces $\mathcal{L}$ and $\left[  M\mathcal{L}%
\right]  $ is then invariant under $\pi\left(  \alpha\right)  $. In
particular, $\left[  M\mathcal{L}\right]  $ is invariant, which is the claim.
\end{proof}

Our next assertion merges Claims \ref{Cla2} and \ref{Cla3} into the following induction:

\begin{claim}
\label{Cla4}For each $n=1,2,\dots$, we have the decomposition%
\[
\mathcal{L}\oplus\left[  M\mathcal{L}\right]  \oplus\dots\oplus\left[
M^{n-1}\mathcal{L}\right]  =\left(  M^{n}\mathcal{H}\right)  ^{\perp},
\]
where the terms in the decomposition are mutually pairwise orthogonal, and
further each of the spaces $\left[  M^{k}\mathcal{L}\right]  $ is invariant
under $\pi\left(  \alpha\right)  $, $\alpha\in L^{\infty}\left(
\mathbb{T}\right)  $, $k=0,1,\dots,n-1$, where we set $\left[  M^{0}%
\mathcal{L}\right]  =\mathcal{L}$.
\end{claim}

\begin{proof}
This is a simple induction which is based on Claims \ref{Cla1}--\ref{Cla3},
and it is left to the reader. Alternatively, we can prove it by using the
earlier claims on the operators $M^{2},M^{3},\dots$.
\end{proof}

\begin{claim}
\label{Cla5}We have the decomposition \textup{(\ref{eqR-is.6})} of the theorem
with the two properties \textup{(\ref{ThmR-is.2(1)})--(\ref{ThmR-is.2(2)}).}
\end{claim}

\begin{proof}
For each $n$, let $Q_{n}$ denote the projection onto $\left[  M^{n}%
\mathcal{H}\right]  $. Since $\left[  M^{n+1}\mathcal{H}\right]
\subset\left[  M^{n}\mathcal{H}\right]  $, this is a decreasing sequence of
projections in $\mathcal{H}$, By Hilbert space theory (see, e.g.,
\cite{SzFo70}), it has a limit $Q_{\infty}$, i.e., $\lim_{n\rightarrow\infty
}\left\|  Q_{n}h-Q_{\infty}h\right\|  =0$ for all $h\in\mathcal{H}$, and
$Q_{\infty}$ is the orthogonal projection onto $\mathcal{H}^{\left(
\infty\right)  }=\bigcap_{n=1}^{\infty}Q_{n}\mathcal{H}$. Recall that
$Q_{n}\mathcal{H}=\left[  M^{n}\mathcal{H}\right]  $, by definition! In fact,%
\[
\left\|  Q_{\infty}h\right\|  =\inf_{n}\left\|  Q_{n}h\right\|  ,\qquad
h\in\mathcal{H}.
\]
But Claim \ref{Cla4} states that $I-Q_{n}$ is the projection onto
$\mathcal{L}\oplus\left[  M\mathcal{L}\right]  \oplus\cdots\oplus\left[
M^{n-1}\mathcal{L}\right]  $, where we write $I$ for the identity operator in
$\mathcal{H}$. Now $I-Q_{n}$ is an increasing family of projections, and its
limit $I-Q_{\infty}$ is the projection onto $\sideset{}{^{\smash{\oplus}}%
}{\sum}\limits_{n=0}^{\infty}\left[  M^{n}\mathcal{L}\right]  $. Let
$\mathcal{B}$ denote this space. Vectors $b$ in $\mathcal{B}$ are
characterized by $\left(  I-Q_{\infty}\right)  b=b$, or equivalently
$Q_{\infty}b=0$, and each $b$ has the unique representation $b=\sum
_{0}^{\infty}b_{n}$, $\left\|  b\right\|  ^{2}=\sum_{0}^{\infty}\left\|
b_{n}\right\|  ^{2}$, $b_{n}\in\left[  M^{n}\mathcal{L}\right]  $,
$n=0,1,\dots$.

Since each of the spaces $\left[  M^{n}\mathcal{L}\right]  $ is invariant
under $\pi\left(  \alpha\right)  $, $\alpha\in L^{\infty}\left(
\mathbb{T}\right)  $, by Claims \ref{Cla3}--\ref{Cla4}, it follows that
$\mathcal{B}$ is also $\pi\left(  \alpha\right)  $-invariant. Since
$\pi\left(  \alpha\right)  ^{\ast}=\pi\left(  \bar{\alpha}\right)  $, it
follows that%
\begin{equation}
\mathcal{B}^{\perp}=\mathcal{H}^{\left(  \infty\right)  }=\bigcap
_{n=1}^{\infty}\left[  M^{n}\mathcal{H}\right]  \label{eqR-is.7}%
\end{equation}
is also $\pi\left(  \alpha\right)  $-invariant.
\end{proof}

The proof of Theorem \ref{ThmR-is.2}(\ref{ThmR-is.2(3)}) will be given after
the next three corollaries.

\begin{corollary}
\label{CorR-is.3}For each $n=1,2,\dots$, the space $\left[  M^{n}%
\mathcal{H}\right]  $ is invariant under $\pi\left(  \alpha\right)  $,
$\alpha\in L^{\infty}\left(  \mathbb{T}\right)  $.
\end{corollary}

\begin{proof}
We showed in Claim \ref{Cla4} that%
\begin{equation}
\left(  M^{n}\mathcal{H}\right)  ^{\perp}=\sideset{}{^{\smash{\oplus}}}{\sum
}\limits_{k=0}^{n-1}\left[  M^{k}\mathcal{L}\right]  , \label{eqR-is.8}%
\end{equation}
and that the right-hand side has the $\pi\left(  \alpha\right)  $-invariance.
Since $\pi\left(  \alpha\right)  ^{\ast}=\pi\left(  \bar{\alpha}\right)  $, we
conclude that $\left(  M^{n}\mathcal{H}\right)  ^{\perp\,\perp}=\left[
M^{n}\mathcal{H}\right]  $ is also $\pi\left(  L^{\infty}\left(
\mathbb{T}\right)  \right)  $-invariant.
\end{proof}

Let $\left(  R,\pi,M,\mathcal{H}\right)  $ be as in the statement of Theorem
\ref{ThmR-is.2}, i.e., $M\colon\mathcal{H}\rightarrow\mathcal{H}$ is an
$\left( R,\pi\right) $-isometry
relative to some Ruelle operator $R$ and representation $\pi$. We
say that a closed subspace $\mathcal{S}\subset\mathcal{H}$ is \emph{double
invariant} if $\mathcal{S}$ is invariant under both $M$ and $M^{\ast}$. It is
then immediate from the definition of $M^{\ast}$ that $\mathcal{S}$ is double
invariant if and only if both $\mathcal{S}$ and $\mathcal{S}^{\perp}$
($=\mathcal{H}\ominus\mathcal{S}$) are invariant under $M$, i.e., $M\left(
\mathcal{S}\right)  \subset\mathcal{S}$, and $M\left(  \mathcal{S}^{\perp
}\right)  \subset\mathcal{S}^{\perp}$.

\begin{corollary}
\label{CorR-is.4}Let $\left(  R,\pi,M,\mathcal{H}\right)  $ be as described,
and let $\mathcal{B}$ \textup{(}$=\smash{\sideset{}{^{\smash{\oplus}}}{\sum
}\limits_{0}^{\infty}}\left[  M^{n}\mathcal{L}\right]  $\textup{)} and
$\mathcal{H}^{\left(  \infty\right)  }$ \textup{(}$=\bigcap_{1}^{\infty
}\left[  M^{n}\mathcal{H}\right]  $\textup{)} be as in Theorem
\textup{\ref{ThmR-is.2}.}

Then both $\mathcal{B}$ and $\mathcal{H}^{\left(  \infty\right)  }$ are double
invariant under $M$.
\end{corollary}

\begin{proof}
From the comment before the statement of the Corollary, it is enough to show
that both $\mathcal{B}$ and $\mathcal{H}^{\left(  \infty\right)  }$ are
invariant under $M$. Recall from Theorem \ref{ThmR-is.2} that $\mathcal{B}%
=\left(  \mathcal{H}^{\left(  \infty\right)  }\right)  ^{\perp}$, and
$\mathcal{H}^{\left(  \infty\right)  }=\mathcal{B}^{\perp}$. But it is clear
from the definition of $\mathcal{H}^{\left(  \infty\right)  }$ that $M\left(
\mathcal{H}^{\left(  \infty\right)  }\right)  \subset\mathcal{H}^{\left(
\infty\right)  }$. Since arbitrary vectors $b$ in $\mathcal{B}$ may be
represented as $b=\sum_{0}^{\infty}M^{n}l_{n}$, $l_{n}\in\mathcal{L}$,
$\left\|  b\right\|  ^{2}=\sum_{0}^{\infty}\left\|  M^{n}l_{n}\right\|  ^{2}$,
it follows that $Mb=\sum_{0}^{\infty}M^{n+1}l_{n}$. If $b$ is encoded with the
sequence $\left(  l_{0},l_{1},l_{2},\dots\right)  $, $l_{n}\in\mathcal{L}$,
then $Mb\sim\left(  0,l_{0},l_{1},\dots\right)  $, i.e.,%
\begin{equation}
M\left(  \left(  l_{0},l_{1},l_{2},\dots\right)  \right)  =\left(
0,l_{0},l_{1},\dots\right)  , \label{eqNewR-is.9}%
\end{equation}
and the $M$-invariance for $\mathcal{B}$ follows from this. By the initial
argument we conclude that both $\mathcal{B}$ and $\mathcal{H}^{\left(
\infty\right)  }$ are double invariant.
\end{proof}

The advantage of the $\left(  l_{0},l_{1},l_{2},\dots\right)  $ representation
of $\mathcal{B}$ is that $M^{\ast}|_{\mathcal{B}}$ takes an especially simple form:

\begin{corollary}
\label{CorR-isMstar}If vectors in $\mathcal{B}$ are represented in the form
$\left(  l_{0},l_{1},l_{2},\dots\right)  $,\linebreak $\sum2^{n}\left\|
l_{n}\right\|  ^{2}<\infty$, $l_{n}\in\mathcal{L}=\ker\left(  M^{\ast}\right)
$, then the action of $M^{\ast}$ on $\mathcal{B}$ is%
\begin{equation}
\left(  l_{0},l_{1},l_{2},\dots\right)  \longmapsto\left(  \pi\left(  \left|
m_{0}\right|  ^{2}\right)  l_{1},\pi\left(  \left|  m_{0}\left(  z^{2}\right)
\right|  ^{2}\right)  l_{2},\pi\left(  \left|  m_{0}\left(  z^{4}\right)
\right|  ^{2}\right)  l_{3},\dots\right)  . \label{eqNewR-is.10}%
\end{equation}
\end{corollary}

\begin{proof}
The proof follows from the following calculation: $M^{\ast}l_{0}=0$, and
$M^{\ast}M^{n}l_{n}=$\linebreak $\pi\left(  \left|  m_{0}\right|  ^{2}\right)  M^{n-1}%
l_{n}=M^{n-1}\pi\left(  \left|  m_{0}\left(  z^{2^{n-1}}\right)  \right|
^{2}\right)  l_{n}$.
\end{proof}

\begin{proof}
[Proof of Theorem \textup{\ref{ThmR-is.2}(\ref{ThmR-is.2(3)})}]Let the
representation $\pi_{0}$ of $L^{\infty}\left(  \mathbb{T}\right)  $ in
$\mathcal{L}$ be given as in Theorem \ref{ThmR-is.2}(\ref{ThmR-is.2(3)}). Let
$m_{0}$, $h$ be as stated at the outset, i.e., $R_{m_{0}}\left(  h\right)
=h$. On the vectors described in (\ref{eqNewR-is.9}), define the Hilbert-space
norm, and corresponding completion, by%
\[
\left\|  \left(  l_{n}\right)  _{0}^{\infty}\right\|  ^{2}:=\sum_{n=0}%
^{\infty}\left\|  \pi_{0}\left(  m_{0}^{\left(  n\right)  }\right)
l_{n}\right\|  ^{2},
\]
and define $M$ as in (\ref{eqNewR-is.9}). A simple computation, using the
corresponding inner product%
\[
\ip{\left( l_{n}^{{}}\right) }{\left( l_{n}^{\prime}\right) }:=\sum
_{n=0}^{\infty}\ip{\pi_{0}^{{}}\left( m_{0}^{\left( n\right) }\right
) l_{n}^{{}}}{\pi_{0}^{{}}\left( m_{0}^{\left( n\right) }\right) l_{n}%
^{\prime}},
\]
and Table \ref{tableZak.2} then yields an adjoint operator $M^{\ast}$ which
turns out to be (\ref{eqNewR-is.10}). It is now a simple matter to verify that
$M$ is the desired sub-isometry.
\end{proof}

\renewcommand{\qed}{}
\end{proof}

\begin{remark}
\label{RemR-is.5}It follows from Corollary \textup{\ref{CorR-is.3}} that each
projection $Q_{n}$ \textup{(}onto the space $\left[  M^{n}\mathcal{H}\right]
$\textup{)} commutes with $\pi\left(  \alpha\right)  $, $\alpha\in L^{\infty
}\left(  \mathbb{T}\right)  $, i.e., $Q_{n}\pi\left(  \alpha\right)
=\pi\left(  \alpha\right)  Q_{n}$, but it is generally \emph{not} the case
that $M$ commutes with $\pi\left(  \alpha\right)  $. However, $M$ may possibly
commute with a special $\pi\left(  \alpha_{0}\right)  $ for some $\alpha
_{0}\in L^{\infty}\left(  \mathbb{T}\right)  $. We have the following simple
result on that.
\end{remark}

\begin{proposition}
\label{ProR-is.6}Let $m_{0}\in L^{\infty}\left(  \mathbb{T}\right)  $ be as
specified in Theorem \textup{\ref{ThmR-is.2},} and let $\pi$ be a
\emph{faithful} representation of $L^{\infty}\left(  \mathbb{T}\right)  $ on a
Hilbert space $\mathcal{H}$. Let $R$ be the Ruelle operator constructed from
$m_{0}$, and let $M$ be a given sub-isometry. Let $\alpha_{0}\in L^{\infty
}\left(  \mathbb{T}\right)  $. If $M\pi\left(  \alpha_{0}\right)  =\pi\left(
\alpha_{0}\right)  M$, then it follows that $R\left(  \alpha_{0}\right)
=\alpha_{0}$, i.e., $\alpha_{0}$ is an eigenvector for $R$ with eigenvalue $1$.
\end{proposition}

\begin{proof}
By (\ref{DefR-is.1(1)}) of Definition \ref{DefR-is.1},
\[
\pi\left(  \left|  m_{0}\right|  ^{2}\alpha_{0}\right)  =M^{\ast}M\pi\left(
\alpha_{0}\right)  =M^{\ast}\pi\left(  \alpha_{0}\right)  M=\pi\left(
R^{\ast}\left(  \alpha_{0}\right)  \right)  .
\]
In the first step, (\ref{DefR-is.1(1)}) is used on $\openone$, then
commutativity is used, and in the last step, (\ref{DefR-is.1(1)}) is used on
$\alpha_{0}$. Since $\pi$ is assumed faithful,%
\begin{equation}
\left|  m_{0}\right|  ^{2}\alpha_{0}=R^{\ast}\left(  \alpha_{0}\right)
=\left|  m_{0}\left(  z\right)  \right|  ^{2}\alpha_{0}\left(  z^{2}\right)  ,
\label{eqR-is.9}%
\end{equation}
and%
\begin{align*}
R\left(  \alpha_{0}\right)  \left(  z\right)   &  =\frac{1}{2}\sum_{w^{2}%
=z}\left|  m_{0}\left(  w\right)  \right|  ^{2}\alpha_{0}\left(  w\right) \\
&  =\frac{1}{2}\sum_{w^{2}=z}\left|  m_{0}\left(  w\right)  \right|
^{2}\alpha_{0}\left(  z\right) \\
&  =\alpha_{0}\left(  z\right)  ,
\end{align*}
where the quadratic property of $m_{0}$ was used in the last step. Hence
$R\left(  \alpha_{0}\right)  =\alpha_{0}$ as claimed.
\end{proof}

\begin{remark}
\label{RemR-is.7}\emph{The converse implication to the one given in
Proposition} \textup{\ref{ProR-is.6}} \emph{is not true.}

Let $m_{0}\left(  z\right)  =\frac{1}{\sqrt{2}\mathstrut}\left(
1+z^{3}\right)  $ \textup{(}see also Section \textup{\ref{Zak}),} and let $R$
and $M$ be the corresponding Ruelle operator and cascade operator. The scaling
function $\varphi$ realized in $L^{2}\left(  \mathbb{R}\right)  $ is
$\varphi=\frac{1}{3}\chi_{\left[  0,3\right]  }^{{}}$, and a little
calculation \textup{(}see Section \textup{\ref{Zak})} shows that%
\[
p_{2}\left(  \varphi\right)  =p\left(  \varphi,\varphi\right)  =\frac{1}%
{9}\left(  z^{-2}+2z^{-1}+3+2z+z^{2}\right)  ,
\]
where $z=e^{-i\omega}$, $\omega\in\mathbb{R}$. But $R\left(  p_{2}\left(
\varphi\right)  \right)  =p_{2}\left(  \varphi\right)  $, so this is an
eigenfunction for $R$. Let $\alpha_{0}:=p_{2}\left(  \varphi\right)  $. We
claim that $\pi\left(  \alpha_{0}\right)  $ does not commute with $M$. In
fact, commutativity with $M$ is equivalent to identity \textup{(\ref{eqR-is.9}%
)} from the proof of Proposition \textup{\ref{ProR-is.6}.} An inspection shows
that \textup{(\ref{eqR-is.9})} is not satisfied in this example. In other
words,%
\begin{equation}
\left|  m_{0}\left(  z\right)  \right|  ^{2}\alpha_{0}\left(  z\right)
\neq\left|  m_{0}\left(  z\right)  \right|  ^{2}\alpha_{0}\left(
z^{2}\right)  . \label{eqR-is.10}%
\end{equation}
Both sides in \textup{(\ref{eqR-is.10})} are polynomials, i.e., in
$\mathbb{C}\left[  z,z^{-1}\right]  $. The right-hand side contains a term
$z^{7}$ whereas the left-hand side does not.
\end{remark}

\begin{remark}
\label{RemR-is.8}The eigenvalue problem for $R$ plays a crucial role for
approximation of wavelets; see, e.g., \textup{\cite{Str96}, \cite{CoDa96},
\cite{Vil94}.}
\end{remark}

In the study of refinement operators, the sub-isometries usually have a
slightly different formulation, and for the particular cascade operator $M$ in
(\ref{eqBack.9}), the alternative formulation is stated in Lemma
\ref{LemZak.1}; see (\ref{eqZak.11}).

We now specialize the general formulation $\left(  \pi,M,\mathcal{H}\right)  $
of the present section. Recall that $\pi$ was a representation of $L^{\infty
}\left(  \mathbb{T}\right)  $ on a Hilbert space $\mathcal{H}$, and
$M\colon\mathcal{H}\rightarrow\mathcal{H}$ was an operator satisfying
(\ref{DefR-is.1(1)})--(\ref{DefR-is.1(2)}) of Definition \ref{DefR-is.1}
relative to some given Ruelle operator $R=R_{m_{0}}$. We will specialize as
follows: $\mathcal{H=H}_{Z}$ (the Hilbert space of Section \ref{Back}), and%
\begin{equation}
\left(  \pi_{Z}\left(  \alpha\right)  H\left(  z,x\right)  \right)
=\alpha\left(  z\right)  H\left(  z,x\right)  ,\qquad\alpha\in L^{\infty
}\left(  \mathbb{T}\right)  ,\;H\in\mathcal{H}_{Z}. \label{eqR-is.11}%
\end{equation}
In this specialized setup, we then have the following.

\begin{proposition}
\label{ProR-is.9}Let the pair $\left(  \pi_{Z},\mathcal{H}_{Z}\right)  $ be as
described, and let $M$ be an operator in $\mathcal{H}_{Z}$. Then the following
two conditions are equivalent.

\begin{enumerate}
\item \label{ProR-is.9(1)}For a.e.\ $z$ in $\mathbb{T}$, we have the
\mbox{identity}
\[
\ip{MH_{1}\left( z\right) }{MH_{2}\left( z\right) }_{L^{2}\left( 0,1\right
) }=R\left( \ip{H_{1}\left( \,\cdot\,\right) }{H_{2}\left( \,\cdot\,\right
) }_{L^{2}\left( 0,1\right) }\right) \left( z\right)
\]%
for all $H_{1},H_{2}\in\mathcal{H}_{Z}$.

\item \label{ProR-is.9(2)}$M^{\ast}\pi_{Z}\left(  \alpha\right)  M=\pi
_{Z}\left(  R^{\ast}\alpha\right)  $ for all $\alpha\in L^{\infty}\left(
\mathbb{T}\right)  $.
\end{enumerate}

\noindent Here the Ruelle operator is defined from an arbitrary $m_{0}$ as usual.
\end{proposition}

\begin{proof}
(\ref{ProR-is.9(1)}) $\Rightarrow$ (\ref{ProR-is.9(2)}). Let $H_{1},H_{2}%
\in\mathcal{H}_{Z}$. Then
\begin{align*}
\ip{H_{1}}{\left( M^{\ast}\pi_{Z}\left( \alpha\right) M\right) H_{2}%
}_{\mathcal{H}_{Z}}  &  =\int_{\mathbb{T}}\ip{MH_{1}\left( z\right) }%
{MH_{2}\left( z\right) }\alpha\left(  z\right)  \,d\mu\left(  z\right) \\
&  =\int_{\mathbb{T}}R\left(  \ip{H_{1}\left( \,\cdot\,\right) }{H_{2}%
\left( \,\cdot\,\right) }\right)  \alpha\left(  z\right)  \,d\mu\left(
z\right) \\
&  =\int_{\mathbb{T}}\ip{H_{1}\left( z\right) }{H_{2}\left( z\right) }\left|
m_{0}\left(  z\right)  \right|  ^{2}\alpha\left(  z^{2}\right)  \,d\mu\left(
z\right) \\
&  =\ip{H_{1}}{\pi_{Z}\left( R^{\ast}\left( \alpha\right) \right) H_{2}%
}_{\mathcal{H}_{Z}},
\end{align*}
and this proves (\ref{ProR-is.9(2)}).

(\ref{ProR-is.9(2)}) $\Rightarrow$ (\ref{ProR-is.9(1)}). If
(\ref{ProR-is.9(2)}) holds, the calculation shows that the second and the
third terms must agree for all $\alpha\in L^{\infty}\left(  \mathbb{T}\right)
$, and, by duality, this means that (\ref{ProR-is.9(1)}) must hold a.e.\ on
$\mathbb{T}$.
\end{proof}

\section[Singular cascade approximations]{\label{SiC}SINGULAR CASCADE APPROXIMATIONS}

We return in this section to the cascade operator $M$ from Section \ref{Intr},
but it will be convenient to state the results for the Hilbert space
$\mathcal{H}_{Z}$ of Section \ref{Back}. We will also need the representation
$\pi_{Z}$ of (\ref{eqBack.12}). From the dictionary in Lemma \ref{LemPoof.1}
(Table \ref{LemPoof.1(table)}) we note that $\mathbb{Z}$-translations in
$L^{2}\left(  \mathbb{R}\right)  $, $h\left(  \,\cdot\,\right)  \mapsto
h\left(  \,\cdot\,+n\right)  $, $n\in\mathbb{Z}$, correspond to $\pi
_{Z}\left(  e_{n}\right)  $ in $\mathcal{H}_{Z}$ where $e_{n}\left(  z\right)
=z^{n}$. Also we need the fact that the operation $h\mapsto\alpha\ast h$ on
$L^{2}\left(  \mathbb{R}\right)  $, $h\in L^{2}\left(  \mathbb{R}\right)  $,
$\alpha\in L^{\infty}\left(  \mathbb{T}\right)  $, corresponds to
multiplication $\alpha\left(  e^{-i\omega}\right)  \hat{h}\left(
\omega\right)  $ where $\hat{h}\left(  \omega\right)  =\int_{\mathbb{R}%
}e^{-i\omega x}h\left(  x\right)  \,dx$. Hence the result when stated in
$\mathcal{H}_{Z}$ can easily be translated to either $L^{2}\left(
\mathbb{R}\right)  $ or $\widehat{L^{2}\left(  \mathbb{R}\right)  }$.

Let $H_{i}\in\mathcal{H}_{Z}$. We shall need the sesquilinear form $p$ from
Section \ref{Back},%
\[
p\left(  H_{1},H_{2}\right)  \left(  z\right)  =\int_{0}^{1}\overline
{H_{1}\left(  z,x\right)  }H_{2}\left(  z,x\right)  \,dx=\ip{H_{1}%
\left( z\right) }{H_{2}\left( z\right) }_{L^{2}\left(  0,1\right)  }.
\]
If $H=H_{1}=H_{2}$, we introduce the abbreviation $p_{2}\left(  H\right)
=p\left(  H,H\right)  $. An important property of a starting vector $h\in
L^{2}\left(  \mathbb{R}\right)  $ for the cascade algorithm is orthogonality
of the translates $\left\{  h\left(  \,\cdot\,-n\right)  \mid n\in
\mathbb{Z}\right\}  $, i.e.,%
\begin{equation}
\int_{-\infty}^{\infty}\overline{h\left(  x-n\right)  }h\left(  x\right)
\,dx=\delta_{n,0}\left\|  h\right\|  _{2}^{2}. \label{eqSiC.orthogonality}%
\end{equation}
The following is immediate from Proposition \ref{ProZak.4} (Section \ref{Zak}):

\begin{lemma}
\label{LemSiC.1}Let $h\in L^{2}\left(  \mathbb{R}\right)  $, and set $H=Zh$.
The following are equivalent:

\begin{enumerate}
\item \label{LemSiC.1(1)}$h$ satisfies the orthogonality condition
\textup{(\ref{eqSiC.orthogonality});} and

\item \label{LemSiC.1(2)}$p_{2}\left(  H\right)  \left(  z\right)
\equiv\left\|  H\right\|  _{\mathcal{H}_{Z}}^{2}$ a.e.\ on $\mathbb{T}$.
\end{enumerate}

\noindent We shall refer to this as the orthogonality condition, meaning
orthogonal $\mathbb{Z}$-translates in $L^{2}\left(  \mathbb{R}\right)  $.
\end{lemma}

\begin{remark}
\label{RemSiC.2}We note similarly that when the Strang--Fix condition%
\begin{equation}
\sum_{n\in\mathbb{Z}}h\left(  x+n\right)  \equiv1 \label{eqSiC.StrangFix}%
\end{equation}
makes sense for $h\in L^{2}\left(  \mathbb{R}\right)  $, e.g., if $h$ is of
compact support, then this corresponds to the following condition on $H=Zh$:%
\begin{equation}
H\left(  z=1,x\right)  \equiv1\text{\qquad a.a.\ }x\in\left[  0,1\right]  .
\label{eqSiC.Haax}%
\end{equation}

It is well known that the cascade approximation must start with the
orthogonality condition of Lemma \textup{\ref{LemSiC.1}.} Suppose
$H\in\mathcal{H}_{Z}$ satisfies \textup{(\ref{eqSiC.Haax}),} and%
\begin{equation}
H=H_{\mathcal{B}}+H_{\infty} \label{eqSiC.decomposition}%
\end{equation}
is the Wold decomposition from \textup{(\ref{eqR-is.6})} in Theorem
\textup{\ref{ThmR-is.2},} i.e.,%
\begin{align*}
H_{\mathcal{B}}  &  \in\mathcal{B}=\sideset{}{^{\smash{\oplus}}}{\sum
}\limits_{n=0}^{\infty}\left[  M^{n}\mathcal{L}\right]  ,\\%
\intertext{and}%
H_{\infty}  &  \in\mathcal{H}^{\left(  \infty\right)  }=\bigcap_{n=1}^{\infty
}\left[  M^{n}\mathcal{H}\right]
\end{align*}
relative to a fixed
$\left( R,\pi\right) $-isometry $M$, where $R=R_{m_{0}}$, and $m_{0}$ is
given. Then, if $H$ satisfies the orthogonality
\textup{(\ref{eqSiC.orthogonality}),} or equivalently
\textup{(\ref{LemSiC.1(2)})} in Lemma \textup{\ref{LemSiC.1},} then one of the
two, $H_{\mathcal{B}}$ or $H_{\infty}$, can satisfy the same only if the other
is zero. This follows from the following formula for the norm in
$\mathcal{H}=\mathcal{H}_{Z}$:%
\begin{equation}
\left\|  H\right\|  _{\mathcal{H}_{Z}}^{2}=\int_{\mathbb{T}}p_{2}\left(
H\right)  \left(  z\right)  \,d\mu\left(  z\right)  . \label{eqSiC.1}%
\end{equation}
Applied to the decomposition $H=H_{\mathcal{B}}+H_{\infty}$ in
\textup{(\ref{eqSiC.decomposition}),} we get%
\begin{equation}
\int_{\mathbb{T}}p_{2}\left(  H\right)  \left(  z\right)  \,d\mu\left(
z\right)  =\int_{\mathbb{T}}p_{2}\left(  H_{\mathcal{B}}\right)  \left(
z\right)  \,d\mu\left(  z\right)  +\int_{\mathbb{T}}p_{2}\left(  H_{\infty
}\right)  \left(  z\right)  \,d\mu\left(  z\right)  , \label{eqSiC.2}%
\end{equation}
from which the claim is clear. If $H$ is orthogonal, then
\[
1=\left\|  H\right\|  ^{2}=\int_{\mathbb{T}}p_{2}\left(  H\right)  \left(
z\right)  \,d\mu\left(  z\right)  ,
\]
so if for example $p_{2}\left(  H_{\mathcal{B}}\right)  \equiv1$, then%
\[
1=1+\int_{\mathbb{T}}p_{2}\left(  H_{\infty}\right)  \left(  z\right)
\,d\mu\left(  z\right)  ,
\]
so $H_{\infty}=0$.
\end{remark}

While Remark \ref{RemSiC.2} above shows that the $\mathcal{B}$-decomposition
is an obstruction to cascade approximation, it is still the case that
solutions $H$ to $M\left(  H\right)  =H$ (i.e., $H=Z\varphi$ for some scaling
function $\varphi\in L^{2}\left(  \mathbb{R}\right)  $) yield

\begin{proposition}
\label{ProSiC.3}Let $m_{0}$ be a low-pass filter, and let $R$ be the
corresponding Ruelle operator. Let $\pi_{Z}$ be the representation
\textup{(\ref{eqBack.12})} or \textup{(\ref{eqZak.21})} of $L^{\infty}\left(
\mathbb{T}\right)  $ in $\mathcal{H}_{Z}$.\renewcommand{\theenumi}%
{\alph{enumi}}

\begin{enumerate}
\item \label{ProSiC.3(1)}Let $M$ be a sub-isometry. Then%
\begin{equation}
R\left(  p_{2}\left(  H\right)  \right)  =p_{2}\left(  MH\right)  ,
\label{eqSiC.3}%
\end{equation}
so if $M\left(  H\right)  =H$, then $p=p_{2}\left(  H\right)  $ solves%
\begin{equation}
R\left(  p\right)  =p. \label{eqSiC.4}%
\end{equation}

\item \label{ProSiC.3(2)}In general, the Fourier series for $p_{2}\left(
H\right)  $ is
\begin{equation}
p_{2}\left(  H\right)  \left(  z\right)  =\sum_{n\in\mathbb{Z}}z^{n}\ip
{\pi\left( e_{n}\right) H}{H}_{\mathcal{H}_{Z}}, \label{eqSiC.5}%
\end{equation}
where $e_{n}\left(  z\right)  =z^{n}$, and%
\[
\left\|  H\right\|  _{\mathcal{H}_{Z}}^{2}=\sum_{n\in\mathbb{Z}}\left|
\ip{\pi\left( e_{n}\right) H}{H}\right|  ^{2}.
\]
\end{enumerate}
\end{proposition}

\begin{proof}
Immediate from Lemma \ref{LemSiC.1} and Proposition \ref{ProR-is.9} in the
previous section. The calculation of the expansion (\ref{eqSiC.5}) follows
from checking that the Fourier coefficients of $p_{2}\left(  H\right)  \left(
z\right)  $ are as stated, i.e., that
\[
\int_{\mathbb{T}}z^{-n}p_{2}\left(  H\right)  \left(  z\right)  \,d\mu\left(
z\right)  =\int_{\mathbb{T}}p\left(  \pi_{Z}\left(  e_{n}\right)  H,H\right)
\left(  z\right)  \,d\mu\left(  z\right)  =\ip{\pi_{Z}\left( e_{n}\right
) H}{H}_{\mathcal{H}_{Z}}.\settowidth{\qedskip}{$\displaystyle
\int_{\mathbb{T}}z^{-n}p_{2}\left(  H\right)  \left(  z\right)  \,d\mu\left(
z\right)  =\int_{\mathbb{T}}p\left(  \pi_{Z}\left(  e_{n}\right)  H,H\right)
\left(  z\right)  \,d\mu\left(  z\right)  =\ip{\pi_{Z}\left( e_{n}\right
) H}{H}_{\mathcal{H}_{Z}}.$}\addtolength{\qedskip}{-\textwidth}\rlap{\hbox
to-0.5\qedskip{\hfil\qed}}%
\]
\renewcommand{\qed}{}
\end{proof}

We show next that $\lambda=1$ is the only point on the unit circle which is an
eigenvalue of $R$, at least when we restrict to continuous eigenfunctions.

\begin{theorem}
\label{ThmSiC.4}Let $m_{0}$ be a continuous low-pass filter, and let $R$ be
the corresponding Ruelle operator in $L^{\infty}\left(  \mathbb{T}\right)  $.
Suppose that the eigenspace $\left\{  \xi\in C\left(  \mathbb{T}\right)  \mid
R\left(  \xi\right)  =\xi\right\}  $ is one-dimensional, i.e., $\mathbb{C}%
\openone$. Then if $\left|  \lambda\right|  \geq1$, and $\lambda\neq1$, the
eigenvalue problem%
\begin{equation}
R\left(  \alpha\right)  =\lambda\alpha\label{eqSiC.6}%
\end{equation}
has no nonzero solution in $C\left(  \mathbb{T}\right)  $.
\end{theorem}

\begin{proof}
First note that if $m_{0}$ is assumed continuous, then $R$ maps $C\left(
\mathbb{T}\right)  $ into itself. We see this by first approximating $m_{0}$
with finite sums $\sum_{k}a_{k}z^{k}$. For each such finite sum, the
corresponding $R$ maps $\mathbb{C}\left[  z,z^{-1}\right]  $ ($=$ the ring of
finite Fourier series) into itself. Since the norm of $R$, as a operator in
$L^{\infty}\left(  \mathbb{T}\right)  $, is one, i.e., $\left\|  R\right\|
_{\infty,\infty}=1$, we conclude that $R$ maps the norm closure of
$\mathbb{C}\left[  z,z^{-1}\right]  $ into itself. This norm-closure is
$C\left(  \mathbb{T}\right)  $ by the Stone--Weierstrass theorem.

Since $R$ leaves $C\left(  \mathbb{T}\right)  $ invariant, it dualizes to a
map on the measures $M\left(  \mathbb{T}\right)  =C\left(  \mathbb{T}\right)
^{\ast}$, and we claim that $R^{\ast}\left(  \delta_{1}\right)  =\delta_{1}$
where $\delta_{1}$ denotes the Dirac point measure at $1$. Indeed, let $\xi\in
C\left(  \mathbb{T}\right)  $; then%
\[
\left(  R\xi\right)  \left(  1\right)  =\frac{1}{2}\sum_{w^{2}=1}\left|
m_{0}\left(  w\right)  \right|  ^{2}\xi\left(  w\right)  =\frac{1}{2}\left(
\left|  m_{0}\left(  1\right)  \right|  ^{2}\xi\left(  1\right)  +\left|
m_{0}\left(  -1\right)  \right|  ^{2}\xi\left(  -1\right)  \right)
=\xi\left(  1\right)
\]
since $m_{0}\left(  1\right)  =\sqrt{2}$ and $m_{0}\left(  -1\right)  =0$.
(This is the low-pass property.) It follows that $\alpha\left(  1\right)  =0$
when $\alpha$ is the eigenfunction in (\ref{eqSiC.6}).

Suppose we did have a solution $\alpha$ to (\ref{eqSiC.6}) as stated. Let
$c\in\mathbb{C}$, and set $\beta_{c}=\openone+c\alpha$. Then $\beta_{c}\left(
1\right)  =1$. Pick a scaling function $H_{0}$ for the cascade operator $M$
corresponding to $m_{0}$, i.e., $M=M_{m_{0}}$, and $M\left(  H_{0}\right)
=H_{0}$. Since $\mathcal{E}_{1}=\left\{  \xi\in C\left(  \mathbb{T}\right)
\mid R\xi=\xi\right\}  $ is one-dimensional, and $p_{2}\left(  H_{0}\right)
\in\mathcal{E}_{1}$ by Lemma \ref{LemSiC.1}(\ref{LemSiC.1(1)}), we conclude
that $p_{2}\left(  H_{0}\right)  =C\openone$. Lemma \ref{LemSiC.1}%
(\ref{LemSiC.1(2)}) then implies that $C=1$, so $p_{2}\left(  H_{0}\right)
=\openone$. We then have
\begin{align}
R^{n}\left(  \beta_{c}\right)  \left(  z\right)   &  =R^{n}\left(  \beta
_{c}p_{2}\left(  H_{0}\right)  \right)  \left(  z\right)  =R^{n}\left(
p\left(  H_{0},\pi\left(  \beta_{c}\right)  H_{0}\right)  \right)  \left(
z\right)  \label{eqSiC.7}\\
&  =p\left(  H_{0},M^{n}\pi_{Z}\left(  \beta_{c}\right)  H_{0}\right)  \left(
z\right)  ,\nonumber
\end{align}
using the argument from the proof of Proposition \ref{ProSiC.3}. The left-hand
side of (\ref{eqSiC.7}) is $\openone+c\lambda^{n}\alpha$, which is not
convergent when $c\neq0$, $\lambda\neq1$, $\left|  \lambda\right|  \geq1$, and
$\alpha\neq0$ in $L^{\infty}\left(  \mathbb{T}\right)  $.

On the right-hand side in (\ref{eqSiC.7}), some more analysis is needed. Our
next claim is that
\[
\lim_{n\rightarrow\infty}\left(  M^{n}\pi_{Z}\left(  \beta\right)
H_{0}\right)  \left(  z,x\right)  =\beta\left(  1\right)  H_{0}\left(
z,x\right)
\]
whenever $\beta\in C\left(  \mathbb{T}\right)  $, and the limit is in
$\mathcal{H}_{Z}$, or in $L^{2}\left(  \mathbb{R}\right)  $ after a
translation of the result via the Zak transform. Using Lemma \ref{LemPoof.1},
we get
\[
M^{n}\pi_{Z}\left(  \beta\right)  H_{0}\left(  z,\,\cdot\,\right)
=2^{-\frac{n}{2}}\sum_{w^{2^{n\mathstrut}}=z}m_{0}\left(  w\right)
m_{0}\left(  w^{2}\right)  \cdots m_{0}\left(  w^{2^{n-1}}\right)
\beta\left(  w\right)  H_{0}\left(  w,2^{n}x\right)  .
\]
Let $m_{0}^{\left(  n\right)  }\left(  w\right)  :=m_{0}\left(  w\right)
m_{0}\left(  w^{2}\right)  \cdots m_{0}\left(  w^{2^{n-1}}\right)  $. Since
$M\left(  H_{0}\right)  =H_{0}$, we get for the difference (picking $\beta$
such that $\beta\left(  1\right)  =1$)%
\[
H_{0}-M^{n}\pi_{Z}\left(  \beta\right)  H_{0}=2^{-\frac{n}{2}}\sum
_{w^{2^{n\mathstrut}}=z}m_{0}^{\left(  n\right)  }\left(  w\right)
\cdot\left(  1-\beta\left(  w\right)  \right)  H\left(  w,2^{n}x\right)  .
\]
The $\mathcal{H}_{Z}$-norm is that of $L^{2}\left(  \mathbb{T}\times\left[
0,1\right]  \right)  $; we split up the integral $\frac{1}{2\pi}\int_{-\pi
}^{\pi}\cdots\,d\omega$ over $\mathbb{T}$ into two regions, one $-\delta$
$\leq\omega\leq\delta$, and the other the union of the intervals $-\pi
\leq\omega\leq-\delta$ and $\delta\leq\omega\leq\pi$. We pick $\delta$ such
that $\left|  1-\beta\left(  e^{-i\omega}\right)  \right|  \leq\varepsilon$
when $\left|  \omega\right|  \leq\delta$. After estimating the two separate
contributions, and using Lemma \ref{LemSiC.1} and Corollary \ref{CorZak.5},
the desired result follows. When it is applied to the right-hand side in
(\ref{eqSiC.7}), we get%
\[
p\left(  H_{0},M^{n}\pi_{Z}\left(  \beta_{c}\right)  H_{0}\right)  \left(
z\right)  \underset{n\rightarrow\infty}{\longrightarrow}\beta_{c}\left(
1\right)  p_{2}\left(  H_{0}\right)  \left(  z\right)  =p_{2}\left(
H_{0}\right)  \left(  z\right)
\]
in the $L^{1}\left(  \mathbb{T}\right)  $-norm. But $\beta_{c}\left(
1\right)  =1$. Comparing the two results for the limits of the left- and
right-hand sides of (\ref{eqSiC.7}), we arrive at the desired contradiction,
and conclude that the eigenvalue problem (\ref{eqSiC.6}) does not have
eigenvectors as stated.
\end{proof}

\section[Singular vectors]{\label{SiV}SINGULAR VECTORS}

We now return to the subspace $\mathcal{B}$ of the Wold decomposition in
Theorem \ref{ThmR-is.2}. But we will specialize to the Hilbert space
$\mathcal{H}=\mathcal{H}_{Z}$, although the operators on $\mathcal{H}_{Z}$
correspond to (unitarily equivalent) versions, on $L^{2}\left(  \mathbb{R}%
\right)  $ and $\widehat{L^{2}\left(  \mathbb{R}\right)  }$, via the inverse
Zak transform $Z^{\ast}$ and the Fourier transform, respectively.

The object is to understand when the singular space $\mathcal{B}$ is present
in the decomposition $\mathcal{H}_{Z}=\mathcal{B}\oplus\mathcal{H}^{\left(
\infty\right)  }$. This is important, as we showed in the previous section
that the cascade approximation picks up divergences when $\mathcal{B}\neq0$.
Since
\begin{equation}
\mathcal{B}=\sideset{}{^{\smash{\oplus}}}{\sum}\limits_{n=0}^{\infty}\left[
M^{n}\mathcal{L}\right]  , \label{eqSiV.1}%
\end{equation}
where $\mathcal{L}=\ker\left(  M^{\ast}\right)  $, we see that the question is
decided by the question of when $\mathcal{L}\neq0$. Hence we need to study
$M^{\ast}$ more closely, and the second relation (\ref{ThmR-is.2(2)}) from
Theorem \ref{ThmR-is.2} is crucial for that.
While there is some connection between the results in this section and those
of \cite{DaLa}, there are several differences as well: the present approach is
general and applies equally well to higher dimensions, matrix dilations,
general $N$-to-$1$ maps (onto) in metric spaces, and, for the wavelets, even
to the case when the Hilbert space is different from the standard one, i.e.,
$L^{2}\left(  \mathbb{R}\right)  $.

Let $h\in L^{2}\left(  \mathbb{R}\right)  $, and set $\hat{h}\left(
\omega\right)  =\int_{\mathbb{R}}e^{-i\omega x}h\left(  x\right)  \,dx$. The
inverse Fourier transform will be denoted
\[
f\spcheck\left(  x\right)  =\frac{1}{2\pi}\int_{\mathbb{R}}e^{i\omega
x}f\left(  \omega\right)  \,d\omega.
\]
If $\mathcal{S}\subset L^{2}\left(  \mathbb{R}\right)  $ is a closed subspace,
we set%
\begin{equation}
\mathcal{S}\spcheck=\left\{  f\spcheck\mid f\in\mathcal{S}\right\}
.\label{eqSiV.5}%
\end{equation}
If $E\subset\mathbb{R}$ is a measurable subset, consider the projection
$P_{E}$ given by
\[
\left(  P_{E}f\right)  \left(  \omega\right)  =\chi_{E}^{{}}\left(
\omega\right)  f\left(  \omega\right)  ,\qquad\omega\in\mathbb{R},
\]
and define%
\begin{equation}
L^{2}\left(  E\right)  :=P_{E}\left(  L^{2}\left(  \mathbb{R}\right)  \right)
.\label{eqSiV.6}%
\end{equation}

We now turn to the closed subspaces in the decomposition (\ref{eqR-is.6}).
Considering the $L^{2}\left(  \mathbb{R}\right)  $ model, we show that there
are measurable subsets $E_{k}\left(  m_{0}\right)  \subset\mathbb{R}$ such
that%
\begin{equation}
\left[  M^{k}\mathcal{L}\right]  =L^{2}\left(  E_{k}\left(  m_{0}\right)
\right)  \spcheck\text{,\qquad the }L^{2}\left(  \mathbb{R}\right)  \text{
picture, }k=0,1,\dots,\label{eqSiV.2}%
\end{equation}
where $\left(  \,\cdot\,\right)  \spcheck$ denotes the inverse Fourier
transform in $L^{2}\left(  \mathbb{R}\right)  $. In fact it follows from
conclusion (\ref{ThmR-is.2(2)}) of Theorem \ref{ThmR-is.2}, i.e., our
Wold-type decomposition theorem (Section \ref{R-is}) that the spaces $\left[
M^{k}\mathcal{L}\right]  \subset L^{2}\left(  \mathbb{R}\right)  $ must have
the stated form (\ref{eqSiV.2}), but the object in the present section is to
find the sets. To see this, recall the formula%
\begin{equation}
\left(  \pi\left(  \xi\right)  h\right)  \sphat\left(  \omega\right)
=\xi\left(  e^{-i\omega}\right)  \hat{h}\left(  \omega\right)  ,\qquad\xi\in
L^{\infty}\left(  \mathbb{T}\right)  ,\;h\in L^{2}\left(  \mathbb{R}\right)
,\;\omega\in\mathbb{R},\label{eqSiV.3}%
\end{equation}
for the $L^{\infty}\left(  \mathbb{T}\right)  $-representation $\pi$. Using
Theorem \ref{ThmR-is.2} a second time, we conclude that there is also a
measurable $E_{\infty}\left(  m_{0}\right)  \subset\mathbb{R}$ such that
$\mathcal{H}^{\left(  \infty\right)  }=\left(  L^{2}\left(  E_{\infty}\left(
m_{0}\right)  \right)  \right)  \spcheck$. Using finally the orthogonality
part of the conclusion in Theorem \ref{ThmR-is.2}, we note that%
\begin{equation}
\left(  \bigcup_{k=0}^{\infty}E_{k}\left(  m_{0}\right)  \right)  \cup
E_{\infty}\left(  m_{0}\right)  =\mathbb{R},\label{eqSiV.4}%
\end{equation}
which is to say that they form a tiling of $\mathbb{R}$. Since the spaces in
the decomposition are mutually orthogonal, the sets $E_{k}\left(
m_{0}\right)  $, $k=0,1,\dots,k=\infty$, must be pairwise non-overlapping up
to measure zero in $\mathbb{R}$, i.e.,%
\[
E_{k}\left(  m_{0}\right)  \cap E_{l}\left(  m_{0}\right)
\]
has Lebesgue measure zero when $k\neq l$.

When the measurable subset $E\subset\mathbb{R}$ is given, we use the notation
$2E$ for%
\begin{equation}
2E=\left\{  2\omega\mid\omega\in E\right\}  ,\label{eqSiV.7}%
\end{equation}
and similarly, for $a\in\mathbb{R}$, set%
\begin{equation}
E+a=\left\{  \omega+a\mid\omega\in E\right\}  .\label{eqSiV.8}%
\end{equation}

The mapping $\omega\mapsto e^{-i\omega}$ is a measurable bijection of $\left[
-\pi,\pi\right\rangle $ onto $\mathbb{T}$. Let $c_{0}$ be its inverse (also
measurable). If $\log$ denotes the corresponding principal branch of the
complex logarithm, then%
\begin{equation}
c_{0}\left(  z\right)  =i\log\left(  z\right)  ,\qquad z\in\mathbb{T}.
\label{eqSiV.9}%
\end{equation}

Let $m_{0}$ be a subband filter satisfying conditions (\ref{IntrAxiom(1)}%
)--(\ref{IntrAxiom(3)}) in the Introduction, and let
\begin{equation}
N\left(  m_{0}\right)  :=\left\{  z\in\mathbb{T}\mid m_{0}\left(  z\right)
=0\right\}  . \label{eqSiV.10}%
\end{equation}
Let $M$ denote the corresponding cascade operator, i.e.,%
\begin{equation}
\left(  Mh\right)  \left(  x\right)  =\sqrt{2}\sum_{n\in\mathbb{Z}}%
a_{n}h\left(  2x-n\right)  ,\qquad h\in L^{2}\left(  \mathbb{R}\right)  ,
\label{eqSiV.11}%
\end{equation}
where $m_{0}\left(  z\right)  =\sum_{n\in\mathbb{Z}}a_{n}z^{n}$. Finally, let
$E\left(  m_{0}\right)  $ be given by%
\begin{equation}
E\left(  m_{0}\right)  =\bigcup_{n\in\mathbb{Z}}2c_{0}\left(  N\left(
m_{0}\right)  \right)  +4\pi n. \label{eqSiV.12}%
\end{equation}

\begin{lemma}
\label{LemSiV.1}%
\[
\ker\left(  M^{\ast}\right)  =L^{2}\left(  E\left(  m_{0}\right)  \right)
\spcheck=\left\{  f\spcheck\mid f\in L^{2}\left(  E\left(  m_{0}\right)
\right)  \right\}  .
\]
\end{lemma}

\begin{proof}
We have the formula%
\begin{equation}
\left(  M^{\ast}h\right)  \sphat\left(  \omega\right)  =\overline{m_{0}\left(
e^{-i\omega}\right)  }\cdot\hat{h}\left(  2\omega\right)  ,\qquad\omega
\in\mathbb{R}, \label{eqSiV.13}%
\end{equation}
directly from a Fourier transform of $M^{\ast}$. It follows that $\hat
{h}\left(  2\omega\right)  =0$ if $z=e^{-i\omega}\in\mathbb{T}\setminus
N\left(  m_{0}\right)  $. The lemma now follows from (\ref{eqSiV.12}) and
(\ref{eqSiV.13}).
\end{proof}

An immediate consequence of this lemma and (\ref{eqSiV.13}) is the following result:

\begin{proposition}
\label{ProSiV.2}We have the equivalence%
\[
\ker\left(  M^{\ast}\right)  =\left\{  0\right\}  \Longleftrightarrow
\mu\left(  N\left(  m_{0}\right)  \right)  =0,
\]
where $\mu$ denotes the Haar measure on $\mathbb{T}$.
\end{proposition}

We shall need the following second consequence of the lemma:

\begin{proposition}
\label{ProSiV.3}If $h\in\ker\left(  M^{\ast}\right)  $, then $\hat{h}\equiv0$
on $E\left(  m_{0}\right)  +2\pi$.
\end{proposition}

\begin{proof}
We have%
\[
E\left(  m_{0}\right)  +2\pi=\bigcup_{n\in\mathbb{Z}}2\cdot\left(
c_{0}\left(  N\left(  m_{0}\right)  \right)  +\pi\right)  +4\pi n=\bigcup
_{n\in\mathbb{Z}}2\cdot c_{0}\left(  -N\left(  m_{0}\right)  \right)  +4\pi
n.
\]
But if $z\in N\left(  m_{0}\right)  $, then $\left|  m_{0}\left(  -z\right)
\right|  =\sqrt{2}$ (from property (\ref{IntrAxiom(3)}) of $m_{0}$ in the
Introduction), so $-z\in\mathbb{T}\setminus N\left(  m_{0}\right)  $. Formula
(\ref{eqSiV.12}) then implies that $\hat{h}\left(  \omega+2\pi\right)  =0$ if
$\omega\in E\left(  m_{0}\right)  $, which is the desired conclusion.
\end{proof}

We now turn to $\ker\left(  M^{\ast\,k}\right)  $, $k>1$. We find measurable
sets $F_{k}\left(  m_{0}\right)  \subset\mathbb{R}$ such that%
\begin{equation}
\ker\left(  M^{\ast\,k}\right)  =\left(  L^{2}\left(  F_{k}\left(
m_{0}\right)  \right)  \right)  \spcheck. \label{eqSiV.14}%
\end{equation}
Let $m_{0}^{\left(  k\right)  }\left(  z\right)  :=m_{0}\left(  z\right)
m_{0}\left(  z^{2}\right)  \cdots m_{0}\left(  z^{2^{k-1}}\right)  $. Then%
\begin{align}
N\left(  m_{0}^{\left(  k\right)  }\right)   &  =\left\{  z\in\mathbb{T}%
\bigm| m_{0}^{\left(  k\right)  }\left(  z\right)  =0\right\}
\label{eqSiV.15}\\
&  =N\left(  m_{0}\right)  \cup\left\{  z\bigm| z^{2}\in N\left(
m_{0}\right)  \right\}  \cup\dots\cup\left\{  z\bigm| z^{2^{k-1}}\in N\left(
m_{0}\right)  \right\}  .\nonumber
\end{align}
By Claim \ref{Cla4} in the proof of Theorem \ref{ThmR-is.2}, we have%
\begin{equation}
\ker\left(  M^{\ast\,k+1}\right)  \ominus\ker\left(  M^{\ast\,k}\right)
=\left[  M^{k}\mathcal{L}\right]  , \label{eqSiV.16}%
\end{equation}
so we will get the components $\left[  M^{k}\mathcal{L}\right]  $ from the
following lemma.

\begin{lemma}
\label{LemSiV.4}The sets $F_{k}\left(  m_{0}\right)  $ from formula
\textup{(\ref{eqSiV.14})} for $\ker\left(  M^{\ast\,k}\right)  $ are
\begin{equation}
F_{k}\left(  m_{0}\right)  =\bigcup_{n\in\mathbb{Z}}2^{k}\cdot c_{0}\left(
N\left(  m_{0}^{\left(  k\right)  }\right)  \right)  +2^{k+1}\cdot\pi n.
\label{eqSiV.17}%
\end{equation}
\end{lemma}

\begin{proof}
This is the same argument as the one used in Lemma \ref{LemSiV.1} above, and
it is based on%
\[
\left(  M^{\ast\,k}h\right)  \sphat\left(  \omega\right)  =\overline
{m_{0}^{\left(  k\right)  }\left(  e^{-i\omega}\right)  }\cdot\hat{h}\left(
2^{k}\cdot\omega\right)  ,\qquad h\in L^{2}\left(  \mathbb{R}\right)
,\;\omega\in\mathbb{R}.\settowidth{\qedskip}{$\displaystyle\left(  M^{\ast
\,k}h\right)  \sphat\left(  \omega\right)  =\overline
{m_{0}^{\left(  k\right)  }\left(  e^{-i\omega}\right)  }\cdot\hat{h}\left(
2^{k}\cdot\omega\right)  ,\qquad h\in L^{2}\left(  \mathbb{R}\right)
,\;\omega\in\mathbb{R}.$}\addtolength{\qedskip}{-\textwidth}\rlap{\hbox
to-0.5\qedskip{\hfil\qed}}%
\]
\renewcommand{\qed}{}
\end{proof}

\begin{figure}[ptb]
\mbox{\psfig
{figure=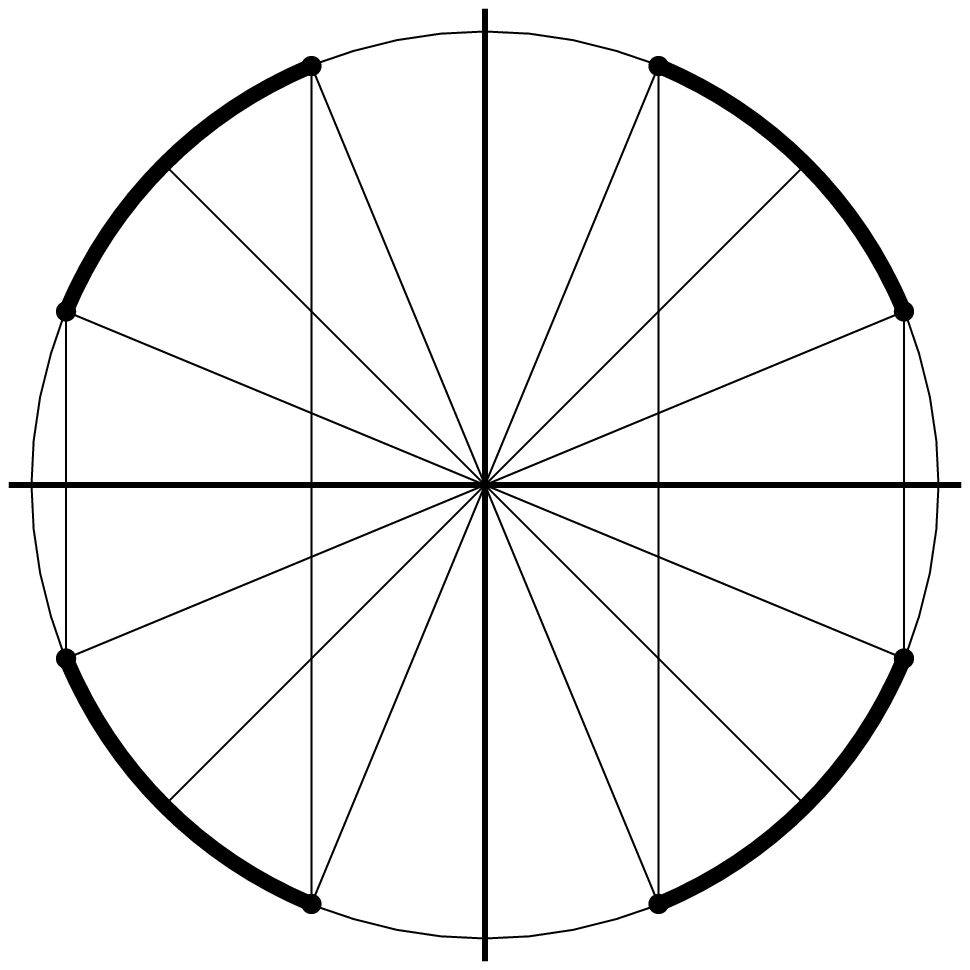,bbllx=0bp,bblly=0bp,bburx=288bp,bbury=288bp,width=288pt}%
}\llap{\setlength{\unitlength}{144pt}\begin{picture}(2,2)(-1,-1)
\put(0.9239,0.3827){\makebox(0,0)[l]{$\frac{\pi}{8}$}}
\put(0.3827,0.9239){\makebox(0,0)[b]{$\frac{3\pi}{8}$}}
\put(-0.3827,0.9239){\makebox(0,0)[b]{$\frac{5\pi}{8}$}}
\put(-0.9239,0.3827){\makebox(0,0)[r]{$\frac{7\pi}{8}$}}
\put(-0.9239,-0.3827){\makebox(0,0)[r]{$-\frac{7\pi}{8}$}}
\put(-0.3827,-0.9239){\makebox(0,0)[t]{$-\frac{5\pi}{8}$}}
\put(0.3827,-0.9239){\makebox(0,0)[t]{$-\frac{3\pi}{8}$}}
\put(0.9239,-0.3827){\makebox(0,0)[l]{$-\frac{\pi}{8}$}}
\end{picture}}\caption{$\sigma^{-2}\left(  N\left(  m_{0}\right)  \right)  $}%
\label{FigFourArcs}%
\end{figure}

In using formula (\ref{eqSiV.14}), the following observation on $N\left(
m_{0}^{\left(  k\right)  }\right)  $ is useful. Let $\sigma\left(  z\right)
=z^{2}$ be the square map of $\mathbb{T}$, and let $N\subset\mathbb{T}$ be a
subset. Let
\[
\sigma^{-1}\left(  N\right)  =\left\{  z\in\mathbb{T}\mid\sigma\left(
z\right)  \in N\right\}  .
\]
To understand the dynamical picture of $\sigma$ and its multivalued inverse
$\sigma^{-1}$, i.e., their iterations, we have included a graphical
illustration in Figures \ref{FigFourArcs} and \ref{FigErgoMap}. Then
\[
N\left(  m_{0}^{\left(  n\right)  }\right)  =N\left(  m_{0}\right)  \cup
\sigma^{-1}\left(  N\left(  m_{0}\right)  \right)  \cup\dots\cup
\sigma^{-\left(  n-1\right)  }\left(  N\left(  m_{0}\right)  \right)
\]
and
\[
N\left(  m_{0}^{\left(  n+1\right)  }\right)  \bigm\backslash N\left(
m_{0}^{\left(  n\right)  }\right)  =\sigma^{-n}\left(  N\left(  m_{0}\right)
\right)  =\left\{  z\in\mathbb{T}\bigm|z^{2^{n}}\in N\left(  m_{0}\right)
\right\}  .
\]

For the example which follows (Example \ref{ExaSiV.6} below), we will have
\begin{align*}
N\left(  m_{0}\right)   &  =\left\{  e^{-i\omega}\bigm|\frac{\pi}{2}%
\leq\left|  \omega\right|  \leq\pi\right\}  ,\\
\sigma^{-1}\left(  N\left(  m_{0}\right)  \right)   &  =\left\{  e^{-i\omega
}\bigm|\frac{\pi}{4}\leq\left|  \omega\right|  \leq\frac{3\pi}{4}\right\}
\text{,\qquad and}\\
\sigma^{-2}\left(  N\left(  m_{0}\right)  \right)   &  =\left\{  e^{-i\omega
}\bigm|\frac{\pi}{8}\leq\left|  \omega\right|  \leq\frac{3\pi}{8}\right\}
\cup\left\{  e^{-i\omega}\bigm|\frac{5\pi}{8}\leq\left|  \omega\right|
\leq\frac{7\pi}{8}\right\}  ,
\end{align*}
the last made up of four arc segments on the unit circle (see Figure
\ref{FigFourArcs}); and the induction is clear for the general case
$\sigma^{-n}\left(  N\left(  m_{0}\right)  \right)  $.

When Lemma \ref{LemSiV.4} is combined with formula (\ref{eqSiV.16}), we get%
\begin{equation}
\left[  M^{k}\mathcal{L}\right]  =\left(  L^{2}\left(  F_{k+1}\left(
m_{0}\right)  \setminus F_{k}\left(  m_{0}\right)  \right)  \right)  \spcheck,
\label{eqSiV.18}%
\end{equation}
again with $\left(  \,\cdot\,\right)  \spcheck$ denoting inverse Fourier
transform. We now combine the results above into a proposition.

\begin{proposition}
\label{ProSiV.5}We have
\begin{equation}
\left[  M^{k}\mathcal{L}\right]  =L^{2}\left(  E_{k}\left(  m_{0}\right)
\right)  \spcheck, \label{eqSiV.19}%
\end{equation}
where%
\[
E_{k}\left(  m_{0}\right)  =F_{k+1}\left(  m_{0}\right)  \setminus
F_{k}\left(  m_{0}\right)  ,
\]
starting with%
\[
E_{0}\left(  m_{0}\right)  =F_{1}\left(  m_{0}\right)  =E\left(  m_{0}\right)
=\bigcup_{n\in\mathbb{Z}}2c_{0}\left(  N\left(  m_{0}\right)  \right)  +4\pi
n.
\]
\end{proposition}

\begin{proof}
Contained in the previous argument.
\end{proof}

\begin{example}
\label{ExaSiV.6}\begin{figure}[ptb]
\vspace*{12pt}\mbox{\hspace{6pt}\psfig
{figure=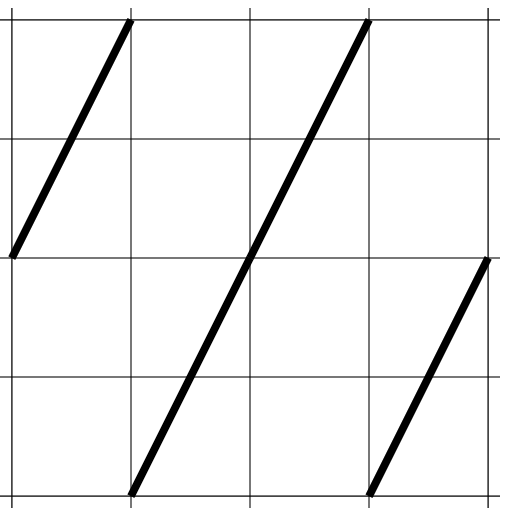,bbllx=0bp,bblly=0bp,bburx=144bp,bbury=144bp,width=144pt}%
\llap{\setlength{\unitlength}{72pt}\begin{picture}(2,2)(-1,-1)
\put(1,0){\makebox(0,0)[lb]{$\phantom{(}\pi$}}
\put(0,1){\makebox(0,0)[rb]{$\pi\phantom{)}$}}
\put(-1,0){\makebox(0,0)[rb]{$-\pi\phantom{)}$}}
\put(0,-1){\makebox(0,0)[rt]{$-\pi\phantom{)}$}}
\put(0,0){\makebox(0,0)[rb]{$0\phantom{)}$}}
\end{picture}}\hspace{42pt}\psfig
{figure=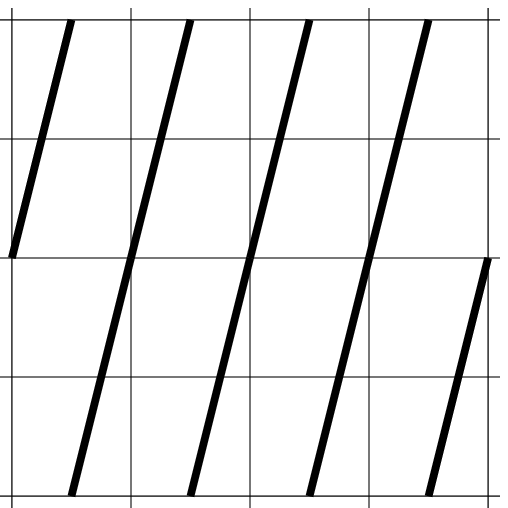,bbllx=0bp,bblly=0bp,bburx=144bp,bbury=144bp,width=144pt}%
\llap{\setlength{\unitlength}{72pt}\begin{picture}(2,2)(-1,-1)
\put(1,0){\makebox(0,0)[lb]{$\phantom{(}\pi$}}
\put(0,1){\makebox(0,0)[rb]{$\pi\phantom{)}$}}
\put(-1,0){\makebox(0,0)[rb]{$-\pi\phantom{)}$}}
\put(0,-1){\makebox(0,0)[rt]{$-\pi\phantom{)}$}}
\put(0,0){\makebox(0,0)[rb]{$0\phantom{)}$}}
\end{picture}}}\\[12pt]\mbox{\hspace{6pt}\makebox[144pt]{$\sigma\colon
z\mapsto z^{2}$}\hspace{42pt}\makebox[144pt]{$\sigma^{2}\colon z\mapsto
z^{4}$}}\caption{An ergodic map on $\mathbb{T}$}%
\label{FigErgoMap}%
\end{figure}

It is known that the function $\varphi\in L^{2}\left(  \mathbb{R}\right)  $
given by%
\begin{equation}
\varphi\left(  x\right)  =\left.  \chi_{\left[  -\pi,\pi\right]  }^{{}%
}\right.  \spcheck\left(  x\right)  =\frac{\sin\pi x}{\pi x}\label{eqSiV.20}%
\end{equation}
is a scaling function. The subband filter $m_{0}$ is%
\[
m_{0}\left(  e^{-i\omega}\right)  =\sqrt{2}\chi_{\left[  -\frac{\pi}{2}%
,\frac{\pi}{2}\right]  }^{{}}\left(  \omega\right)  .
\]
Hence%
\[
c_{0}\left(  N\left(  m_{0}\right)  \right)  =\left[  -\pi,-\frac{\pi}%
{2}\right\rangle \cup\left[  \frac{\pi}{2},\pi\right\rangle
\]
and%
\[
\mathcal{L}=\ker\left(  M^{\ast}\right)  =L^{2}\left(  E\left(  m_{0}\right)
\right)  \spcheck,
\]
where
\begin{equation}
E\left(  m_{0}\right)  =\bigcup_{n\in\mathbb{Z}}\left[  \pi,3\pi\right\rangle
+4\pi n.\label{eqSiV.21}%
\end{equation}
Writing $m_{0}\left(  \frac{\omega}{2}\right)  =m_{0}\left(  e^{-i\frac
{\omega}{2}}\right)  $, we see that
\[
m_{0}\left(  \frac{\omega}{2}+\pi\right) =\sqrt{2}\chi_{E\left(
m_{0}\right)  }^{{}}\left(  \omega\right) 
\text{\qquad and\qquad}
m_{0}\left(  \frac{\omega}{2}\right) =\sqrt{2}\chi_{G\left(
m_{0}\right)  }^{{}}\left(  \omega\right)  ,
\]
where $G\left(  m_{0}\right)  =\bigcup_{n\in\mathbb{Z}}\left[  -\pi
,\pi\right\rangle +4\pi n$. We have $E\left(  m_{0}\right)  \cup G\left(
m_{0}\right)  =\mathbb{R}$, and
\[
G\left(  m_{0}\right)  =E\left(  m_{0}\right)  +2\pi,
\]
which yields the following formula for the Ruelle operator:%
\begin{equation}
\left(  R\xi\right)  \left(  \omega\right)  =\chi_{G\left(  m_{0}\right)
}^{{}}\left(  \omega\right)  \xi\left(  \frac{\omega}{2}\right)
+\chi_{E\left(  m_{0}\right)  }^{{}}\left(  \omega\right)  \xi\left(
\frac{\omega}{2}+\pi\right)  ,\qquad\xi\in L^{\infty}\left(  \mathbb{T}%
\right)  ,\;\omega\in\mathbb{R},\label{eqSiV.22}%
\end{equation}
where we use $\mathbb{T}\simeq\mathbb{R}\diagup2\pi\mathbb{Z}$. It is easy to
check directly from \textup{(\ref{eqSiV.22})} that
\[
\left(  R\xi\right)  \left(  \omega+2\pi\right)  =\left(  R\xi\right)  \left(
\omega\right)  \text{,\qquad for all }\omega\in\mathbb{R}.
\]
Similarly the scaling identity for $\varphi$, i.e.,
\[
\sqrt{2}\hat{\varphi}\left(  \omega\right)  =m_{0}\left(  e^{-i\frac{\omega
}{2}}\right)  \hat{\varphi}\left(  \frac{\omega}{2}\right)  ,
\]
takes the form%
\[
\hat{\varphi}\left(  \omega\right)  =\chi_{G\left(  m_{0}\right)  }^{{}%
}\left(  \omega\right)  \hat{\varphi}\left(  \frac{\omega}{2}\right)
,\qquad\omega\in\mathbb{R},
\]
which in addition to $\hat{\varphi}_{1}=\chi_{\left[  -\pi,\pi\right\rangle
}^{{}}$ has the solution $\hat{\varphi}_{2}=\chi_{G\left(  m_{0}\right)  }%
^{{}}$, but, of course, $\chi_{G\left(  m_{0}\right)  }^{{}}$ is not in
$L^{2}\left(  \mathbb{R}\right)  $. Nonetheless, the inverse transform,
$\varphi_{2}=\left.  \chi_{G\left(  m_{0}\right)  }^{{}}\right.  \spcheck$,
makes sense as a distribution. The corresponding quadratic expression,
$p_{2}\left(  \varphi\right)  =p\left(  \varphi,\varphi\right)  $, satisfies
\begin{align*}
p_{2}\left(  \varphi_{1}\right)  \left(  e^{-i\omega}\right)   &  \equiv1,\\%
\intertext{while}%
p_{2}\left(  \varphi_{2}\right)  \left(  e^{-i\omega}\right)   &  =\sum
_{n\in\mathbb{Z}}\chi_{G\left(  m_{0}\right)  }^{{}}\left(  \omega+2\pi
n\right)  \equiv\infty.
\end{align*}
Since
\[
m_{0}^{\left(  k\right)  }\left(  \omega\right)  =2^{\frac{k}{2}}\chi_{\left[
-\frac{\pi}{2},\frac{\pi}{2}\right\rangle }^{{}}\left(  \left[  \omega\right]
\right)  \chi_{\left[  -\frac{\pi}{2},\frac{\pi}{2}\right\rangle }^{{}}\left(
\left[  2\omega\right]  \right)  \cdots\chi_{\left[  -\frac{\pi}{2},\frac{\pi
}{2}\right\rangle }^{{}}\left(  \left[  2^{k-1}\omega\right]  \right)
,\qquad\omega\in\mathbb{R},
\]
where $\left[  \omega\right]  :=c_{0}\left(  e^{-i\omega}\right)  $, we see
that the higher cases result from an iteration of an ergodic map on
$\mathbb{T}$ as follows: $c_{0}\left(  z\right)  \mapsto c_{0}\left(
z^{2}\right)  $, or $\omega\mapsto2\omega\operatorname{mod}2\pi$. Starting
with
\begin{align*}
c_{0}\left(  N\left(  m_{0}\right)  \right)   &  =\left[  -\pi,-\frac{\pi}%
{2}\right\rangle \cup\left[  \frac{\pi}{2},\pi\right\rangle ,\\%
\intertext{we get}%
c_{0}\left(  N\left(  m_{0}^{\left(  2\right)  }\right)  \right)   &  =\left[
-\pi,-\frac{\pi}{4}\right\rangle \cup\left[  \frac{\pi}{4},\pi\right\rangle
,\\%
\intertext{and by induction}%
c_{0}\left(  N\left(  m_{0}^{\left(  k\right)  }\right)  \right)   &  =\left[
-\pi,-\frac{\pi}{2^{k\mathstrut}}\right\rangle \cup\left[  \frac{\pi
}{2^{k\mathstrut}},\pi\right\rangle
\end{align*}
\textup{(}see \textup{(\ref{eqSiV.9}) and Figure \ref{FigErgoMap}} for the map
$c_{0}$\textup{),} and%
\begin{align*}
F_{k}\left(  m_{0}^{\left(  k\right)  }\right)   &  =\bigcup_{n\in\mathbb{Z}%
}\left[  -2^{k}\pi,-\pi\right\rangle \cup\left[  \pi,2^{k}\pi\right\rangle
+2^{k+1}\pi n,\\
\left[  M^{k}\mathcal{L}\right]   &  =L^{2}\left(  E_{k}\left(  m_{0}\right)
\right)  \spcheck,\\%
\intertext{where}%
E_{k}\left(  m_{0}\right)   &  =F_{k+1}\left(  m_{0}\right)  \setminus
F_{k}\left(  m_{0}\right)  ,
\end{align*}
starting with%
\begin{align*}
E_{0}\left(  m_{0}\right)   &  =E\left(  m_{0}\right)  =F_{1}\left(
m_{0}\right)  =\bigcup_{n\in\mathbb{Z}}\left[  \pi,3\pi\right\rangle +4\pi n\\%
\intertext{and}%
E_{1}\left(  m_{0}\right)   &  =F_{2}\left(  m_{0}\right)  \setminus
F_{1}\left(  m_{0}\right)  =\bigcup_{n\in\mathbb{Z}}\left[  2\pi
,6\pi\right\rangle +8\pi n.
\end{align*}
\end{example}

\begin{remark}
\label{RemOnExaSiV.6}Example \textup{\ref{ExaSiV.6}} is a special case of a
classification from \cite{BrJo97b} of filters $m_{0}$ for which $\left|
m_{0}\right|  $ takes on only the two values $\sqrt{2}$ and $0$. However, the
present analysis is focused on convergence questions which are not addressed there.

There is one more feature which sets this class of examples apart from those
where $\mathcal{L}=\ker\left(  M^{\ast}\right)  =0$, i.e., those for which the
complement of the support of $m_{0}$ in $\mathbb{T}$ has measure zero
\textup{(}see Proposition \textup{\ref{ProSiV.2}).} When $\ker\left(  M^{\ast
}\right)  \neq0$, then%
\begin{equation}
\mathcal{E}_{1}=\left\{  \xi\in L^{\infty}\left(  \mathbb{T}\right)  \mid
R\left(  \xi\right)  =\xi\right\}  \label{eqSiV.eigenspace}%
\end{equation}
may be infinite-dimensional.

In Example \textup{\ref{ExaSiV.6},} recall $\mathcal{H}^{\left(
\infty\right)  }=L^{2}\left(  -\pi,\pi\right)  \spcheck$, and%
\begin{align}
p_{2}\left(  f\right)  \left(  e^{-i\omega}\right)   &  =\chi_{\left[
-\pi,\pi\right\rangle }^{{}}\left(  \omega\right)  \left|  \hat{f}\left(
\omega\right)  \right|  ^{2}\label{eqSiV.quadratic}\\
&  =P_{\left[  -\pi,\pi\right\rangle }\left(  \left|  \smash{\hat{f}}\right|
^{2}\right)  \left(  \omega\right)  ,\qquad f\in\mathcal{H}^{\left(
\infty\right)  },\;\omega\in\mathbb{R}.\nonumber
\end{align}
\end{remark}

Since $M\varphi=\varphi$, it is clear that $\varphi\in\mathcal{H}^{\left(
\infty\right)  }$, where, in this case,%
\[
\mathcal{H}^{\left(  \infty\right)  }=\bigcap_{n=1}^{\infty}M^{n}\left(
L^{2}\left(  \mathbb{R}\right)  \right)  =L^{2}\left(  \mathbb{R}\right)
\ominus\left\{  M^{k}\mathcal{L}\mid k=0,1,\dots\right\}  ,
\]
where $\mathcal{L}:=\ker M^{\ast}$. Since $\mathcal{H}^{\left(  \infty\right)
}=L^{2}\left(  E_{\infty}\right)  \spcheck$, it follows that $\left[  -\pi
,\pi\right\rangle \subset E_{\infty}$, or equivalently $L^{2}\left(  -\pi
,\pi\right)  \spcheck\subset L^{2}\left(  E_{\infty}\right)  \spcheck$. But if
iteration of (\ref{eqSiV.13}) for $\left(  M^{\ast\,n}h\right)  \sphat$ is
combined with the above analysis of the example, we conclude that, in fact,
the inclusion is equality, i.e., $\left[  -\pi,\pi\right\rangle =E_{\infty}$,
up to Lebesgue measure zero in $\mathbb{R}$, and, therefore, $\mathcal{H}%
^{\left(  \infty\right)  }=L^{2}\left(  -\pi,\pi\right)  \spcheck$. The main
ingredient in this argument is the known ergodicity of the map $z\overset
{\sigma}{\longmapsto}z^{2}$ of $\mathbb{T}$; see Figure \ref{FigErgoMap}, and
\cite{Kea72} for details on $N$-to-$1$ endomorphisms of measure space.

We now consider the cascade approximation $\lim_{n\rightarrow\infty}M^{n}h$ in
$L^{2}\left(  \mathbb{R}\right)  $. The object is to pick $h$ such that the
limit is the scaling function $\varphi$. By Theorem \ref{ThmR-is.2}, the
optimal choice dictates $h\in\mathcal{H}^{\left(  \infty\right)  }$, so that
excludes (for the example) the standard choice for starting point of cascades
$h\rightarrow Mh\rightarrow M^{2}h$ $\rightarrow\cdots$, which is
$h=\chi_{\left[  0,1\right\rangle }^{{}}$. Clearly $\hat{h}\left(
\omega\right)  =\widehat{\chi_{\left[  0,1\right\rangle }^{{}}}\left(
\omega\right)  =e^{-i\frac{\omega}{2}}\frac{\sin\frac{\omega}{2\mathstrut}%
}{\frac{\omega}{2}}$ is not in $L^{2}\left(  -\pi,\pi\right)  $, so $\left\|
M^{n}\chi_{\left[  0,1\right\rangle }^{{}}-\varphi\right\|  _{L^{2}\left(
\mathbb{R}\right)  }\!\rightarrow0$ as $n\rightarrow\infty$ is at best a slow approximation.

Let $h=h_{\mathcal{B}}+h_{\infty}$ be the $\mathcal{B}\oplus\mathcal{H}%
^{\left(  \infty\right)  }$ decomposition of Theorem \ref{ThmR-is.2}. From
Corollary \ref{CorR-is.4} we have $M^{n}h_{\mathcal{B}}\in\mathcal{B}$ and
$M^{n}h_{\infty}\in\mathcal{H}^{\left(  \infty\right)  }$, so%
\[
M^{n}h-\varphi=\underset{\mathcal{H}^{\left(  \infty\right)  }}{\underbrace
{M^{n}h_{\infty}-\varphi}}+\underset{\mathcal{B}}{\underbrace{M^{n}%
h_{\mathcal{B}}}}%
\]
and%
\begin{equation}
\left\|  M^{n}h-\varphi\right\|  ^{2}=\left\|  M^{n}h_{\infty}-\varphi
\right\|  ^{2}+\left\|  M^{n}h_{\mathcal{B}}\right\|  ^{2}. \label{eqSiV.25}%
\end{equation}

\begin{proposition}
\label{ProSiV.8}Let $m_{0}(  e^{-i\omega})  =\sqrt{2}\chi^{{}%
}_{\left[  -\frac{\pi}{2},\frac{\pi}{2}\right\rangle }(  [
\omega]  )  $ as a function on $\mathbb{T}=\mathbb{R}\diagup
2\pi\mathbb{Z}$, and let $R$, $M$ be the respective Ruelle and cascade
operators. Let $\varphi=\left.  \chi^{{}}_{\left[  -\pi,\pi\right\rangle
}\right.  \spcheck$ and $h=\chi^{{}}_{\left[  0,1\right\rangle }$. Then
$M\varphi=\varphi$ and $h=h_{\mathcal{B}}+h_{\infty}$, where%
\[
\widehat{h_{\mathcal{B}}}\left(  \omega\right) =\chi^{{}}_{\mathbb{R}
\setminus\left[  -\pi,\pi\right\rangle }\left(  \omega\right)  e^{-i\frac
{\omega}{2}}\frac{\sin\left(  \frac{\omega}{2}\right)  }{\frac{\omega}{2}}
\text{\qquad and\qquad}
\widehat{h_{\infty}}\left(  \omega\right) =\chi^{{}}_{\left[  -\pi
,\pi\right\rangle }\left(  \omega\right)  e^{-i\frac{\omega}{2}}\frac
{\sin\left(  \frac{\omega}{2}\right)  }{\frac{\omega}{2}},
\]
with the approximations:

\begin{enumerate}
\item \label{ProSiV.8(1)}$\left\|  M^{n}h_{\infty}-\varphi\right\|
_{2}\underset{n\rightarrow\infty}{\longrightarrow}0$, and

\item \label{ProSiV.8(2)}$\left\|  M^{n}h_{\mathcal{B}}\right\|  _{2}%
\underset{n\rightarrow\infty}{\longrightarrow}0$.
\end{enumerate}
\end{proposition}

\begin{proof}
Note that, when (\ref{ProSiV.8(1)})--(\ref{ProSiV.8(2)}) are combined with
(\ref{eqSiV.25}), we get $\left\|  M^{n}h-\varphi\right\|  _{2}\underset
{n\rightarrow\infty}{\longrightarrow}0$, but (\ref{ProSiV.8(1)}) is a better approximation.

We first prove (\ref{ProSiV.8(1)}) by checking that\renewcommand{\theenumi
}{\alph{enumi}}

\begin{enumerate}
\item \label{ProSiV.8proof(1)}$\left\|  M^{n}h_{\infty}\right\|
_{2}\rightarrow1$, and

\item \label{ProSiV.8proof(2)}$\ip{M^{n}h_{\infty}}{\varphi}\rightarrow1$, as
$n\rightarrow\infty$.
\end{enumerate}

\noindent For the first term (\ref{ProSiV.8proof(1)}) we have%
\begin{equation}
\left\|  M^{n}h_{\infty}\right\|  _{2}^{2}=\int_{\mathbb{T}}R^{n}\left(
p_{2}\left(  h_{\infty}\right)  \right)  \left(  z\right)  \,d\mu\left(
z\right)  , \label{eqSiV.26}%
\end{equation}
where
\[
p_{2}\left(  h_{\infty}\right)  \left(  e^{-i\omega}\right)  =\sum
_{n\in\mathbb{Z}}\left|  \widehat{h_{\infty}}\left(  \omega+2\pi n\right)
\right|  ^{2}.
\]
We take $-\pi\leq\omega<\pi$, and recall the formula $\widehat{h_{\infty}%
}\left(  \omega\right)  =\chi_{\left[  -\pi,\pi\right\rangle }^{{}}\left(
\omega\right)  e^{-i\frac{\omega}{2}}\smash{\frac{\sin\left( \frac{\omega}%
{2}\right) }{\frac{\omega}{2}}}$. Hence $p_{2}\left(  h_{\infty}\right)
\left(  e^{-i\omega}\right)  =\left|  \frac{\sin\left(  \frac{\omega}%
{2}\right)  }{\frac{\omega}{2}}\right|  ^{2}$. Since this function is
continuous, it follows from a theorem of Meyer and Paiva \cite{MePa93} that%
\[
\frac{1}{2\pi}\int_{-\pi}^{\pi}R^{n}\left(  p_{2}\left(  h_{\infty}\right)
\right)  \left(  e^{-i\omega}\right)  \,d\omega\underset{n\rightarrow\infty
}{\longrightarrow}p_{2}\left(  h_{\infty}\right)  \left(  \omega=0\right)
=1.
\]
In view of (\ref{eqSiV.26}), this proves (\ref{ProSiV.8proof(1)}). The
application of the Meyer--Paiva theorem
(see \cite{MePa93} and line 3 in Table \ref{tableZak.2} above)
requires that $p_{2}\left(
\varphi\right)  \equiv1$ on $\mathbb{T}$, which clearly holds for the present
$\varphi$.

The argument for (\ref{ProSiV.8proof(2)}) is similar:
\[
\ip{M^{n}h_{\infty}}{\varphi}_{L^{2}\left(  \mathbb{R}\right)  }%
=\int_{\mathbb{T}}R^{n}\left(  p\left(  h_{\infty},\varphi\right)  \right)
\left(  z\right)  \,d\mu\left(  z\right)  \underset{n\rightarrow\infty
}{\longrightarrow}p\left(  h_{\infty},\varphi\right)  \left(  z=1\right)  =1.
\]
The last step is based on the formula ($\left|  \omega\right|  <\pi$)%
\begin{align*}
p\left(  h_{\infty},\varphi\right)  \left(  e^{-i\omega}\right)   &
=\smash{\sum_{n\in\mathbb{Z}}}\overline{\widehat{h_{\infty}}\left(
\omega+2\pi n\right)  }\hat{\varphi}\left(  \omega+2\pi n\right)  \\
&  =\left(  n=0\text{ term}\right)  =\hat{\varphi}\left(  \omega\right)
=e^{-i\frac{\omega}{2}}\frac{\sin\left(  \frac{\omega}{2}\right)  }%
{\frac{\omega}{2}}.
\end{align*}
It remains to prove (\ref{ProSiV.8(2)}). We have%
\[
\left\|  M^{n}h_{\mathcal{B}}\right\|  _{2}^{2}=\int_{\mathbb{T}}R^{n}\left(
p_{2}\left(  h_{\mathcal{B}}\right)  \right)  \left(  z\right)  \,d\mu\left(
z\right)  ,
\]
and we check that the theorem from \cite{MePa93} applies to that term as well:
by the above argument, we have (for $\left|  \omega\right|  <\pi$)%
\[
p_{2}\left(  h_{\mathcal{B}}\right)  \left(  e^{-i\omega}\right)  =\sum
_{n\in\mathbb{Z}\setminus\left\{  0\right\}  }\frac{\sin^{2}\left(
\frac{\omega}{2}\right)  }{\left(  \frac{\omega}{2}+n\pi\right)
^{2\mathstrut}}.
\]
By a standard summation formula (see, e.g., \cite[p.\ 152]{Car95}), we have%
\[
\sum_{n\in\mathbb{Z}\setminus\left\{  0\right\}  }\frac{1}{\left(
\frac{\omega}{2}+n\pi\right)  ^{2\mathstrut}}=\frac{1}{\sin^{2\mathstrut
}\left(  \frac{\omega}{2}\right)  }-\frac{1}{\left(  \frac{\omega}{2}\right)
^{2\mathstrut}}.
\]
Hence%
\[
p_{2}\left(  h_{\mathcal{B}}\right)  \left(  e^{-i\omega}\right)
=1-\frac{\sin^{2}\left(  \frac{\omega}{2}\right)  }{\left(  \frac{\omega}%
{2}\right)  ^{2\mathstrut}},
\]
which is continuous and equals $0$ at $\omega=0$, so (\ref{ProSiV.8(2)}) also
follows from an application of \cite{MePa93}.
\end{proof}

In the course of the proof, we established the following result.

\begin{proposition}
\label{ProSiV.9}Let $\varphi:=\left.  \chi_{\left[  -\pi,\pi\right\rangle
}^{{}}\right.  \spcheck$ \textup{(}inverse Fourier transform\textup{),} and
let $f\in L^{2}\left(  \mathbb{R}\right)  $ be given with $L^{2}$-Fourier
transform $\hat{f}$. Then we have the following a.e.\ identity on $\mathbb{T}%
$\textup{:}%
\[
p\left(  \varphi,f\right)  \left(  e^{-i\omega}\right)  =\chi_{\left[
-\pi,\pi\right\rangle }^{{}}\left(  \left[  \omega\right]  \right)  \hat
{f}\left(  \omega\right)  ,
\]
where $p$ is defined on $L^{2}\left(  \mathbb{R}\right)  \times L^{2}\left(
\mathbb{R}\right)  $ via the Zak transform,%
\[
p\left(  \varphi,f\right)  \left(  z\right)  =\int_{0}^{1}\overline{\left(
Z\varphi\right)  \left(  z,x\right)  }\left(  Zf\right)  \left(  z,x\right)
\,dx
\]
\textup{(}again for a.a.\ $z\in\mathbb{T}$\textup{).}
\end{proposition}

\section[Approximation results]{\label{Appr}APPROXIMATION RESULTS}

In this final section we prove a general approximation theorem which is based
on some of the same ideas which went into Example \ref{ExaSiV.6} in the
previous section. It uses the approximation kernel%
\begin{equation}
D_{n}\left(  z\right)  :=\left|  m_{0}\left(  z\right)  m_{0}\left(
z^{2}\right)  \cdots m_{0}\left(  z^{2^{n-1}}\right)  \right|  ^{2}%
=R^{\ast\,n}\left(  \openone\right)  \left(  z\right)  . \label{eqAppr.1}%
\end{equation}
This kernel may be constructed from any filter $m_{0}$ subject to conditions
(\ref{IntrAxiom(1)})--(\ref{IntrAxiom(3)}) in the Introduction. If the
corresponding scaling function $\varphi$ with $\hat{\varphi}\left(  0\right)
=1$ satisfies $p_{2}\left(  \varphi\right)  \equiv1$, then by a theorem of
Meyer and Paiva \cite{MePa93}, $D_{n}\left(  \,\cdot\,\right)  $ is an
approximate Dirac delta function $\delta_{1}$ on $C\left(  \mathbb{T}\right)
$. In Example \ref{ExaSiV.6} above, the kernel $D_{n}$ computes out as:%
\[
D_{n}\left(  e^{-i\omega}\right)  =2^{\frac{n}{2}}\chi_{\left[  -\frac{\pi}%
{2},\frac{\pi}{2}\right\rangle }^{{}}\left(  \omega\right)  \chi_{\left[
-\frac{\pi}{2},\frac{\pi}{2}\right\rangle }^{{}}\left(  \left[  2\omega
\right]  \right)  \cdots\chi_{\left[  -\frac{\pi}{2},\frac{\pi}{2}%
\right\rangle }^{{}}\left(  \left[  2^{n-1}\omega\right]  \right)  ,
\]
where, for $s\in\mathbb{R}$, we set $\left[  s\right]  =s-k2\pi$ if
$k\in\mathbb{Z}$ is such that $\left(  2k-1\right)  \pi\leq s<\left(
2k+1\right)  \pi$, or equivalently $\left[  s\right]  =c_{0}\left(
e^{-is}\right)  $, $s\in\mathbb{R}$.

\begin{definition}
\label{DefAppr.1}Let $\xi\in L^{2}\left(  \mathbb{T}\right)  $, and
$c\in\mathbb{C}$. While it is not then possible to assign a value to
$\xi\left(  1\right)  $, i.e., $\xi$ evaluated at $z=1$, we will say that
$E\xi\left(  1\right)  =c$ if
$\xi$ is continuous at $z=1$ and assumes the value $c$ there.
\end{definition}

\begin{theorem}
\label{ThmAppr.2}Let $m_{0}$ satisfy conditions \textup{(\ref{IntrAxiom(1)}%
)--(\ref{IntrAxiom(3)})} in the Introduction, and let $R$, $M$ be the
corresponding Ruelle and cascade operators, where we view $M$ as a operator in
$\mathcal{H}_{Z}\simeq L^{2}\left(  \mathbb{R}\right)  $. Let $H_{0}%
\in\mathcal{H}_{Z}$ be a regular scaling function, i.e., $\left\|
H_{0}\right\|  =1$, and $M\left(  H_{0}\right)  =H_{0}$. Let $F\in
\mathcal{H}_{Z}$ satisfy $p_{2}\left(  F\right)  =\openone$. Then the
following two conditions are equivalent:\renewcommand{\theenumi}{\alph{enumi}}

\begin{enumerate}
\item \label{ThmAppr.2(1)}$\left\|  H_{0}-M^{n}F\right\|  _{\mathcal{H}_{Z}%
}\underset{n\rightarrow\infty}{\longrightarrow}0$,

\item \label{ThmAppr.2(2)}$Ep\left(  H_{0},F\right)  \left(  1\right)  =1$,
\end{enumerate}

\noindent where the evaluation at $z=1$ in \textup{(\ref{ThmAppr.2(2)})} is
specified in Definition \textup{\ref{DefAppr.1}.}
\end{theorem}

\begin{proof}
Since $\left\|  H_{0}\right\|  =\openone$, it follows from Lemma
\ref{LemSiC.1} that $p_{2}\left(  H_{0}\right)  =\openone$, i.e., we have both
of the families $\left\{  \pi\left(  e_{n}\right)  H_{0}\mid n\in
\mathbb{Z}\right\}  $ and $\left\{  \pi\left(  e_{n}\right)  F\mid
n\in\mathbb{Z}\right\}  $ orthogonal in $\mathcal{H}_{Z}$. For $F$, we also
have%
\[
\left\|  F\right\|  _{\mathcal{H}_{Z}}^{2}=\int_{\mathbb{T}}p_{2}\left(
F\right)  \left(  z\right)  \,d\mu\left(  z\right)  =1.
\]
Similarly,%
\begin{align*}
\left\|  M^{n}F\right\|  _{\mathcal{H}_{Z}}^{2}  &  =\int_{\mathbb{T}}%
p_{2}\left(  M^{n}F\right)  \left(  z\right)  \,d\mu\left(  z\right) \\
&  =\int_{\mathbb{T}}R^{n}\left(  p_{2}\left(  F\right)  \right)  \left(
z\right)  \,d\mu\left(  z\right) \\
&  =\int_{\mathbb{T}}R^{n}\left(  \openone\right)  \,d\mu\\
&  =\int_{\mathbb{T}}\openone\,d\mu=1,
\end{align*}
so, for the norm-difference in (\ref{ThmAppr.2(1)}), we have%
\begin{equation}
\left\|  H_{0}-M^{n}\left(  F\right)  \right\|  _{\mathcal{H}_{Z}}%
^{2}=2-2\operatorname{Re}\ip{H_{0}}{M^{n}\left( F\right) }_{\mathcal{H}_{Z}},
\label{eqAppr.2}%
\end{equation}
which means that condition (\ref{ThmAppr.2(1)}) is equivalent to%
\begin{equation}
\lim_{n\rightarrow\infty}\operatorname{Re}\ip{H_{0}}{M^{n}\left( F\right
) }_{\mathcal{H}_{Z}}=1. \label{eqAppr.3}%
\end{equation}
This last term computes out as follows: first introduce the sequence%
\[
D_{n}\left(  z\right)  :=\left|  m_{0}\left(  z\right)  m_{0}\left(
z^{2}\right)  \cdots m_{0}\left(  z^{2^{n-1}}\right)  \right|  ^{2}%
\]
(see (\ref{eqAppr.1}) above), and note that%
\[
R^{\ast\,n}\left(  \openone\right)  =D_{n}\left(  z\right)  .
\]
Since $p_{2}\left(  H_{0}\right)  =\openone$, we know from a theorem of Meyer
and Paiva \cite{MePa93} that $D_{n}\left(  z\right)  \rightarrow0$ for $z$ in
the complement of any neighborhood in $\mathbb{T}$ of $z=1$. But the
$\ip{\,\cdot\,}{\,\cdot\,}_{\mathcal{H}_{Z}}$ term in (\ref{eqAppr.3}) is%
\begin{align*}
\ip{H_{0}}{M^{n}F}_{\mathcal{H}_{Z}}  &  =\ip{M^{n}H_{0}}{M^{n}F}%
_{\mathcal{H}_{Z}}\\
&  =\int_{\mathbb{T}}R^{n}\left(  p\left(  H_{0},F\right)  \right)  \left(
z\right)  \,d\mu\left(  z\right) \\
&  =\int_{\mathbb{T}}R^{\ast\,n}\left(  \openone\right)  \left(  z\right)
\cdot p\left(  H_{0},F\right)  \left(  z\right)  \,d\mu\left(  z\right) \\
&  =\int_{\mathbb{T}}D_{n}\left(  z\right)  p\left(  H_{0},F\right)  \left(
z\right)  \,d\mu\left(  z\right)  .
\end{align*}
If $p\left(  H_{0},F\right)  $ were continuous, the theorem from \cite{MePa93}
would simply give%
\begin{equation}
\int_{\mathbb{T}}D_{n}p\left(  H_{0},F\right)  \,d\mu\underset{n\rightarrow
\infty}{\longrightarrow}p\left(  H_{0},F\right)  \left(  1\right)  ,
\label{eqAppr.4}%
\end{equation}
and the equivalence of (\ref{ThmAppr.2(1)}) and (\ref{ThmAppr.2(2)}) would be
clear from this and (\ref{eqAppr.3}). If (\ref{eqAppr.4}) does not hold, we
still have $p\left(  H_{0},F\right)  \in L^{2}\left(  \mathbb{T}\right)  $ by
Corollary \ref{CorZak.5}. But%
\[
\xi_{n}\left(  z\right)  :=\int_{\mathbb{T}}D_{n}\left(  zw^{-1}\right)
p\left(  H_{0},F\right)  \left(  w\right)  \,d\mu\left(  w\right)
\]
is continuous by the conditions on $m_{0}$, and $\xi_{n}\rightarrow p\left(
H_{0},F\right)  $ in $L^{2}\left(  \mathbb{T}\right)  $ by \cite{MePa93} and
Corollary \ref{CorZak.5}. Suppose now that (\ref{ThmAppr.2(1)}) is given: the
argument from above then shows that (\ref{ThmAppr.2(2)}) holds in the sense of
$Ep\left(  H_{0},F\right)  \left(  1\right)  =1$.
\end{proof}

For the general approximation problem $M^{n}h\rightarrow\varphi$, for
$h,\varphi\in\mathcal{H}$, $M\varphi=\varphi$, we have the following simple result.

\begin{proposition}
\label{ProAppr.3}Consider a system $\left(  M,R,\mathcal{H}\right)  $ as in
Theorem \textup{\ref{ThmR-is.2}} with $M$ a given sub-isometry and
$\mathcal{H}=\mathcal{B}\oplus\mathcal{H}^{\left(  \infty\right)  }$ the
canonical decomposition. Let $h=h_{\mathcal{B}}+h_{\infty}$, $h_{\mathcal{B}%
}\in\mathcal{B}$, $h_{\infty}\in\mathcal{H}^{\left(  \infty\right)  }$ be given.

Then%
\begin{equation}
\left\|  \varphi-M^{n}h\right\|  \geq\left\|  \varphi-M^{n}h_{\infty}\right\|
, \label{eqAppr.5}%
\end{equation}
so that $h_{\infty}$ always gives a better approximation to $\varphi=M\varphi$.
\end{proposition}

\begin{proof}
From Corollary \ref{CorR-is.4}, we have $M^{n}h_{\mathcal{B}}\in\mathcal{B}$
and $M^{n}h_{\infty}\in\mathcal{H}^{\left(  \infty\right)  }$ for all
$n=0,1,2,\dots$, so%
\begin{equation}
\varphi-M^{n}h=\underset{\mathcal{H}^{\left(  \infty\right)  }}{\underbrace
{\varphi-M^{n}h_{\infty}}}+\underset{\mathcal{B}}{\underbrace{M^{n}%
h_{\mathcal{B}}}}, \label{eqAppr.6}%
\end{equation}
and the result is immediate from this.
\end{proof}

When Proposition \ref{ProAppr.3} is applied to the scaling function
$\varphi:=\left.  \chi_{\left[  -\pi,\pi\right]  }^{{}}\right.  \spcheck$ from
Example \ref{ExaSiV.6}, we arrive at the following necessary condition on a
function $h\in L^{2}\left(  \mathbb{R}\right)  $ for $\left\|  M^{n}%
h-\varphi\right\|  _{2}\underset{n\rightarrow\infty}{\longrightarrow}0$ to hold.

\begin{corollary}
\label{CorAppr.4}Let $\varphi:=\left.  \chi_{\left[  -\pi,\pi\right]  }^{{}%
}\right.  \spcheck$, and let $h\in L^{2}\left(  \mathbb{R}\right)  $ be given.
Assume that the Fourier transform $\hat{h}$ is continuous on $\mathbb{R}$.
Then a necessary condition for $\left\|  M^{n}h-\varphi\right\|  _{2}%
\underset{n\rightarrow\infty}{\longrightarrow}0$ is that%
\begin{equation}
\sum_{n\in\mathbb{Z}\setminus\left\{  0\right\}  }\left|  \hat{h}\left(  2\pi
n\right)  \right|  ^{2}=0, \label{eqAppr.7}%
\end{equation}
i.e., $\hat{h}\equiv0$ on $\left\{  2\pi n\mid n\in\mathbb{Z}\setminus\left\{
0\right\}  \right\}  $. Specifically, if the sum in \textup{(\ref{eqAppr.7})}
is $\neq0$, then the cascade algorithm diverges on $h$.
\end{corollary}

\begin{proof}
In the discussion of Example \ref{ExaSiV.6}, we showed that%
\[
\left\|  M^{n}h_{\mathcal{B}}\right\|  _{2}^{2}=\int_{\mathbb{T}}R^{n}\left(
p_{2}\left(  h_{\mathcal{B}}\right)  \right)  \left(  z\right)  \,d\mu\left(
z\right)  ,
\]
and, as $n\rightarrow\infty$, the limit of that term is $p_{2}\left(
h_{\mathcal{B}}\right)  \left(  z=1\right)  $. Let $z=e^{-i\omega}$, $-\pi
\leq\omega<\pi$. Then we saw that
\begin{equation}
p_{2}\left(  h_{\mathcal{B}}\right)  \left(  e^{-i\omega}\right)  =\sum
_{n\in\mathbb{Z}\setminus\left\{  0\right\}  }\left|  \hat{h}\left(
\omega+2\pi n\right)  \right|  ^{2}. \label{eqAppr.8}%
\end{equation}
This follows from the identity%
\[
\widehat{h_{\mathcal{B}}}\left(  \omega\right)  =\chi_{\mathbb{R}%
\setminus\left[  -\pi,\pi\right\rangle }^{{}}\left(  \omega\right)  \hat
{h}\left(  \omega\right)  ,\qquad\omega\in\mathbb{R},
\]
from Example \ref{ExaSiV.6}. Hence%
\[
\lim_{n\rightarrow\infty}\left\|  M^{n}h_{\mathcal{B}}\right\|  _{2}^{2}%
=\sum_{l\in\mathbb{Z}\setminus\left\{  0\right\}  }\left|  \hat{h}\left(  2\pi
l\right)  \right|  ^{2},
\]
and it is clear from (\ref{eqAppr.6}) that $M^{n}h\nrightarrow\varphi$ in
$L^{2}\left(  \mathbb{R}\right)  $ if this term is nonzero.
\end{proof}

\begin{remark}
\label{RemAppr.5}A simple calculation shows that the sum $\sum_{n\neq0}$ in
\textup{(\ref{eqAppr.7})} is nonzero if we take $h=h_{b}=\frac{1}{\sqrt
{b}\mathstrut}\chi_{\left[  0,b\right\rangle }^{{}}$ when $b\in\mathbb{R}_{+}%
$, $b$ \emph{irrational.} \textup{(}The argument can be done with the Poisson
summation formula.\textup{)}

Another example where the sum \textup{(\ref{eqAppr.7})} is nonzero is a
Gaussian, e.g., $h\left(  x\right)  =\pi^{-\frac{1}{2}}e^{-\frac{x^{2}}{2}}$.
The nonzero sum $\sum_{n\in\mathbb{Z}\setminus\left\{  0\right\}  }\left|
\hat{h}\left(  2\pi n\right)  \right|  ^{2}$ can be worked out from Gauss's
lattice sum formula; see, e.g., \cite[p.\ 140]{DyMc72} for details.
\end{remark}

We finally note the following application to the cascade operator $M$ in
$L^{2}\left(  \mathbb{R}\right)  $ associated with some given filter $m_{0}$
satisfying (\ref{IntrAxiom(1)})--(\ref{IntrAxiom(3)}) in the Introduction. For
$f,h\in L^{2}\left(  \mathbb{R}\right)  $, recall the form%
\[
p\left(  f,h\right)  \left(  z\right)  =\sum_{n\in\mathbb{Z}}z^{n}%
\int_{\mathbb{R}}\overline{f\left(  x-n\right)  }h\left(  x\right)  \,dx.
\]
The results in Section \textup{\ref{Poof}} show that%
\[
p\colon L^{2}\left(  \mathbb{R}\right)  \times L^{2}\left(  \mathbb{R}\right)
\longrightarrow L^{1}\left(  \mathbb{T}\right)
\]
is well defined and bounded, i.e.,%
\[
\left\|  p\left(  f,h\right)  \right\|  _{1}\leq\left\|  f\right\|  _{2}%
\cdot\left\|  h\right\|  _{2}%
\]
for all $f,h\in L^{2}\left(  \mathbb{R}\right)  $.

Returning to the approximation problem%
\[
M^{n}h\longrightarrow\varphi,\qquad h=h_{\mathcal{B}}+h_{\infty},\qquad
M\varphi=\varphi,
\]
suppose $z\mapsto p\left(  f,h\right)  \left(  z\right)  $ is continuous when
$f,h$ (as a pair) is picked from $\varphi,h_{\mathcal{B}},h_{\infty}$. We have
the following:

\begin{proposition}
\label{ProAppr.6}If $\varphi$ satisfies $p_{2}\left(  \varphi\right)
=p\left(  \varphi,\varphi\right)  \equiv1$, then
\[
p\left(  \varphi-h_{\infty},h_{\mathcal{B}}\right)  \left(  z=1\right)  =0.
\]
\end{proposition}

\begin{proof}
We have%
\begin{align*}
0  &  =\ip{\varphi-M^{n}h_{\infty}}{M^{n}h_{\mathcal{B}}}_{L^{2}\left(
\mathbb{R}\right)  }\\
&  =\int_{\mathbb{T}}p\left(  \varphi-M^{n}h_{\infty},M^{n}h_{\mathcal{B}%
}\right)  \left(  z\right)  \,d\mu\left(  z\right) \\
&  =\int_{\mathbb{T}}p\left(  M^{n}\varphi,M^{n}h_{\mathcal{B}}\right)
\,d\mu\left(  z\right)  -\int_{\mathbb{T}}p\left(  M^{n}h_{\infty}%
,M^{n}h_{\mathcal{B}}\right)  \,d\mu\left(  z\right) \\
&  =\int_{\mathbb{T}}R^{n}\left(  p\left(  \varphi,h_{\mathcal{B}}\right)
\right)  \left(  z\right)  \,d\mu\left(  z\right)  -\int_{\mathbb{T}}%
R^{n}\left(  p\left(  h_{\infty},h_{\mathcal{B}}\right)  \right)  \left(
z\right)  \,d\mu\left(  z\right) \\
&  \underset{n\rightarrow\infty}{\longrightarrow}p\left(  \varphi
,h_{\mathcal{B}}\right)  \left(  1\right)  -p\left(  h_{\infty},h_{\mathcal{B}%
}\right)  \left(  1\right)  ,
\end{align*}
and the claim follows.
\end{proof}

The last step was based on the idea from the proof of Theorem \ref{ThmAppr.2}
above, and the Meyer--Paiva \cite{MePa93} result. The latter states that, if
$M$ and $\varphi$ are such that $p_{2}\left(  \varphi\right)  \equiv1$, then,
for all continuous functions $\xi$ on $\mathbb{T}$, we have%
\[
\lim_{n\rightarrow\infty}\int_{\mathbb{T}}\left(  R^{n}\xi\right)  \left(
z\right)  \,d\mu\left(  z\right)  =\xi\left(  1\right)  .
\]
An alternative approach to this limit problem is also given in \cite{Kea72}
and \cite{Rue90}.

\begin{acknowledgements}
Helpful conversations with co-authors of earlier papers (Ola Bratteli, David
E. Evans, Steen Pedersen, and Reinhard Werner) are gratefully acknowledged, as
is the excellent typesetting by Brian Treadway.
\end{acknowledgements}

\vspace{\betweenpartsskip}

\baselineskip=12pt
\ifx\undefined\bysame
\newcommand{\bysame}{\leavevmode\hbox to3em{\hrulefill}\,}
\fi

\vspace{\betweenpartsskip}

\baselineskip=14pt
\noindent Department of Mathematics\\
The University of Iowa\\
14 MacLean Hall\\
Iowa City, IA 52242-1419\\
U.S.A.\\[\baselineskip]
E-mail: \texttt{jorgen@math.uiowa.edu}
\vspace{\betweenpartsskip}

\noindent \subjclassname : Primary 46L60, 47D25, 42A16, 43A65; Secondary 46L45, 42A65, 41A15


\end{document}